\RequirePackage{fix-cm}

\documentclass{article}
\usepackage{graphicx}
\usepackage{csquotes}
\usepackage{amsmath,mathtools}
\usepackage{url}
\usepackage{algorithm}
\usepackage{algorithmic}
\usepackage{enumerate}
\usepackage{multirow}
\usepackage{enumitem}
\usepackage{hyperref}
\usepackage{placeins}
\usepackage{amsfonts}
\usepackage{array}
\usepackage{overpic}
\usepackage[round]{natbib}
\usepackage{xspace}

\DeclareMathOperator*{\argmin}{arg\,min}
\newcommand{\trace}{\ensuremath{\mathrm{tr}}}

\renewcommand{\deg}[1]{\ensuremath{\mathrm{deg}\left(#1\right)}}
\newcommand{\pr}[2][]{\ensuremath{\mathrm{Pr}_{#1}\!\left(#2\right)}}
\renewcommand{\vec}[1]{\boldsymbol{#1}}
\newcommand{\cov}[2][]{\ensuremath{\mathrm{Cov}\left(#1,#2\right)}}
\newcommand{\var}[1]{\ensuremath{\mathrm{Var}\left(#1\right)}}

\newcommand{\R}{\mathbb{R}} 
\newcommand{\perm}{\mathfrak{S}} 
\newcommand{\spd}[1]{\ensuremath{\mathcal{S}^+_{#1}}} 
\newcommand{\distrib}[1]{\ensuremath{\mathcal{D}\left(#1\right)}} 
\newcommand{\orthant}{\mathcal{O}}
\newcommand{\preunrooted}[1]{\ensuremath{U_{\!#1}}}
\newcommand{\unrooted}[1]{\ensuremath{\mathcal{U}_{#1}}}
\newcommand{\bhv}[1]{\ensuremath{\mathrm{BHV}_{\!#1}}}
\newcommand{\bhvspace}{BHV tree space\xspace}
\newcommand{\wald}{wald space\xspace}
\newcommand{\w}[1]{\ensuremath{\mathcal{W}_{#1}}}
\newcommand{\prew}[1]{\ensuremath{W_{#1}}}
\newcommand{\bhvmetric}{\ensuremath{d_\mathrm{BHV}}}

\newcommand{\covmetric}{\ensuremath{d_\mathrm{cov}}}
\newcommand{\inducedcovmetric}{\ensuremath{d^*_\mathrm{cov}}}
\newcommand{\jsmetric}{\ensuremath{d_{JS}}}

\newcommand{\ia}{\ensuremath{i}}
\newcommand{\ib}{\ensuremath{j}}
\newcommand{\ic}{\ensuremath{k}}
\newcommand{\id}{\ensuremath{l}}
\newcommand{\ileaf}{\ensuremath{u}}
\newcommand{\jleaf}{\ensuremath{v}}
\newcommand{\ea}{{\ensuremath{e}}}
\newcommand{\eb}{{\ensuremath{\tilde{e}}}}

\newtheorem{Prop}{Proposition}[section]
\newtheorem{lemma}[Prop]{Lemma}

\newtheorem{theorem}[Prop]{Theorem}
\newtheorem{remark}[Prop]{Remark}

\begin{document}

\title{Information geometry for phylogenetic trees
}
%
%
%
%
%

\author{M. K. Garba\footnote{School of Mathematics, Statistics and Physics, Newcastle University, UK, and Department of Mathematical Sciences, Bayero University, Kano, Nigeria,
               \tt{m.k.garba1@ncl.ac.uk}} \and T. M. W. Nye\footnote{School of Mathematics, Statistics and Physics, Newcastle University, UK, \tt{tom.nye@ncl.ac.uk}} \and J. Lueg\footnote{Felix-Bernstein-Institute for Mathematical Statistics in the Biosciences, Georg-August-Universit\"at at G\"ottingen, Germany, \tt{jonas.lueg@stud.uni-goettingen.de}} \and S. F. Huckemann\footnote{Felix-Bernstein-Institute for Mathematical Statistics in the Biosciences, Georg-August-Universit\"at at G\"ottingen, Germany, \tt{huckeman@math.uni-goettingen.de}}
}

\maketitle

\begin{abstract}

We propose a new space of phylogenetic trees which we call \emph{wald space}. 
The motivation is to develop a space suitable for statistical analysis of phylogenies, but with a geometry based on more biologically principled assumptions than existing spaces: in wald space, trees are close if they induce similar distributions on genetic sequence data.  
As a point set, wald space contains the previously developed Billera-Holmes-Vogtmann (BHV) tree space; it also contains disconnected forests, like the edge-product (EP) space but without certain singularities of the EP space. 
We investigate two related geometries on wald space.  
The first is the geometry of the Fisher information metric of character distributions induced by the two-state symmetric Markov substitution process on each tree.  
Infinitesimally, the metric is proportional to the Kullback-Leibler divergence, or equivalently, as we show, any to $f$-divergence. 
The second geometry is obtained analogously but using a related continuous-valued Gaussian process on each tree, and it can be viewed as the trace metric of the affine-invariant metric for covariance matrices. 
We derive a gradient descent algorithm to project from the ambient
space of covariance matrices to wald space. 
For both geometries we derive computational methods to compute geodesics in polynomial time and show numerically that the two information geometries (discrete and continuous) are very similar. 
In particular geodesics are approximated extrinsically. 
Comparison with the BHV geometry shows that our canonical and biologically motivated space is substantially different.

\end{abstract}

\section{Introduction}

Evolutionary relationships between species are represented by phylogenetic trees, in which the leaves represent present-day species, internal vertices represent speciation events, and edge lengths represent the degree of evolutionary divergence between species \citep{Semple2003}. 
Evolutionary relationships are often subject to a high degree of uncertainty, and so it is natural to consider the space of all possible relationships and probability distributions on this space. 
\citet{BHV2001} were the first to construct a space of all phylogenetic trees on a fixed set of leaves. 
This space, known as Billera-Holmes-Vogtmann or \bhvspace, has a very rich geometry: in particular there is a unique geodesic, or shortest length path, between any two points in the space. 
\bhvspace is a so-called $CAT(0)$ space \citep{BHV2001}, meaning it has globally non-positive curvature, and many of its attractive geometric properties follow from this condition. 
A polynomial time algorithm for computing geodesics and their lengths was subsequently developed \citep{Owen2011}. 
A number of statistical methods for analysing samples of phylogenetic trees have been established, which rely fundamentally on the geometry of \bhvspace by transferring conventional multivariate statistical methods into the new geometrical context. 
Algorithms have been developed for computing sample means \citep{Bacak2014,MOP2015}, for constructing confidence regions for the population mean \citep{Willis2019}, and for performing principal component analysis \citep{Nye2011,Feragen2013,Nye2014,Nye2017}. 
An alternative geometry for phylogenetic trees, known as the tropical tree space \citep{Speyer2004,Lin2019}, arises from regarding phylogenetic trees as distance matrices between the species at the leaves. 
Statistical methods such as calculation of sample means \citep{Lin2018} and principal component analysis \citep{Yoshida2019}, have also been developed in tropical tree space. 

In the BHV and tropical tree spaces, trees are regarded primarily as geometric or algebraic objects, without specific consideration to how phylogenetic trees are estimated or interpreted. 
Phylogenetic trees are typically inferred from genetic sequence data via Markov models of sequence evolution over the edges of the tree \citep{Yang2006}, and we are only concerned with such trees. 
Each phylogenetic tree can therefore be regarded as a probability model for genetic sequence data, and a space of all tree-like probability models can be constructed.   
This idea was first considered by \citet{Kim2000}, and then developed more formally in subsequent papers \citep{Moulton2004,Gill2008}. 
The space is known as the \emph{phylogenetic orange space} or \emph{edge-product} space. 
While the space has been studied from the viewpoint of algebraic geometry \citep{Zwiernik2012,engstrom2013toric}, metric geometry on the space has received little attention. 
Recently, methods for approximately computing `probabilistic' metrics on the edge-product space have been developed \citep{Garba2018}. 
These metrics are defined by mapping each tree to its associated distribution on sequence data, and using a metric between these probability distributions. 
Specifically, each tree represents a distribution on characters, where a character is a map from the $N$ leaves of the tree to some alphabet of letters $\Omega$. 
The Hellinger and Jensen-Shannon metrics are defined between distributions on $\Omega^N$ and are pulled back to give metrics between trees. 
Exact calculation of these metrics involves summation over all possible characters, and so when $N$ is large, \citet{Garba2018} use a simulation procedure to estimate the distance between any pair of trees. 
The probabilistic metrics have substantially different properties than the BHV and tropical metrics. 
For example, if all the edge lengths in a pair of given trees are scaled up linearly, then the BHV and tropical distance between the trees both scale in the same way, while in contrast the probabilistic metrics typically tend to zero. 
This is because the letters at the leaves of the trees become independent from one another as the edge lengths increase, due to genetic saturation. 
The distributions on characters represented by the two trees therefore converge to one another as the edge lengths are scaled up, and the distance tends to zero (see Figure 2 in \citet{Garba2018}). 

The metrics studied by \citet{Garba2018} arise from embedding tree space into the larger `ambient' space of all distributions on characters. 
They are obtained from the lengths of `chordal paths' in the ambient space which do not generally lie within the embedded tree space, and are hence called \emph{extrinsic} metrics. 
In contrast, the BHV metric is an \emph{intrinsic} metric, since it is obtained from the lengths of paths lying within tree space. 
For trees sharing a common branching pattern the BHV metric agrees with the corresponding \emph{extrinsic} metric obtained via an embedding into Euclidean space. 
The statistical methods developed on \bhvspace rely heavily on the intrinsic nature of the metric, and this motivated us to seek intrinsic analogs of the probabilistic metrics. 



The aim of this article is to realize intrinsic metrics and their associated geodesics in a new space of forests, the \emph{\wald}\footnote{This space was first discussed by the authors at the Oberwolfach  1804 meeting ``Statistics for Data with Geometric Structure''
 in the \emph{Schwarzwald} (Black Forest)}, that is related to the edge-product tree space (for the subtle, yet essential differences see the discussion in Section \ref{sec:discussion}), when the underlying assumptions are similar to those for the probabilistic metrics. 
We assume that the infinitesimal distance between two trees is measured using the Fisher information matrix. 
We show that this is equivalent to assuming the infinitesimal squared distance is the Kullback-Leibler divergence, or equivalently, any $f$-divergence.  
Our approach uses ideas from information geometry, which is the study of Riemannian differential geometry on spaces of probability distributions. 
The purpose of developing this geometry on this space of forests is with the ultimate aim of obtaining statistical methods analogous to those on other tree spaces. 
The probabilistic metrics and the information geometry have an important advantage over the BHV and tropical geometries: they have by definition a direct biological interpretation in terms of the evolution of genetic sequences. 
In the information geometry, two trees are close when they determine similar distributions of characters, and as a result they would be potentially indistinguishable if inferred from experimental samples of sequence data. 
Conversely, trees are distant in the information geometry when they induce substantially different distributions. 
In contrast, the BHV and tropical metrics are defined more abstractly without reference to evolutionary models or processes.  
Examples of the biological interpretation of the probabilistic metrics were given by \citet{Garba2018}. 

Our approach has two main parts. 
First, we consider geodesics in the information geometry when the model associated with each phylogenetic tree is the two-state symmetric Markov process. 
This is the simplest discrete Markov model of sequence evolution, for which there are two letters in the alphabet, $\Omega=\{0,1\}$. This model is introduced in Section \ref{scn:waldBHV} along with a formal definition of the \wald and a brief review of BHV space.
The thesis of \cite{phdthesis} contains some comparisons of results obtained using the two-state model versus models with the DNA alphabet.
Geodesics in \wald are constructed locally by numerically integrating a certain differential equation determined by the assumptions on the Riemannian metric. 
We explore geodesics on the space of unrooted trees with $5$ leaves, for which visualization is relatively straightforward, and compare the results with those for \bhvspace. 
This forms Section~\ref{sec:twostateinfogeom} of the paper. 
Secondly, in order to improve computational tractability, we consider an alternative continuous-valued model of evolution on each tree. 
This consists of a Gaussian process which approximates the two-state Markov process by matching its moments. 
The continuous random variables at the leaves of the tree have a multivariate normal distribution with zero mean, for which the covariance matrix is related to the matrix of path lengths between the leaves. 
Numerically solving the differential equations for geodesics is much faster under this set of assumptions, and the geodesics closely resemble those for the two-state model. 
However, solutions are still restricted to trees sharing a common branching pattern, or topology. 
The definition of the Gaussian process on trees and numerical solution of geodesics in the corresponding information geometry are described in Section~\ref{sec:gaussianprocess}.   
The information geometry of multivariate normal distributions with zero mean corresponds to a certain geometry on the space of symmetric positive definite matrices, known as the Fisher-Rao or affine-invariant geometry, and the map from the \wald to covariance matrices is an isometric embedding in this space. 
The geometry on the space of symmetric positive definite matrices is analytically tractable, and geodesics can be computed in polynomial time. 
The embedding therefore gives  
intrinsic and extrinsic metrics on the \wald. 
We describe a projection algorithm from the space of symmetric positive definite matrices into the embedded \wald. 
We then use this algorithm to project geodesics in the ambient space down into \wald in various ways to obtain approximate geodesics between trees with different topologies. 
The embedding in the space of symmetric positive definite matrices and the associated geometry is described in Section~\ref{sec:spdembedding}.
We conclude in Section \ref{sec:discussion} with a detailed discussion of the promises and challenges of our new \wald. 

\section{Background and the new \wald}\label{scn:waldBHV}

\subsection{Phylogenetic trees}\label{sec:trees}

For $N=2,3,\ldots$ we define $\preunrooted{N}$ to be the set of unrooted phylogenetic trees on $N$ taxa.  
More specifically, a tree $T$ is an element of $\preunrooted{N}$ if it satisfies the following conditions. 
First, $T$ contains exactly $N$ vertices with degree 1, which are called \emph{leaves}, and these are bijectively labelled $1,\ldots,N$. 
Secondly, $T$ must contain no vertices with degree $2$. 
Thirdly, each edge $e$ in $T$ is assigned a \emph{length} $\ell^e\geq0$ with $\ell^e\neq0$ if $e$ contains a leaf.
An edge in a tree is called a \emph{pendant edge} if it contains a leaf; otherwise it is called an \emph{internal edge}. 
Similarly, the vertices which are not leaves are called \emph{internal vertices}. 

The edge lengths $\ell^e$ on any given tree $T\in \preunrooted{N}$ define a path length distance between any pair of leaves. 
The path length on $T$ between $\ileaf,\jleaf\in\{1,\ldots,N\}$ will be denoted $\ell_{\ileaf\jleaf}$.

Each tree $T\in \preunrooted{N}$ contains at most $2N-3$ edges, in which case the tree is called \emph{fully resolved} or \emph{bifurcating}, and all internal vertices have degree $3$. 
Trees with fewer edges are called \emph{unresolved}, and for $N>3$, these contain at least one vertex with degree $4$ or more. 
Trees which contain only the $N$ pendant edges joined at a single degree-$N$ internal vertex are called \emph{star trees}. 

A tree $T$ is \emph{rooted} when some internal point $\rho\in T$ is labelled as being the root. 
This is conveniently achieved by adding an additional taxon labelled $0$ to the tree via a pendant edge of length zero. 
It follows that the set of rooted phylogenetic trees satisfies the same conditions as $\preunrooted{N}$, except the leaves are bijectively labelled $0,1,\ldots,N$, and the pendant edge containing taxon $0$ has zero length. 
We will work with unrooted trees, but our results are easily transferred to the space of rooted trees via this relationship. 
 
Every fully resolved tree will correspond to a fully resolved \emph{BHV-tree} (reviewed in Section \ref{scn:BHV}) and to a fully resolved \emph{wald}, as introduced below in Section \ref{sec:waldspacedefinition}. 
In both \bhvspace~and in \wald, unresolved trees will be identified with other trees with certain internal edges having zero length, so that conceptually a missing edge is the same as a zero length edge.

\subsection{Billera-Holmes-Vogtmann tree space}\label{scn:BHV}

\begin{figure}
\begin{center}
\begin{overpic}[width=0.7\textwidth]{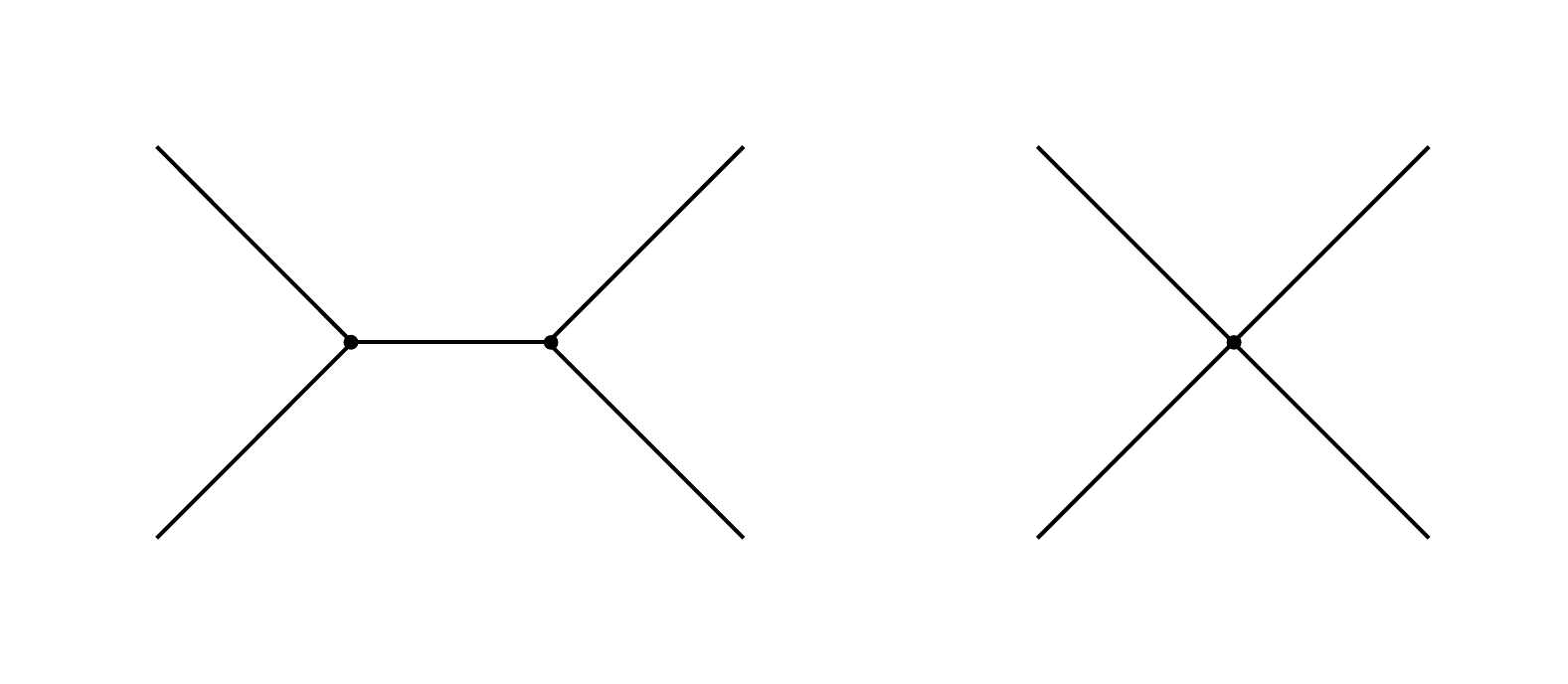}
\put (6,6) {$A$}
\put (6,34) {$B$}
\put (49,6) {$C$}
\put (49,34) {$D$}
\put (62,6) {$A$}
\put (62,34) {$B$}
\put (93,6) {$C$}
\put (93,34) {$D$}
\put (54,20){\scalebox{2.0}{$\sim$}}
\put (93,34) {$D$}
\put (28,23){$0$}
\put (16,11) {$\ell^A$}
\put (17,30) {$\ell^B$}
\put (38,11) {$\ell^C$}
\put (36,30) {$\ell^D$}
\put (71,10) {$\ell^A$}
\put (73,31) {$\ell^B$}
\put (82,10) {$\ell^C$}
\put (81,31) {$\ell^D$}

\end{overpic}
\begin{caption}{\label{fig:BHVequiv}
Two trees in $\preunrooted{N}$ are equivalent under the relation $\sim$ when they are identical after internal edges with length zero are removed, and the vertices at the end of every such edge are merged.
$A,B,C,D$ represent different subtrees joined by edges of length $\ell^A,\ell^B,\ell^C,\ell^D$ to an internal edge with length $\ell=0$ on the left. 
The Markov process $X(t)$ cannot change state on any edge with length zero, so the distribution on $X_1,\ldots,X_N$ is unchanged by removing such edges in this way.}
\end{caption}
\end{center}
\end{figure}

\citet{BHV2001} defined a space of phylogenetic trees, subsequently known as \bhvspace, and described its geometry. 
\bhvspace can be described via an embedding in $\R^d$ for dimension $d$ which increases exponentially with the number of leaves. 
%
%
%
However, we have chosen to describe \bhvspace in a way different from the original authors, and we define it as a quotient space.  
As a result, the \emph{\wald} introduced in the next section is a superset 
of \bhvspace when the spaces are regarded simply as sets, clarifying 
the relationship between the two spaces. 
Importantly, we allow internal edges on trees to have length zero, and under the quotient these are equivalent to trees with those edges missing. 
A second difference is that while \citet{BHV2001} worked with rooted trees, we work with unrooted trees. 
As described in Section~\ref{sec:MarkovProc}, the distribution on binary characters determined by a tree does not depend on the root position under the two-state symmetric model, and so unrooted trees are more natural to work with. 

\bhvspace is defined using the notion of splits, where a \emph{split} is a bipartition of the leaf labels $1,\ldots,N$ into two disjoint sets. 
Cutting an edge of a tree induces such a bipartition of the leaves, and so each edge on a tree corresponds to a split, and the terms \emph{split} and \emph{edge} can be used interchangeably. 
The set of splits represented by a tree is called its \emph{topology}. 

Arbitrary sets of splits do not typically determine valid tree topologies: the splits of a tree must satisfy a compatibility condition. 
For example, the splits $\{1,2\},\{3,4,\ldots,N\}$ and $\{1,3\},\{2,4,\ldots,N\}$ are incompatible, since leaf $1$ cannot be grouped next to both $2$ and $3$ on the same tree.
For any topology $\tau$ with $k$ internal edges, $0\leq k\leq N-3$, the set of trees in $\preunrooted{N}$ with that topology is bijectively parametrized by $\R_{> 0}^N\times \orthant_\tau$ where the first term in the product parametrizes the pendant edge lengths that, by definition, are strictly positive, and $\orthant_\tau = \R^{k}_{\geq 0}$ parametrizes the internal edge lengths. 



The set $\orthant_\tau$ is called the \emph{orthant} associated with topology $\tau$, and  we identify the set of all trees with topology $\tau$ with $\R_{> 0}^N\times \orthant_\tau$. 
Under this identification, the set of all trees $\preunrooted{N}$, as defined in Section~\ref{sec:trees}, is the disjoint union 
\begin{equation*}
\preunrooted{N} = \R_{> 0}^N\times \bigsqcup_\tau \orthant_\tau
\end{equation*}
where the disjoint union is taken over all possible topologies $\tau$. 

The unrooted \bhvspace $\unrooted{N}$ is obtained by taking the quotient of $\preunrooted{N}$ with respect to an equivalence relation: 
\begin{equation*}
\unrooted{N} =  \preunrooted{N} \Big{/} \ \sim. 
\end{equation*}
Two trees in $\preunrooted{N}$ are equivalent under $\sim$ if and only if they are identical modulo the presence of internal splits with zero length, as shown in Figure~\ref{fig:BHVequiv}.
The quotient space factorizes as
\begin{equation*}
\unrooted{N} = \R_{> 0}^N\times \bhv{N}
\end{equation*}
where the first term parametrizes the lengths of the pendant edges and the space $\bhv{N}$ parametrizes the topology and internal edge lengths of the BHV-trees.
When $\tau$ is fully resolved, $\orthant_\tau$ is called a \emph{maximal orthant}. 
Unresolved trees correspond to points on the boundaries of maximal orthants; they can be obtained from fully resolved trees by shrinking internal edge lengths down to zero.

\begin{figure} 
\begin{center}
\begin{overpic}[width=0.65\textwidth]{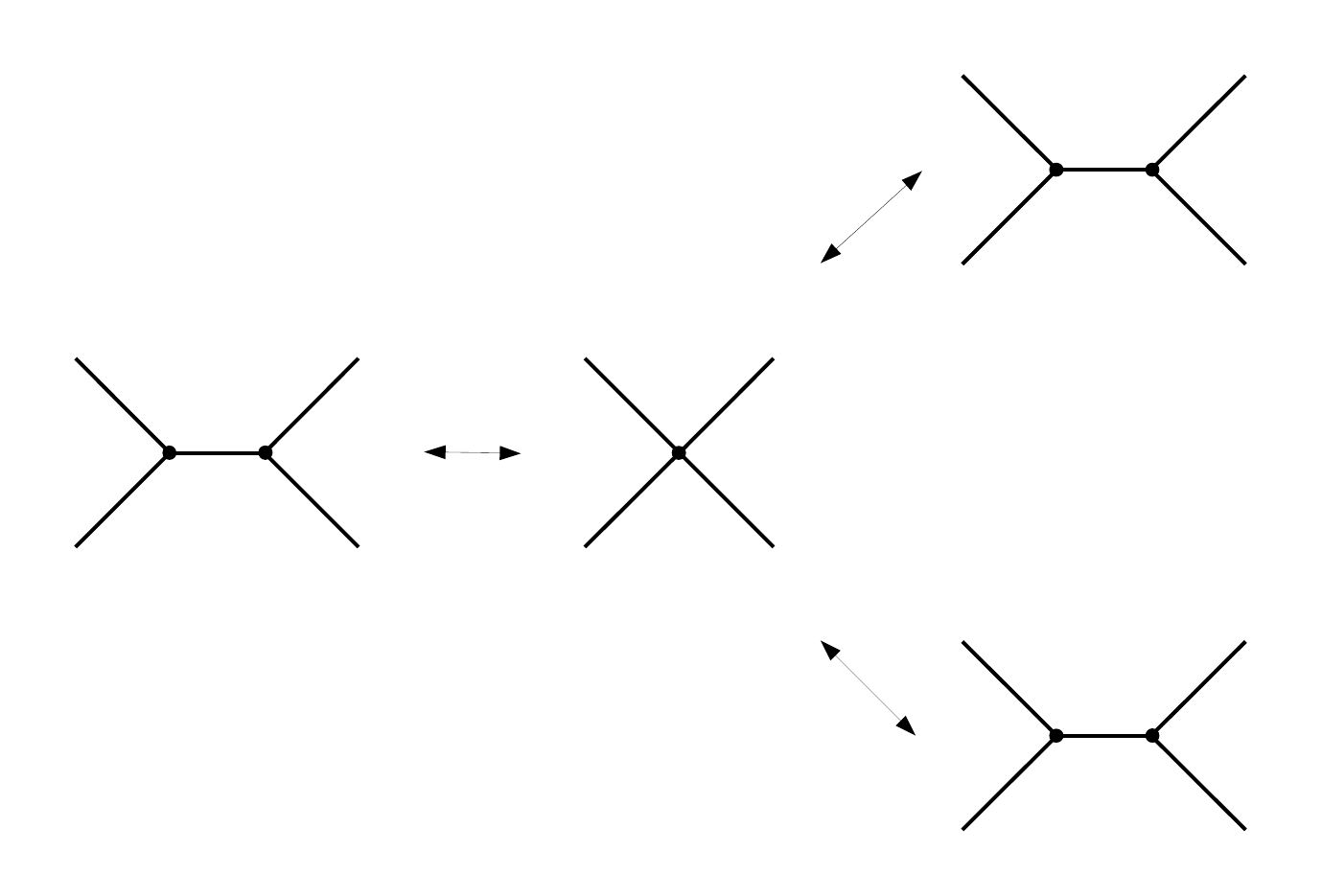}
\put (3,22) {$A$}
\put (3,42) {$B$}
\put (26,22) {$C$}
\put (26,42) {$D$}

\put (69,1) {$A$}
\put (69,21) {$C$}
\put (94,1) {$B$}
\put (94,21) {$D$}

\put (69,44) {$A$}
\put (69,63) {$D$}
\put (94,44) {$B$}
\put (94,63) {$C$}

\put (41,22) {$A$}
\put (41,42) {$B$}
\put (59,22) {$C$}
\put (59,42) {$D$}

\end{overpic}
\end{center}
\begin{caption}{\label{fig:NNI}
When an internal edge from a fully resolved topology is contracted down to length zero (left to centre), there are two fully resolved topologies which can be obtained by expanding out an alternative edge (right). 
$A,B,C,D$ represent subtrees. 
The operation of contracting an internal edge and expanding out an alternative edge is called \textit{nearest neighbour interchange}. 
It follows that at each codimension-$1$ boundary, three maximal orthants are glued together. }
\end{caption}
\end{figure}


Since there are $(2N-5)!!$ fully resolved unrooted topologies, $\bhv{N}$ can be thought of as being constructed by gluing this number of maximal orthants together along their boundaries, where two points are identified if they correspond to the same tree. 
For example, when $N=4$, there are three fully resolved topologies, each of which contains a single internal edge. 
The space $\bhv{4}$ therefore consists of three copies of $\R_{\geq 0}$ glued together at the origin.
The origin corresponds to the star trees, while the location along each of the three copies of $\R_{\geq 0}$ gives the length of the internal edge in each of the three possible fully resolved topologies. 
For $N=5$ there are $15$ possible unrooted tree topologies, each of which contains two internal edges. 
It follows that $\bhv{5}$ consists of $15$ copies of $\R^2_{\geq 0}$ glued along their boundaries.
At each codimension-1 boundary, three maximal orthants are joined together. 
This is because when a single internal edge is contracted to length zero, a degree 4 vertex is obtained, and there are 3 possible ways to add in an edge, including the original edge, in order to obtain a fully resolved topology, as illustrated by Figure~\ref{fig:NNI}. 

The metric on $\bhv{N}$ is constructed as follows. 
The basic idea is that for trees with the same fully-resolved topology but different vectors of internal edge lengths, say $\vec{\ell}_1$ and $\vec{\ell}_2$, the distance is the Euclidean distance $\|\vec{\ell}_1-\vec{\ell}_2\|$, and the corresponding geodesic is the straight line segment in the orthant containing the trees.
\citet{BHV2001} showed that there exists a unique shortest path between any two points in $\bhv{N}$, for which path length is measured using the Euclidean distance in each orthant, and the length of these defines a metric on $\bhv{N}$ which we denote $\bhvmetric$. 
A metric on $\unrooted{N}$, denoted $d_{\unrooted{N}}$, is obtained as the product metric when the metric on pendant edges is taken to be the Euclidean distance.
An algorithm has been developed which constructs geodesics and calculates their lengths in $O(N^4)$ time \citep{Owen2011}.

\subsection{The two-state symmetric Markov model}\label{sec:MarkovProc}

Genetic sequence evolution is typically modelled using discrete-valued continuous-time Markov processes defined over the edges of a tree $T$ \citep{Yang2006,bryant05}. 
DNA sequence evolution is modelled by associating to each point $t\in T$, a random variable $X(t)$ which takes values in an alphabet $\{A,C,G,T\}$. 
In this paper, however, we will consider the two-state symmetric Markov process with alphabet $\Omega=\{0,1\}$.
This simplification is made in order to make the mathematics more tractable and for computational speed. 
Nonetheless, some of the calculations using the two-state symmetric can readily be performed using DNA models.
More details are given in the thesis of \cite{phdthesis} in which simulations show similarity of geometries obtained from the two- and the four-state process.
The transition probability of the symmetric two-state model is defined in terms of the path length $\ell_{t_1t_2}$ between any two points $t_1,t_2\in T$:
\begin{align}
\pr{X(t_2)=X(t_1)} &= \frac{1}{2}\left( 1+e^{-\ell_{t_1t_2}} \right),\quad\text{and} \notag\\
\pr{X(t_2)\neq X(t_1)} &= \frac{1}{2}\left( 1-e^{-\ell_{t_1t_2}} \right). 
\label{equ:twostatetrans}
\end{align}
The stationary distribution of this Markov process is $Bern(1/2)$, and the process is assumed to be in its stationary state over the tree.
As a result, for all $t\in T$, $X(t)$ has a marginal Bernoulli distribution, $X(t)\sim Bern(1/2)$. 
While the random variables $X_1,\ldots,X_N$ at the leaves of the tree have the same marginal distributions, they are not independent since the tree imposes a dependence structure. 
The following lemma determines certain moments of the process $X(t)$ giving insight on the dependence structure of $X_1,\dots,X_N$. The proof is straightforward using the transition probabilities in Equation~$\eqref{equ:twostatetrans}$.

\begin{lemma}\label{lem:twostatemoments}
\begin{enumerate}
\item If $X_1,\ldots,X_N$ are the random variables at the leaves of a tree $T\in \preunrooted{N}$ determined by the discrete Markov process defined above, then $\cov[X_\ileaf]{X_\jleaf}=\frac{1}{4}\exp(-\ell_{\ileaf\jleaf})$ where $\ell_{\ileaf\jleaf}$ is the path length between leaves $\ileaf$ and $\jleaf$. 
\item If $t_1,t_2\in T$ are path length $\ell_{t_1t_2}$ apart, then the conditional distribution of $X(t_2)$ given $X(t_1)=\omega\in\{0,1\}$ has variance $\frac{1}{4}\big(1-\exp(-2\ell_{t_1t_2})\big)$.
\end{enumerate}
\end{lemma}



It is straightforward to simulate realizations of $X(t)$ in the following way. 
First simulate a Poisson process with rate 1 independently on each edge of the tree. 
The positions of the simulated events correspond to points at which $X(t)$ changes parity. 
Secondly, pick any point $t_0\in T$ which is not a change point and sample $X(t_0)$ from $Bern(1/2)$. 
The change points generated from the Poisson process then determine the value of $X(t)$ for all other $t\in T$. 
The distribution obtained is independent of the choice of $t_0$, because the Markov process is reversible. In particular, the Markov process is independent of the choice of $t_0$, which could be considered as a root. 

Under the model, each edge length can be interpreted as the expected number of change points that occur over the edge. 
Internal edges are allowed to have length zero, which means that no change in $X(t)$ occurs over the edge. 
On the other hand, when edges are long, the number of changes is likely to be large, and the letters at either end of the edge are weakly correlated. 
Biologists refer to this effect as \emph{saturation}. 
A fixed change of edge length $\delta\ell$ therefore has more effect on the distribution of characters when applied to a short edge as opposed to a long edge in some given tree. 
For example, an increase of $\delta\ell=0.1$ to an edge with length $\ell=0.1$ approximately doubles the probability that the letters at either end of the edge are different, but the same change to an edge of length $\ell=10$ has almost no effect on this probability, which due to saturation is very close to $1/2$. 
This idea becomes important when we consider defining distances between trees via the information they represent, 
in particular using the probability mass function of the nontrivial distribution of $(X_1,\ldots,X_N)$.

\begin{remark}\label{rmk:polynomial}
1. The probability mass function of $(X_1,\ldots,X_N)$ determined by $T$ is denoted $p_T(s)$ where $s\in\{0,1\}^N$ is called a \emph{binary character}. 
Given any binary character $s$, the values of $p_T(s)$ can be evaluated via a recursive algorithm \citep{bryant05}, described in  Appendix A. 
Appendix A also contains a modified form of the algorithm which is used to compute exactly the derivatives of $p_T(s)$ with respect to the edge lengths. 

2. Also, as shown in Appendix A,  $p_T(s)$ is a multivariate polynomial in $1+e^{\ell^k}$ and $1-e^{\ell^k}$ where $\ell^k$ ranges over the edge lengths in $T$.

3. Crucially, the map $T \mapsto p_T$ from fully resolved trees in $\preunrooted{N}$ to the space of probability mass functions on $\{0,1\}^N$ is injective up to the equivalence relation introduced in Section \ref{scn:BHV} (\citep{rogers1997consistency,Allman2008}). 
This has two implications: first that the probability mass function $p_T$ uniquely characterizes each element of $\unrooted{N}$, and secondly that metrics on distributions of characters pull back to define a metric between fully resolved trees, as described in Section \ref{sec:probmetrics}. 
\end{remark}

\subsection{A new forest space: the \wald}
\label{sec:waldspacedefinition}

The following \wald gives an alternative viewpoint of phylogenetic trees by regarding them as Markov models for sequence evolution \citep{Kim2000,Moulton2004,Gill2008}. 
We will give a description of the \wald that is related to the \emph{edge-product} space of previous authors, by defining it as a quotient space which adds trees with infinitely long edges to \bhvspace. 
As described in the introduction, each probabilistic metric considered by \citet{Garba2018} converges to zero in the limit as all the edge lengths in a pair of trees simultaneously tend to infinity. 
This behaviour indicates that shortest paths might cross a tree with infinitely long edges, which is why we add this \emph{boundary at infinity}. 
Furthermore, we allow pendant edges of length zero under certain conditions.
Thus, in the \wald, edge lengths $\ell^e$ take values in $\R_{\geq 0}\cup\{\infty\}$.
It is convenient to reparametrize to the \emph{$\vec{\lambda}$-parametrization} by defining \emph{weights} $\lambda^e = 1-\exp(-\ell^e)$, so $\lambda^e\in[0,1]$. 
Under this transformation, $\bhv{N}$ becomes a set of unit cubes, rather than orthants, glued along faces where $\lambda^e=0$ for one or more edges.
The \wald is defined by imposing additional gluing rules on faces where $\lambda^e=1$. 

In order to be able to identify trees with infinitely long edges along faces where $\lambda^e=1$, we construct the \wald from forests, that is, disjoint unions of unrooted trees. 
We start with a preliminary space leading to the definition of the \wald further below. Let $\prew{N}$  be the collection of forests satisfying the following necessary and sufficient conditions for each $F\in \prew{N}$.
\begin{enumerate}
\item The forest $F$ contains exactly $N$ labelled vertices, these are called leaves and 
labelled $1,\ldots,N$. 
\item There are no unlabelled vertices in $F$ of degree 0,1 and 2.
\item For every pair of leaves $\ileaf, \jleaf$ in the same tree in $F$, at least one edge $e$ on the unique path from $\ileaf$ to $\jleaf$ satisfies $\lambda^e > 0$.
\end{enumerate}

Clearly, $\preunrooted{N}\subset \prew{N}$. 
The condition on the edge weights ensures that no pair of leaves is coincident and consequently that metrics are always well-defined, as described in Section~\ref{sec:probmetrics}.
We impose an equivalence relation $\sim$ on $\prew{N}$, defined by the following two rules. 

\begin{figure} 
\begin{center}
\begin{overpic}[width=0.8\textwidth]{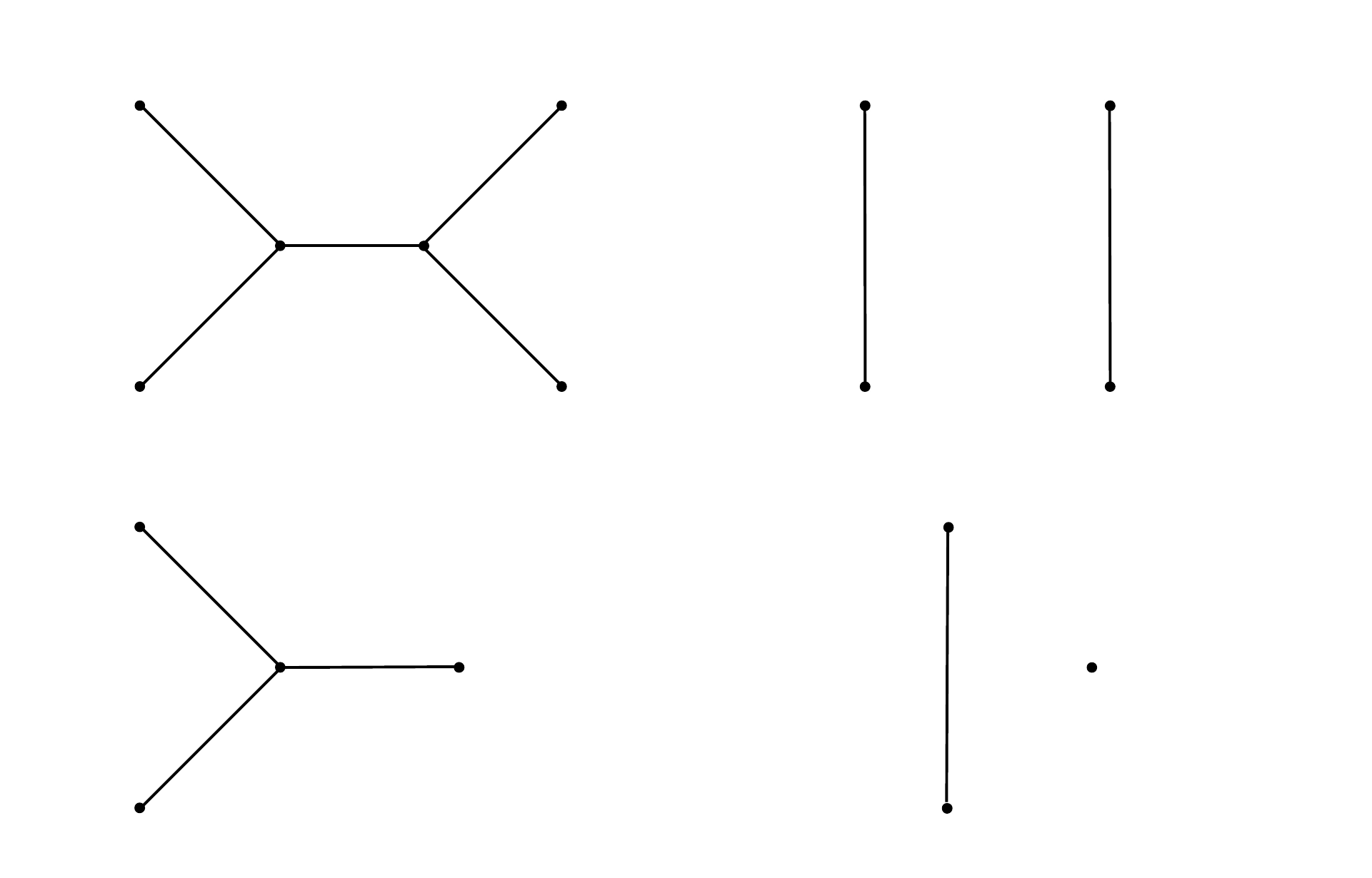}

\put (2,45) {(a)}

\put (7,33) {$A$}
\put (7,57) {$B$}
\put (42,33) {$C$}
\put (42,57) {$D$}

\put (15,37) {$\lambda^A$}
\put (15,54) {$\lambda^B$}
\put (33,37) {$\lambda^C$}
\put (33,54) {$\lambda^D$}
\put (25,48) {$1$}

\put (43,43) {$\bigcup ~ F$}
\put (53,45){\scalebox{2.0}{$\sim$}}

\put (60,33) {$A$}
\put (60,57) {$B$}
\put (82,33) {$C$}
\put (82,57) {$D$}

\put (64,51) {$\lambda'^{AB}$}
\put (83,51) {$\lambda'^{CD}$}

\put (88,43) {$\bigcup F$}

\put (2,15) {(b)}

\put (7,3) {$A$}
\put (7,26) {$B$}
\put (35,14) {$L$}

\put (15,7) {$\lambda^A$}
\put (15,24) {$\lambda^B$}
\put (26,17) {$1$}

\put (43,12) {$\bigcup ~ F$}
\put (53,14){\scalebox{2.0}{$\sim$}}

\put (66,3) {$A$}
\put (66,26) {$B$}
\put (78,11) {$L$}

\put (62,20) {$\lambda'^{AB}$}

\put (88,12) {$\bigcup ~ F$}

\end{overpic}
\end{center}
\begin{caption}{\label{fig:edgeprodequiv}
Illustration of \emph{boundary at infinity} rule used to define $\w{N}$.
In both (a) and (b), the forests on the left are equivalent to the forests on the right. 
$A,B,C,D$ are subtrees. 
(a) Internal edge with weight $1$: the edge is deleted, disconnecting the tree. 
The resulting edges between subtrees $A$ and $B$ are replaced by a single edge with weight $\lambda'^{AB}=\lambda^A+\lambda^B-\lambda^A\lambda^B$ and similarly for $C,D$. 
(b) Pendant edge, where $L$ is a leaf. The pendant edge with weight $1$ is removed, and the resulting edges between $A$ and $B$ are replaced by a single edge with weight $\lambda'^{AB}$.  
The term $F$ in both panels refers to other disconnected components in the forests. }
\end{caption}

\end{figure}

\noindent\textbf{BHV boundary rule:} 
Given $F_1,F_2\in \prew{N}$, suppose all internal edges with $\lambda^e=0$ are removed from the forests, and the vertices at either end of each such edge are merged. If the resulting forests are identical, then $F_1\sim F_2$. 
The rule is the same as that in Figure~\ref{fig:BHVequiv}. 

\noindent\textbf{Boundary at infinity rule:} 
Suppose $F\in \prew{N}$ contains an edge with $\lambda^e=1$, and that $F$ is modified as follows. 
The edge with $\lambda^e=1$ is removed, disconnecting the tree it belongs to. 
If this results in any unlabelled vertex having degree 2, then those vertices are removed. 
If $v$ is such a vertex, and the two adjacent edges $e,\tilde{e}$ have weights $\lambda^e,\lambda^{\tilde{e}}\in[0,1]$, then $e,\tilde{e}$ are replaced by a single edge with weight $\lambda^\ea+\lambda^\eb-\lambda^\ea\lambda^\eb$, as is further explained below. 
Now suppose $F_1,F_2\in \prew{N}$, and this process of modifying unit-weight edges is applied to both forests. 
Then $F_1\sim F_2$ if the resultant forests are identical, as illustrated in Figure~\ref{fig:edgeprodequiv}. 

The \emph{\wald} $\w{N}$ is defined to be the quotient $\prew{N}\big{/}\sim$ and it immediately follows that as sets $\unrooted{N}\subset\w{N}$, but the geometry imposed on $\w{N}$ will be completely different from the geometry of the BHV space.

The boundary rule at infinity requires some explanation. 
The rule declares that edges of weight $\lambda^e=1$ (or equivalently length $\ell^e=\infty$) can be deleted from a forest $F$, but unlike the BHV rule for which the vertices at the ends of the edge are merged, edge removal disconnects a tree in $F$. 
When resulting degree-2 vertices are removed, the edge length is preserved so that the new edge has length $\ell^\ea+\ell^\eb$. 
The corresponding weight $\lambda$ is given by $\lambda=1-\exp(-(\ell^\ea+\ell^\eb))=\lambda^\ea+\lambda^\eb-\lambda^\ea\lambda^\eb$. 
Unlike the BHV boundary rule, in which finitely many trees are identified in each equivalence class, infinitely many combinations of edge weights $\lambda^\ea,\lambda^\eb$ give rise to the same value $\lambda^\ea+\lambda^\eb-\lambda^\ea\lambda^\eb$. 
It follows that an uncountable collection of forests can be identified into a single equivalence class in $\w{N}$. 

In the edge-product space \citep{Moulton2004,Gill2008}, an alternative parametrization is used, defining $\mu^e=1-\lambda^e=\exp(-\ell^e)$ to be the weight of edge $e$. 
This parametrization has the advantage that sums of edge lengths $\ell^1+\cdots+\ell^m$ become products of edge weights $\mu^{1}\times\cdots\times\mu^m$ (hence the name `edge-product'). 
The boundary at infinity rule is simpler under this parametrization: the weights in Figure~\ref{fig:edgeprodequiv} panel (b) become $\mu^\ea$, $\mu^\eb$ and $0$ on the left and $\mu^\ea\mu^\eb$ on the right. 
However, under the $\vec{\mu}$-parametrization, the BHV boundary with $\ell^e= 0$ lies on faces of cubes with $\mu^e=1$, whereas the boundary at infinity has $\mu^e=0$.
We prefer to work with the $\vec{\lambda}$-parametrization since it gives a more intuitive interpretation of the weights, i.e. $\ell^e = 0$ corresponds to $\lambda^e = 0$ and $\ell^e = \infty$ corresponds to $\lambda^e = 1$. 
Forests which contain more than one connected component lie in the faces of cubes with at least one $\lambda^e = 1$.
Since the pendant edges can be expanded out to infinite length, they are also subject to the boundary at infinity rule, and so the representation of pendant edge lengths in $\w{N}$ is not via a product geometry, as it is for \bhvspace. 
While the star trees correspond to all internal edges having zero length, $\w{N}$ also contains a point which consists of $N$ isolated vertices. 

The BHV boundaries enable tree topologies to be changed via nearest neighbour interchange (NNI) operations (as illustrated by Figure~\ref{fig:NNI}). 
The boundary at infinity corresponds to a different topological operation, called tree bisection and reconnection (TBR) \citep{allen2001}. 
Under this operation, an edge $e$ in a tree can be expanded up to the boundary $\lambda^e=1$. 
Removing the edge bisects the tree, and the two components can be reconnected by an edge $\eb$ with $\lambda^\eb=1$ placed arbitrarily between the two trees. 
Reducing the weight $\lambda^\eb$ down from 1 then gives a tree with a topology different from the original tree. 
It follows that there exist continuous paths in the \wald between trees with different topologies, which pass through the boundary at infinity and, as a result, change tree topology via TBR operations. 
This is in contrast to \bhvspace in which paths between trees of different topologies involve only NNI operations, as edges are contracted down to length zero and alternative edges are expanded out.

While the set $\w{N}$ was defined above via an equivalence relation on forests, we also need to understand how it parametrizes Markov models and then characterise its elements again as probability mass functions on $\{0,1\}^N$.
The two-state symmetric Markov process extends from being defined on trees to forests by taking the process on each connected component in a forest to be independent of the other components. 
This defines a distribution $p_F$ on $\{0,1\}^N$ for each $F\in \prew{N}$. 
In fact the distribution uniquely determines the equivalence class of $F$, and vice versa, as the following lemma shows. 

\begin{lemma}\label{lem:edgeprod}
Given $F_1,F_2\in \prew{N}$, then $F_1\sim F_2$ if and only if $p_{F_1}(s)=p_{F_2}(s)$ for all $s$. 
\end{lemma}

A proof is given in the Appendix.

Note that the forest consisting of $N$ isolated vertices corresponds to the random variables $X_1,\ldots,X_N$ being independent, and this can be obtained from any tree by 
expanding all edges one after the other, as, by definition, there is at least one edge between any two leaves. 

\section{Information geometry for the two-state symmetric model}\label{sec:twostateinfogeom}

Information geometry provides methods for constructing metrics and geodesics on parametrized sets of probability distributions. 
In this section we embed \wald in the space of distributions of two-state characters, and investigate the corresponding information geometry analytically and computationally. 

\subsection{Geometry of embeddings}\label{sec:embeddings}

Suppose that $\theta:X\rightarrow Y$ where $(Y,d)$ is a metric space and $\theta$ is injective. 
We will say that $X$ is embedded in $Y$, and refer to $Y$ as the \emph{ambient} space. 
The embedding can be used to construct certain metrics on $X$. 
First, since $\theta$ is injective, $d$ pulls back to define a metric on $X$ which we denote $d_X$:
\begin{equation*}
d_X(x_1,x_2) = d(\theta(x_1),\theta(x_2))
\end{equation*}
for all $x_1,x_2\in X$. 
The pull back metric is often called the induced \emph{extrinsic metric} and, it is simply the restriction of $d$ to $X\subseteq Y$, and so when the context is clear, it is also denoted $d$. 
The probabilistic metrics described in Section~\ref{sec:probmetrics} are constructed in this way. 
A second metric, called the induced \emph{intrinsic metric} and denoted $d^*(x_1,x_2)$, is defined as the infimum of the length of all possible paths in $X$ between $x_1,x_2\in X\subseteq Y$ when path length is measured using the metric $d$. 
If no path with finite length exists between $x_1$ and $x_2$ then $d^*(x_1,x_2)=\infty$, in which case $d^*$ is not a metric. 
Details of this construction of the induced intrinsic metric are given by \citet{BH1999}. 
The metrics $d$ and $d^*$, if well-defined, give
$X$ the structure of a \emph{length space}, which is a space in which the metric between points $x_1,x_2$ is the infimum of the lengths of paths between those points. 
Length spaces are similar to geodesic metric spaces, except that the infimum is not necessarily achieved by a path lying within the space; in a geodesic metric space a minimum length path exists between every pair of points, and so every geodesic metric space is a length space. 
An example of a length space which is not a geodesic metric space is $\R^2$ with the origin removed and the Euclidean metric. 
Points antipodal to the origin cannot be joined by a geodesic, but the distance between them is the infimum of the lengths of paths joining the points. 

In order to illustrate the relationship between $d$ and $d^*$, consider the example of the embedding of the unit sphere $X=S^2$ in $Y=\R^3$ equipped with the Euclidean metric $d$. 
For any two points $x_1,x_2$, $d(x_1,x_2)$ is the length of the straight line segment in the ambient space $\R^3$ joining the points. 
This metric is usually called the \emph{chordal metric} on $S^2$. 
However, when we consider paths between $x_1,x_2$ which are restricted to lie in $S^2$, the shortest paths (with respect to $d$) are great circles, and the induced metric $d^*$ is the arc length metric. 
In fact, $S^2$ is a geodesic metric space, since the infimum of path length is always achieved by a great circle. 

In the following Section \ref{sec:probmetrics}, the \wald $\w{N}$ will be embedded in the space of distributions of characters. Later in Section \ref{sec:gaussianprocess} it will be embedded in  the space of $N\times N$ symmetric positive definite matrices. 
Each embedding will be used to construct metrics on $\w{N}$.   

\subsection{Probabilistic metrics}\label{sec:probmetrics}

Here we briefly describe the probabilistic metrics developed by \citet{Garba2018} since these will be used for comparison with other metrics. 
The Kullback-Leibler divergence is a commonly used measure of the difference between two distributions. 
Given two probability mass functions $p,q$ on characters $\{0,1\}^N$, the Kullback-Leibler divergence from $q$ to $p$ is defined as 
\begin{equation*}
D_{KL}(p;q) = \sum_{s\in\{0,1\}^N}p(s)\log\left( \frac{p(s)}{q(s)} \right) 
\end{equation*}
provided $p(s)=0$ only when $q(s)=0$. 
The Kullback-Leibler divergence is not a metric since it is not symmetric. 
However, metrics can be defined as follows:
the Jensen-Shannon metric $\jsmetric$ is defined by
\begin{equation*}
\jsmetric(p,q)^2 = \frac{1}{2}D_{KL}\left(p;\frac{p+q}{2}\right)+\frac{1}{2}D_{KL}\left(q;\frac{p+q}{2}\right)
\end{equation*}
and the Hellinger metric $d_H$ is defined by
\begin{equation*}
d_H(p,q)^2=\sum_{s\in\{0,1\}^N}\Big(\sqrt{p(s)}-\sqrt{q(s)}\Big)^2.
\end{equation*}
Recently, probabilistic metrics have been developed which are based on distributions of gene trees instead of distributions of characters \citep{adams2019}.

The Kullback-Leibler divergence, squared Jensen-Shannon metric and squared Hellinger metric are all examples of a more general class of distances between probability distributions known as $f$-divergences. 
Given any convex function $f(t)$ such that $f(1)=0$, the $f$-divergence of $p$ from $q$ is defined as
\begin{equation}\label{equ:fdiv}
D_f(p;q) = \sum_{s\in \{0,1\}^N} q(s)f\left( \frac{p(s)}{q(s)} \right).
\end{equation}
The Kullback-Leibler divergence $D_{KL}(p;q)$ is obtained by taking $f(t)=t\log t$, while the reversed divergence $D_{KL}(q;p)$ is obtained with $f(t) = -\log(t)$. 
The squared Jensen-Shannon metric and squared Hellinger metric can also be obtained by using more complicated functions $f$, cf. \cite{sason2016f}. 

Now, let $F\in\prew{N}$ be a forest representative of a wald $[F]\in \w{N}$.
As described in Section \ref{sec:waldspacedefinition}, the distributions at leaves of different trees of $F$ are independent. 
For two leaves in the same tree in $F$, some degree of evolution occurs between them since by definition of $\prew{N}$ no two leaves are coincident. 
Therefore, all characters are possible, giving
\begin{equation}\label{equ:prob-non-zero}
p_F(s)\neq 0\mbox{ for all }s\in \{0,1\}^N\,.
\end{equation}

It follows that the Kullback-Leibler divergence is always well-defined between distributions of characters corresponding to forest representatives from the \wald.
Since by Lemma~\ref{lem:edgeprod} the map $[F] \mapsto p_F$ is injective for $[F]\in \w{N}$ 
the Jensen-Shannon and Hellinger metrics pull back to define extrinsic metrics on the \wald $\w{N}$ (analogously to \citep{Garba2018}). 
As already mentioned in the introduction, statistical methods rely heavily on the intrinsic nature of metrics, and thus we aim for more geometrical structure in the next section by imposing the Fisher information metric (a Riemannian metric) onto the \wald.

\subsection{A two-state process geometry for the \wald}

\bhvspace $\unrooted{N}$ and \wald $\w{N}$ both do not have the structure of a manifold globally, but the interior of each maximal orthant is a manifold parametrized by $\vec{\ell}$ or $\vec{\lambda}$. Therefore we consider first the information geometry on the subspaces of \wald corresponding to a fixed fully resolved tree topology -- here every wald has only one single tree representative, since every wald corresponding to a forest with more than one component, as well as a wald containing a pendant edge with length zero, lies on the boundary of unit cubes corresponding to fully resolved tree topologies. Secondly, we establish global results about the constructed geometry of $\w{N}$.

Thus suppose $\tau$ is a fully resolved tree topology, and that trees with this topology are parametrized by $\vec{\ell}=(\ell^1,\ldots,\ell^{2N-3})\in \R_{>0}^N\times\orthant_\tau$. 
Let $p_{\vec{\ell}}(s)$ be the probability mass function $p_T(s)$ associated with tree $T$ determined by $\tau,\vec{\ell}$. 
Recalling that $p_{\vec{\ell}}(s)>0$ for all $s$, due to (\ref{equ:prob-non-zero}), the \emph{Fisher information matrix} at $\vec{\ell}$ is 
\begin{equation}\label{equ:TwoStateFisherInfo}
g_{\ia\ib}(\vec{\ell}) = \sum_{s\in \{0,1\}^N}p_{\vec{\ell}}(s)\Big(\partial_\ia\log p_{\vec{\ell}}(s)\Big)\Big(\partial_\ib\log p_{\vec{\ell}}(s)\Big)
\end{equation}  
for $\ia,\ib=1,\ldots,2N-3$ where $\partial_\ia = \partial/\partial_{\ell^\ia}$.  
This defines a Riemannian inner product on the tangent space of $\R_{>0}^N\times\orthant_\tau$ at $\vec{\ell}$ (that is a copy of $\R^{2N-3}$) which gives a way to measure the lengths of paths. 
Specifically, if $p_{\vec{\ell}}$ and $p_{\vec{\ell}+\delta\vec{\ell}}$ lie infinitesimally close on a path, then the squared path length between them is defined to be $\sum_{\ia,\ib}\delta\ell^\ia g_{\ia\ib}(\vec{\ell})\delta\ell^\ib$.  
Standard results from Riemannian differential geometry show that if $\vec{\ell}(t)$ is a path in $\R_{>0}^N\times\orthant_\tau$ then it is locally a geodesic (i.e.~it minimizes path length) if it satisfies the differential equation
\begin{equation}\label{equ:geodesic}
\frac{d^2\ell^\ic}{dt^2}+\sum_{i,j}\Gamma_{\ia\ib}^\ic(\vec{\ell})\frac{d\ell^\ia}{dt}\frac{d\ell^\ib}{dt} = 0,\quad \ic=1,\ldots,2N-3
\end{equation}
where $\Gamma_{\ia\ib}^\ic(\vec{\ell})$ are the Christoffel symbols
\begin{equation*}
\Gamma_{\ia\ib}^\ic(\vec{\ell}) = \sum_{\id}\frac{1}{2}g^{\ic\id}\left( \frac{\partial g_{\id\ia}}{\partial\ell^\ib} + \frac{\partial g_{\id\ib}}{\partial\ell^\ia} - \frac{\partial g_{\ia\ib}}{\partial\ell^\id} \right)\,.
\end{equation*}
The matrix $g^{\ia\ib}$ is the inverse of $g_{\ia\ib}$ i.e. $\sum_{\ic}g^{\ia\ic}g_{\ic\ib}=\delta^\ia_{\ib}$ where $\delta^\ia_{\ib}$ is the Kronecker delta. 
It is important to note that the geodesic equation and loci of solutions are invariant under changes of parametrization, and so the equations can be formulated using lengths $\ell^\ia$ or the weights $\lambda^\ia$.
On the boundary, however, this is no longer necessarily true. 

The Riemannian metric defined by the Fisher information matrix is related to the Kullback-Leibler divergence and other $f$-divergences by the following lemma.

\begin{lemma}\label{lem:KLDfromIP} 
Suppose $D_f$ is an $f$-divergence given by some convex function $f$ with $f(1)=0$, as defined by~$\eqref{equ:fdiv}$. 
Consider a small perturbation $\delta\vec{\ell}=(\delta\ell^1,\ldots,\delta\ell^{2N-3})$ of the edge lengths of a tree $(\tau,\vec{\ell})$. 
Then
\begin{equation}\label{equ:perturb}
\sum_{\ia,\ib}\delta\ell^\ia g_{\ia\ib}(\vec{\ell})\delta\ell^\ib =\frac{2}{f''(1)} D_{f}\big(p_{\vec{\ell}+\delta\vec{\ell}}; p_{\vec{\ell}}\big)+O\big(|\delta\vec{\ell}|^3\big)
\end{equation} 
where the error term consists of third-order products of the elements of $\delta\vec{\ell}$ and
\begin{equation}\label{equ:infinitesimal-fdiv}
D_{f}\big(p_{\vec{\ell}+\delta\vec{\ell}}; p_{\vec{\ell}}\big) =\frac{1}{2}f''(1) \sum_s \frac{\big(p_{\vec{\ell}+\delta\vec{\ell}}(s) -p_{\vec{\ell}}(s) \big)^2}{ p_{\vec{\ell}}(s)}+O\big(|\delta\vec{\ell}|^3\big).
\end{equation} 
In other words, the norm of the perturbation, as measured with respect to the Riemannian inner product, is proportional to the $f$-divergence of $p_{\vec{\ell}+\delta\vec{\ell}}$ from $p_{\vec{\ell}}$. 
\end{lemma}

The proof is given in the Appendix. 
Since the lemma applies to an arbitrary $f$-divergence, the term on the right-hand side of Equation~$\eqref{equ:perturb}$ can be the Kullback-Leibler divergence or the squared Jensen-Shannon metric, for example. 
Lemma~\ref{lem:KLDfromIP} gives the fundamental assumption behind the geometries we construct on $\w{N}$: that distances are locally measured by the infinitesimal Kullback-Leibler divergence between probability distributions associated with trees, or equivalently, by any $f$-divergence. 
As a corollary of Lemma~\ref{lem:KLDfromIP}, it follows that the metric defined by Equation~$\eqref{equ:TwoStateFisherInfo}$ is positive definite (i.e.~the metric is not semi-Riemannian). 
This is because the map from trees to distributions of characters is injective, and since $D_f(p;q)>0$ for all $p\neq q$, it follows that the right-hand side of Equation~$\eqref{equ:perturb}$ is strictly positive for all small non-zero perturbations $\delta\vec{\ell}$.

Let $\distrib{\{0,1\}^N}$ be the space of distributions on $\{0,1\}^N$. By Lemma~\ref{lem:edgeprod}, the map from $\w{N}$ to distributions of characters determines an embedding of $\w{N}$ in $\distrib{\{0,1\}^N}$. 
Given a metric $d$ on $\distrib{\{0,1\}^N}$, let $d^\ast$ be the induced intrinsic metric on $\w{N}$, as described in Section~\ref{sec:embeddings}. 


\begin{theorem}\label{thm:edgeprod-metric-space}
Let $d$ and $d_0$ be metrics on $\distrib{\{0,1\}^N}$ which are the square root of an $f$- and $f_0$-divergence, respectively. 
Then for any $[F],[G]\in\w{N}$:
\begin{enumerate}
\item $d^\ast([F],[G]) < \infty$ and thus $d^\ast$ is well-defined.
\item $d^*([F],[G])=c\cdot d^*_0([F],[G])$ for some constant $c>0$. 
\item Any path which realizes the distance $d^*([F],[G])$ is a solution of Equation~$\eqref{equ:geodesic}$ at any point in the interior of a maximal orthant. 
\end{enumerate}
\end{theorem}

The proof is given in the Appendix. 
Valid choices for the metric $d$ in Theorem~\ref{thm:edgeprod-metric-space} include the Jensen-Shannon metric or Hellinger metric.
Theorem~\ref{thm:edgeprod-metric-space} establishes $\w{N}$ with metric $d^*$ as a length space. 
Finiteness of $d^*$ shows, for example, that points in $\w{N}$ corresponding to disconnected forests (or equivalently, trees with infinite edge lengths) are at a finite distance away from orthant interiors. 
The second assertion implies a scaling of the induced intrinsic metric under changes of the function $f$, which, in turn, substantiates our conclusions drawn from Lemma~\ref{lem:KLDfromIP} that the geometry of $\w{N}$ induced by $d$ is invariant under the choice of $f$. 

\subsection{Numerical investigation of the geometry}
The geodesic equation~$\eqref{equ:geodesic}$ can be solved numerically on the interior of any maximal orthant given some initial conditions $\vec{\ell}(0)=\vec{\ell}_0$ and $d\vec{\ell}(0)/dt=\vec{v}_0$. 
As described in Section~\ref{sec:MarkovProc}, the first and second derivatives of $p_{\vec{\ell}}(s)$ with respect to the edge lengths $\ell^\ia$ can be computed analytically. 
Calculation of $g_{\ia\ib}(\vec{\ell})$ consists of a sum over all $2^N$ possible characters involving the derivatives of $p_{\vec{\ell}}(s)$.  
The Christoffel symbols can similarly be calculated as sums over characters. 
A fourth-order Runge-Kutta method was used to integrate the ODEs. 

This numerical scheme was used to construct and visualize geodesics on a single orthant in $\w{5}$.
The particular topology and edge lengths for the orthant are represented by the Newick string $((1:\ell^1,2:\ell^2):\ell^6,3:\ell^3,(4:\ell^4,5:\ell^5):\ell^7)$. 
Parameters $\ell^1,\ldots,\ell^5$ are the pendant edge lengths and $\ell^6,\ell^7$ the lengths of the two internal edges. 
We restricted to $N=5$ leaves in order to enable easy visualization of geodesics. 
Integration was stopped whenever any internal edge was assigned a length $\leq 0$, corresponding to the boundary of the orthant. 
If this occurred for a pendant edge, the pendant edge length was given value zero at that step, and integration was continued. 

\begin{figure}
\begin{tabular}{cc}
\includegraphics[clip, trim= 4.5cm 9.4cm 5cm 9.4cm, scale=0.5]{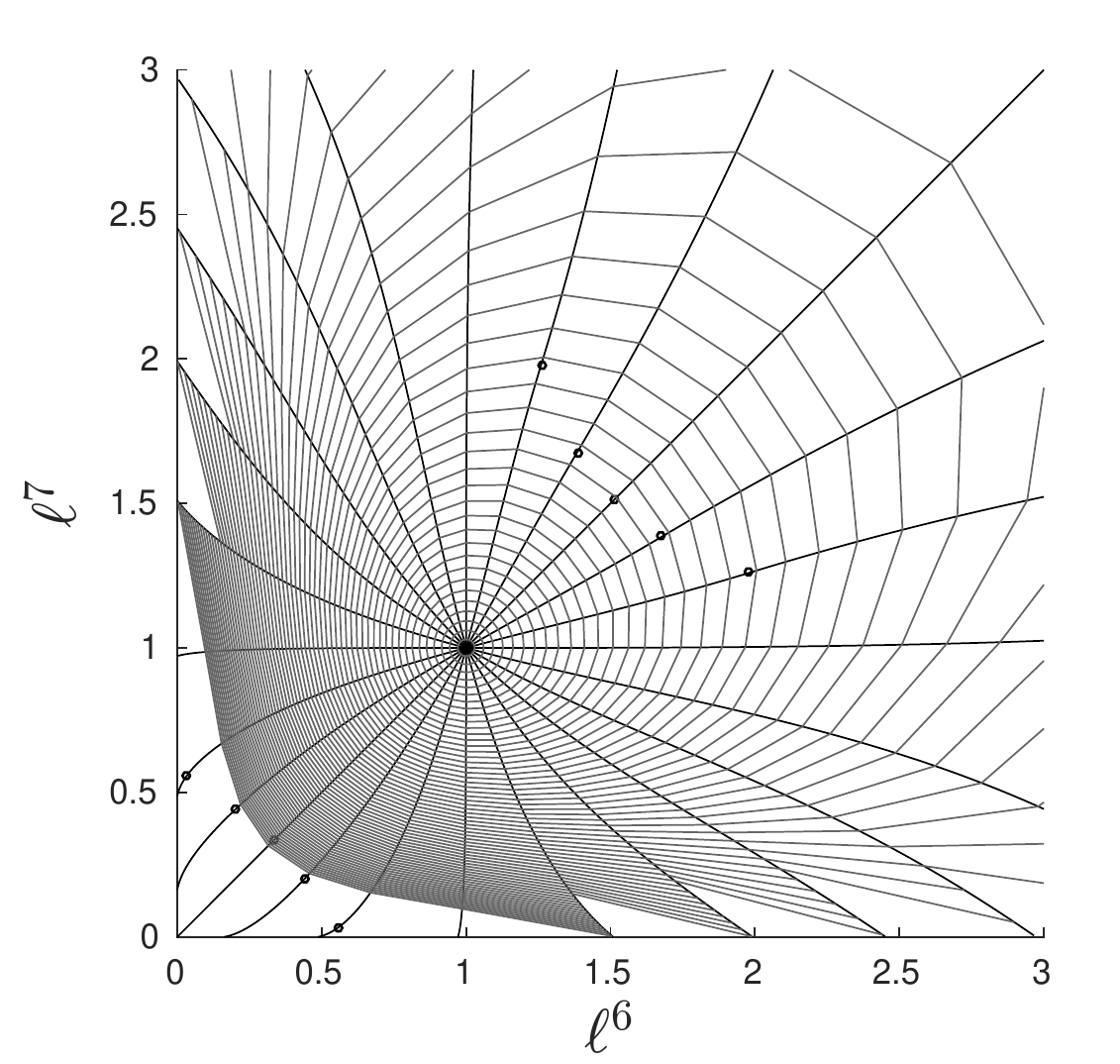} & \includegraphics[clip, trim= 4.5cm 9.4cm 5cm 9.4cm, scale=0.5]{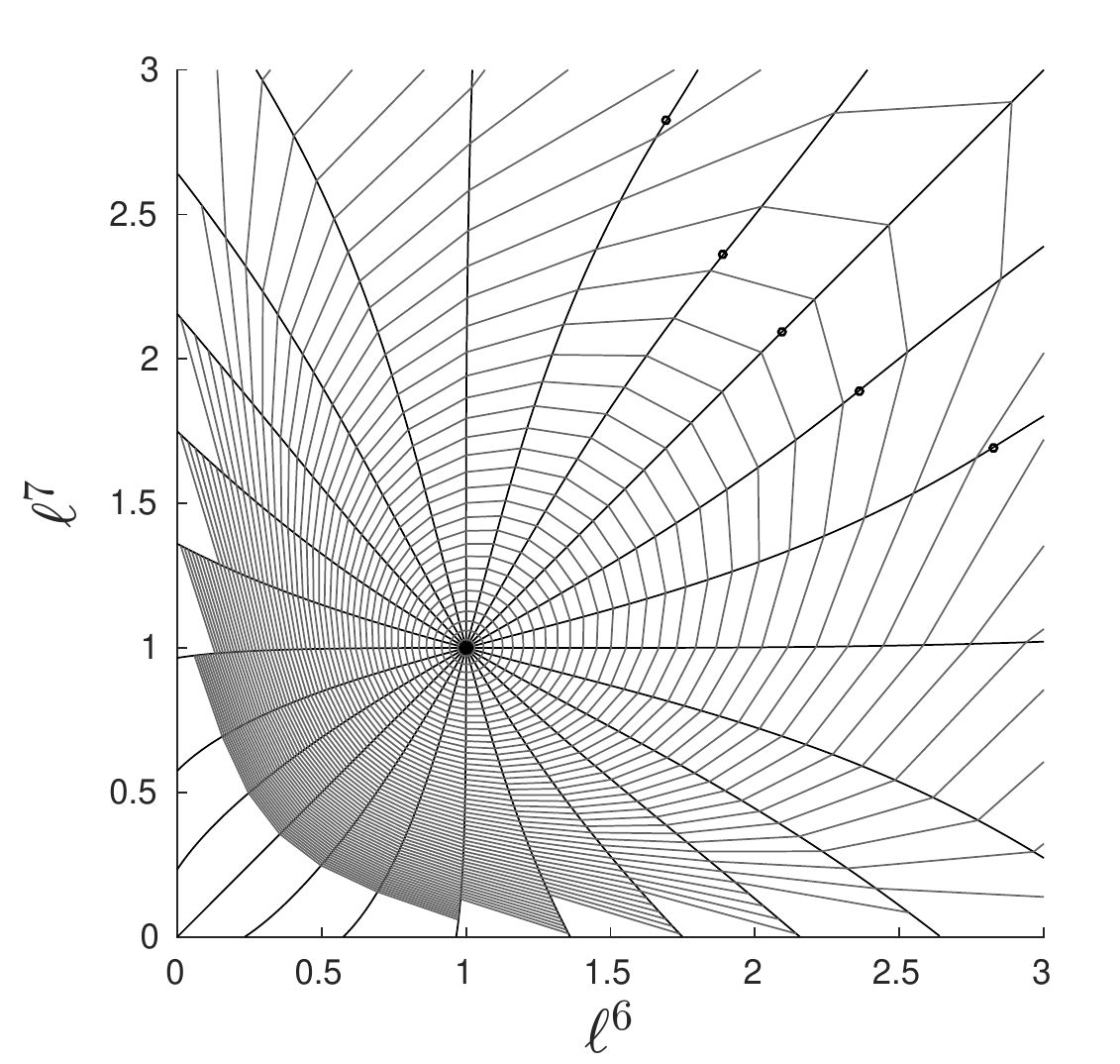} \\
\includegraphics[clip, trim= 4.5cm 9.4cm 5cm 9.4cm, scale=0.5]{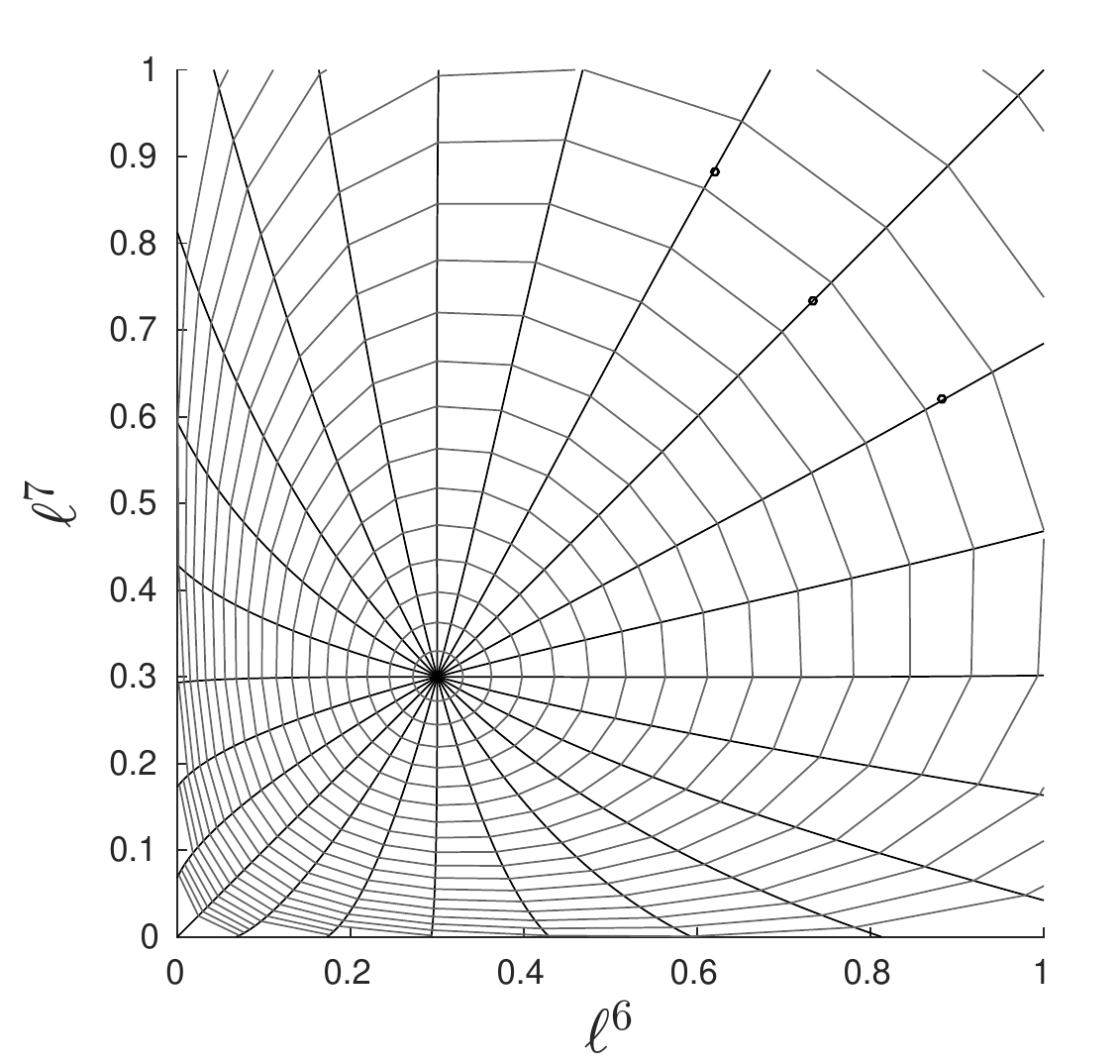} & \includegraphics[clip, trim= 4.5cm 9.4cm 5cm 9.4cm, scale=0.5]{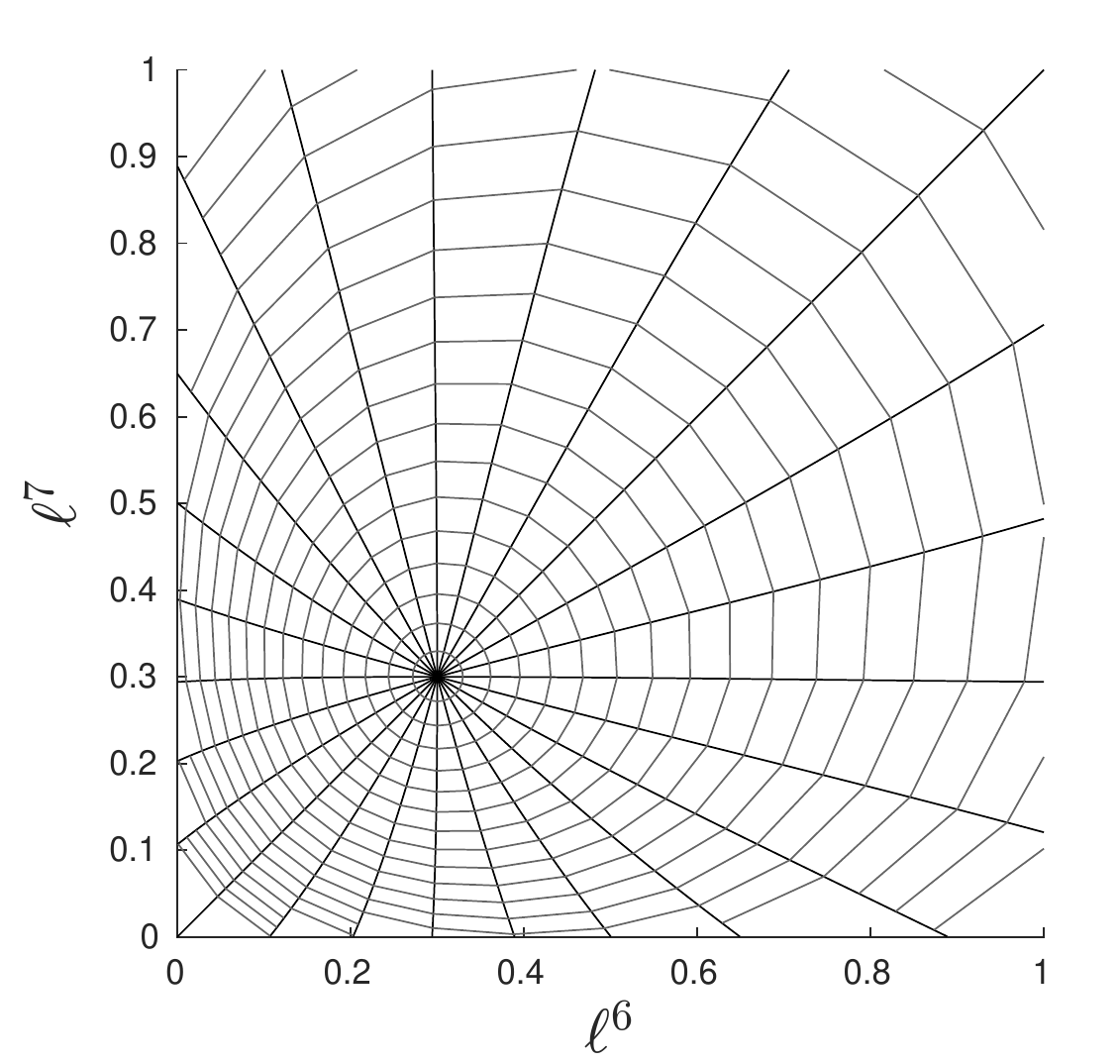} \\
\includegraphics[clip, trim= 4.5cm 9.4cm 5cm 9.4cm, scale=0.5]{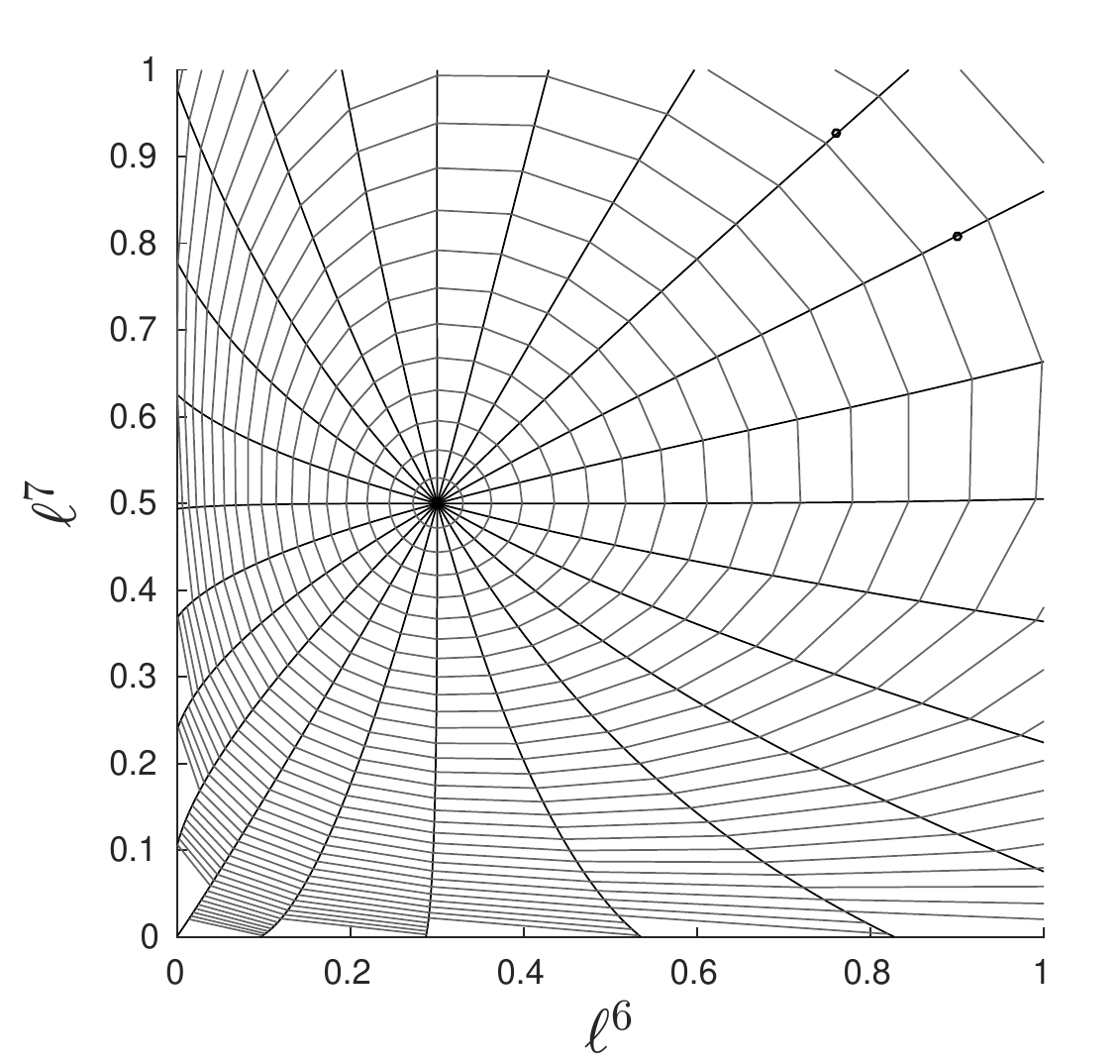} & \includegraphics[clip, trim= 4.5cm 9.4cm 5cm 9.4cm, scale=0.5]{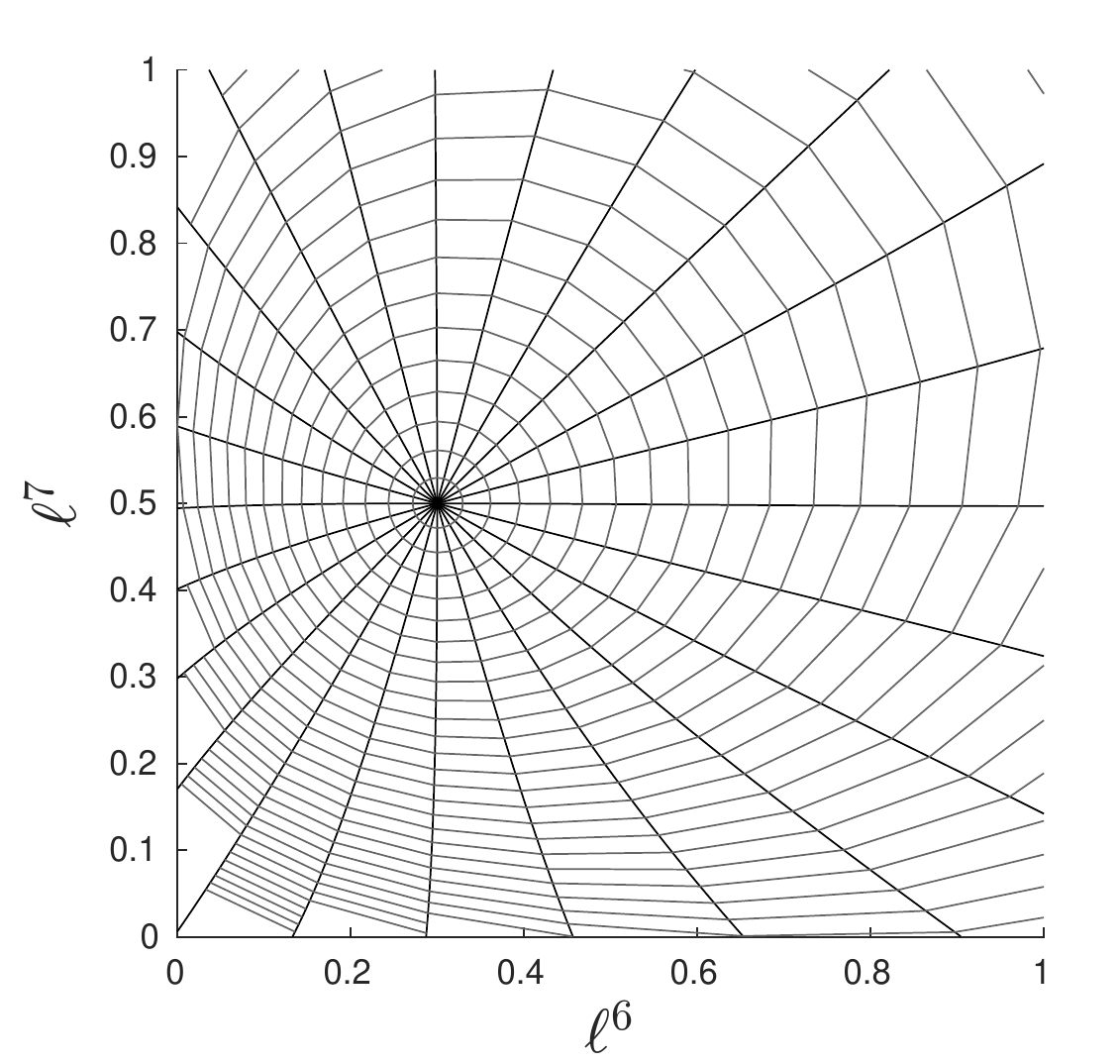} \\
\end{tabular}
\begin{caption}{\label{fig:TwoStateGeodesics}
Solutions to the geodesic equation (black radiating curves) and contours of distance (grey) in a single maximal orthant of $\w{N}$. 
Each panel shows the trajectories for the internal edge lengths $\ell^6,\ell^7$. 
Rows correspond to different initial sets of internal edge lengths. 
Columns correspond to different initial pendant edge lengths: $\ell^\ia=0.1$ for $\ia=1,\ldots,5$ on the left, and $\ell^\ia=0.5$ for $\ia=1,\ldots,5$ on the right. 
The initial velocity on pendant edge lengths was zero in all cases. 
Dots mark points at which a pendant edge was assigned length zero, and at all subsequent points, to avoid negative values. 
Contours near the origin in the top two plots have been removed: they stack up as the origin is approached and obscure the appearance of the geodesics. }
\end{caption}
\end{figure}

Figure~\ref{fig:TwoStateGeodesics} shows typical results. 
The figure shows the orthant representing the two internal edge lengths, with geodesics `fired' from some fixed starting point $\vec{\ell}_0$ in $24$ different directions. 
Also marked on the plots are contours of distance from the starting point. 
Each panel shows results for a different initial tree $\vec{\ell}_0$. 
It is evident that the geodesics are not the same as BHV geodesics, which are straight lines radiating from the initial point with equally spaced circular contours of distance. 
Figure~\ref{fig:TwoStateGeodesics} shows curved geodesics with irregularly spaced contours of distance. 
Contours appear to stack up towards the origin and codimension-1 BHV boundaries, but are more spaced out as geodesics move out towards the boundaries at infinity. 
This is more obvious when initial internal edges are long (top row of Figure~\ref{fig:TwoStateGeodesics}). 
On the other hand, the geodesics are more similar to the BHV geodesics when internal edges are short and pendant edges are long, as in the bottom two rows of the right hand column. In all cases, in contrast to \bhvspace, geodesics seem to be slightly attracted toward the star trees.  
The pendant edges do not behave as they do in \bhvspace: they can change value even when their initial velocity is zero. 

\begin{figure} 
\begin{center}
\begin{tabular}{cc}
\includegraphics[clip, trim=2.5cm 5cm 2.5cm 6cm, width=0.4\textwidth]{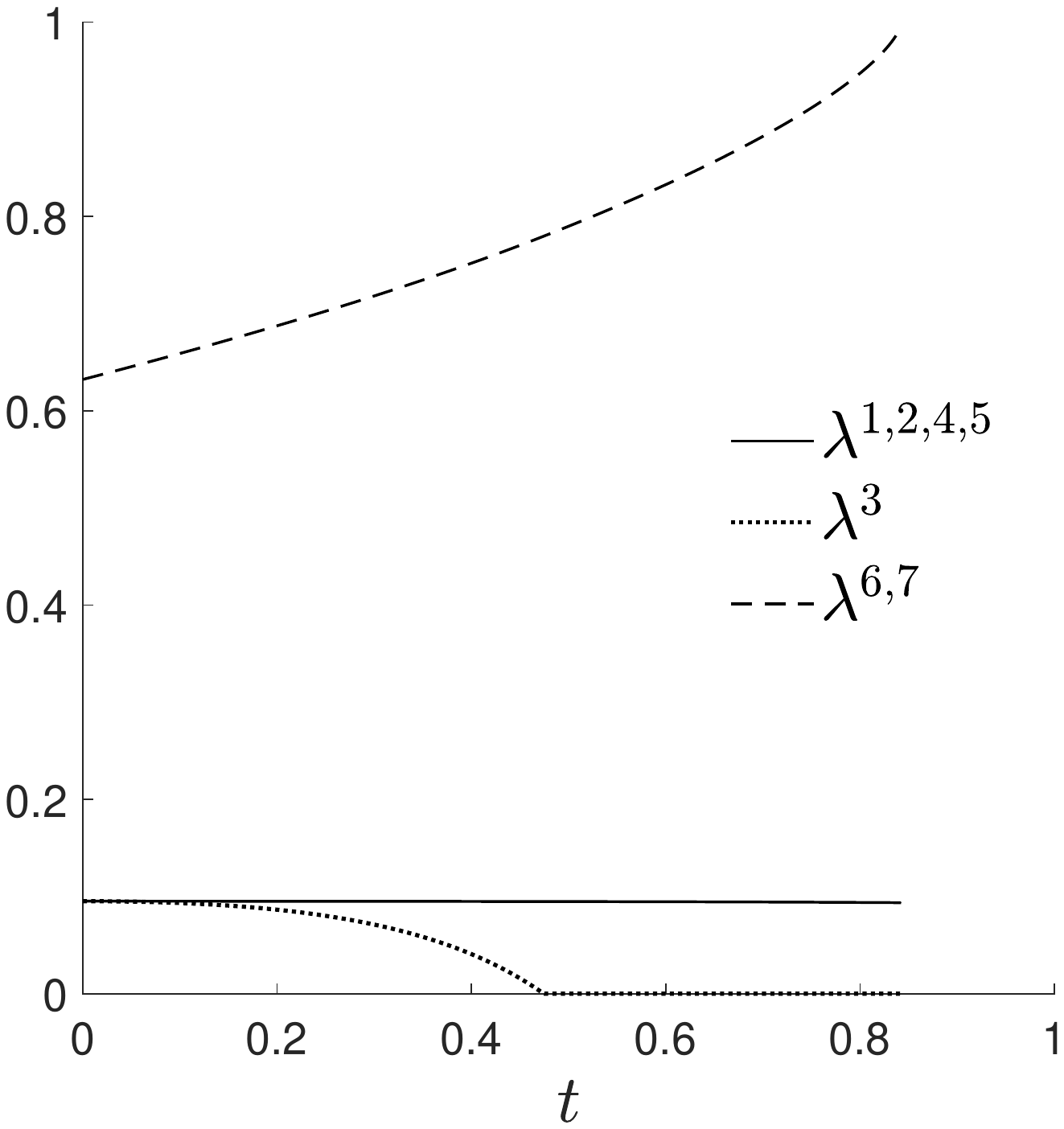} &
\includegraphics[clip, trim= 2.5cm 5cm 2.5cm 6cm,width=0.4\textwidth]{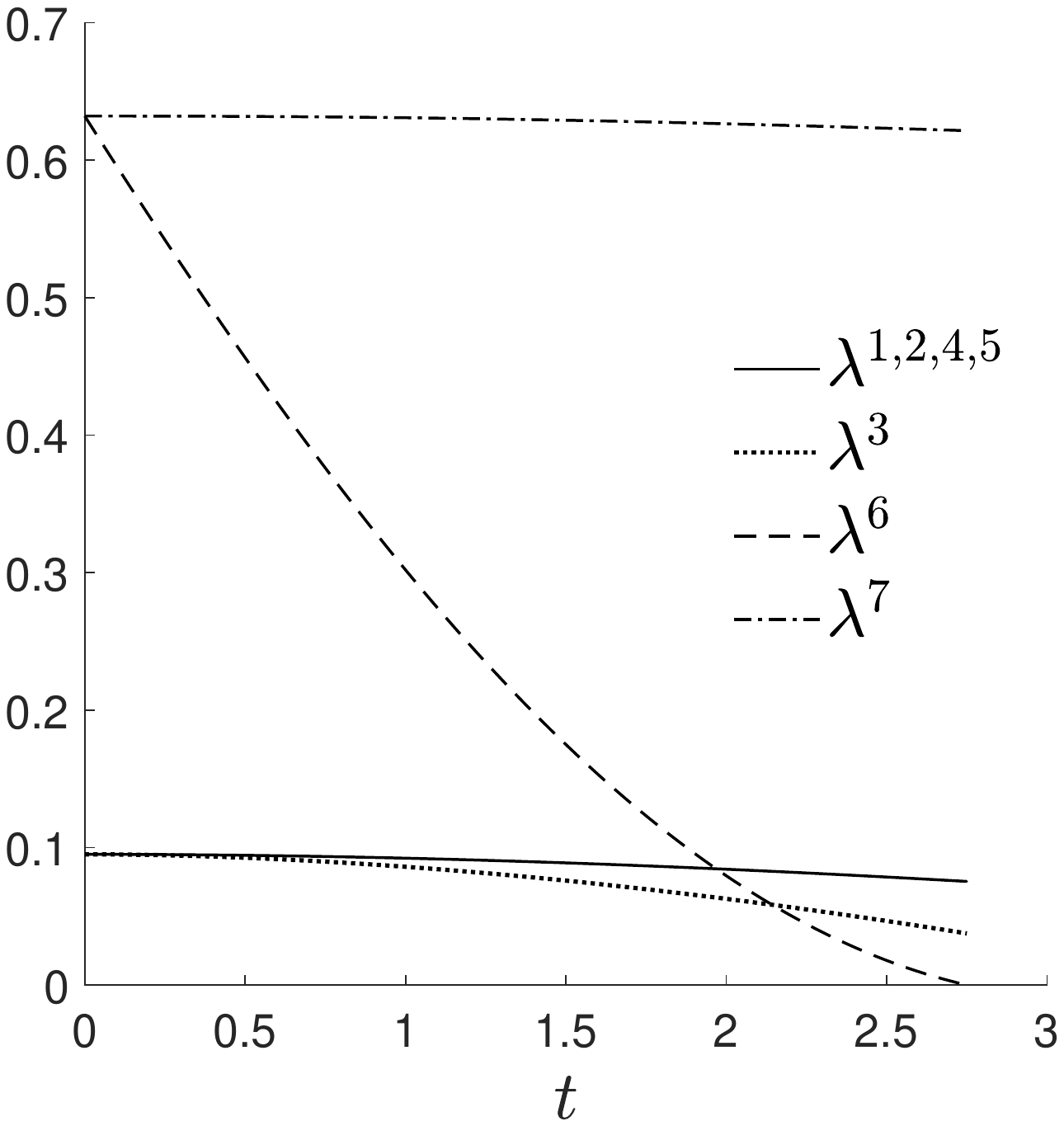} \\
(a) & (b)
\end{tabular}
\end{center}
\begin{caption}{\label{fig:GeodesicGraphs}
Graphs showing edge weights $\lambda^i$ vs time along geodesics in top left panel of Figure~\ref{fig:TwoStateGeodesics}. 
(a) Geodesic heading in North East compass direction. 
(b) Geodesic heading West. }
\end{caption}
\end{figure}

Figure~\ref{fig:GeodesicGraphs} provides more detail for certain geodesics in Figure~\ref{fig:TwoStateGeodesics}. 
The graphs in the figure show each edge weight $\lambda^i$ versus time, and the time is proportional to distance travelled. 
The $\vec{\lambda}$-parametrization was used for these plots since it shows the results most clearly. 
Panel (a) shows that, when contours become more and more widely spaced, the boundary at infinity ($\lambda^6=\lambda^7=1$) can be reached in finite time, rather than asymptotically. 
This shows that points corresponding to trees with infinitely long edges are a finite distance away from the starting point, as established by Theorem~\ref{thm:edgeprod-metric-space}. 
In the $\vec{\ell}$-parametrization, the internal edge lengths rapidly blow up to infinity after time $t=0.8$. 
It follows that the shortest path between two trees with finite edge lengths might involve trees with infinitely long edges. 
On the other hand, for some panels in Figure~\ref{fig:TwoStateGeodesics}, the contours become increasingly close as BHV boundaries are approached, but as panel (b) in Figure~\ref{fig:GeodesicGraphs} shows, points on BHV boundaries are in fact finitely close to orthant interiors, since the boundary is reached in finite time. 

\begin{figure}
\begin{tabular}{cc}
\includegraphics[clip, trim= 4.5cm 9.4cm 5cm 9.4cm, scale=0.5]{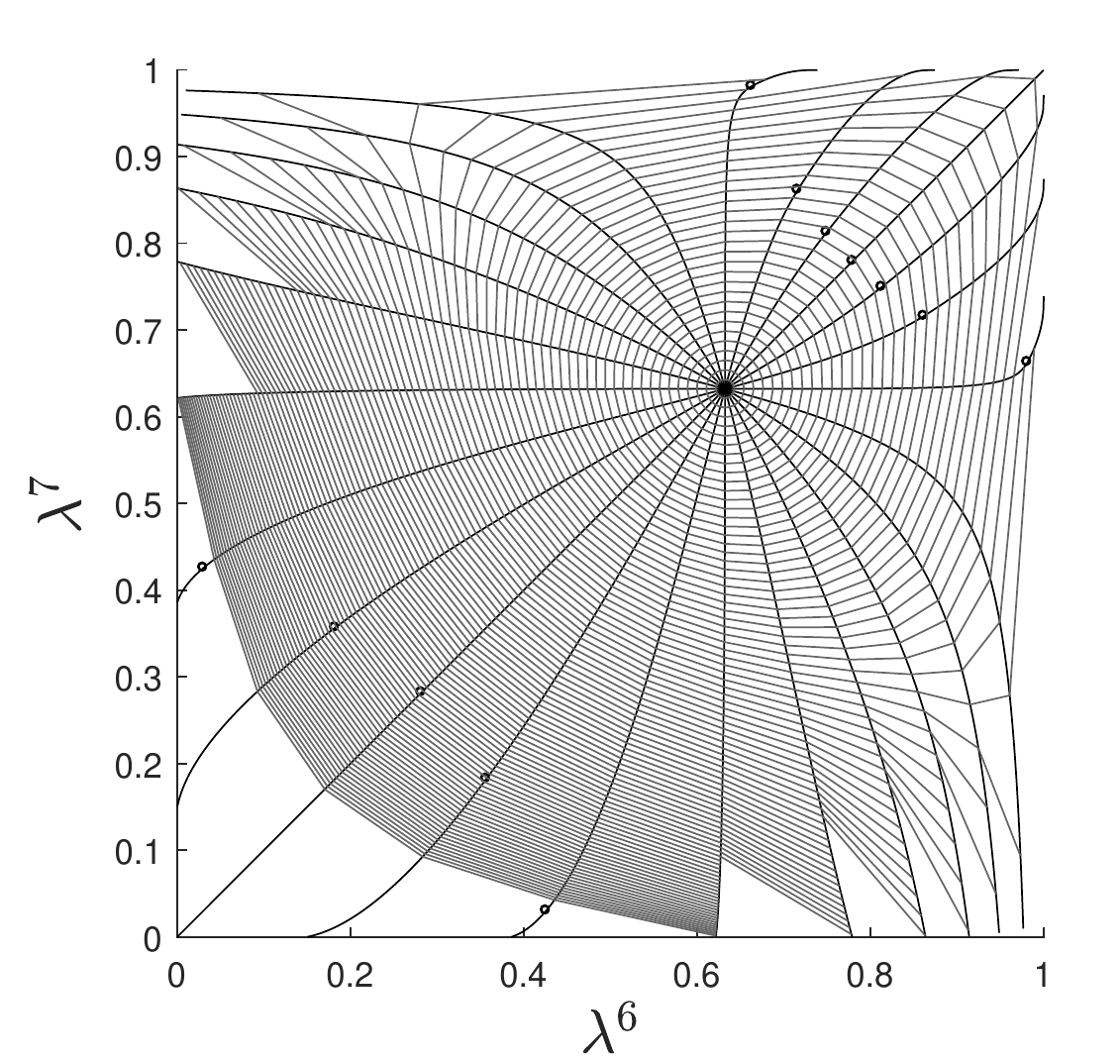} & \includegraphics[clip, trim= 4.5cm 9.4cm 5cm 9.4cm, scale=0.5]{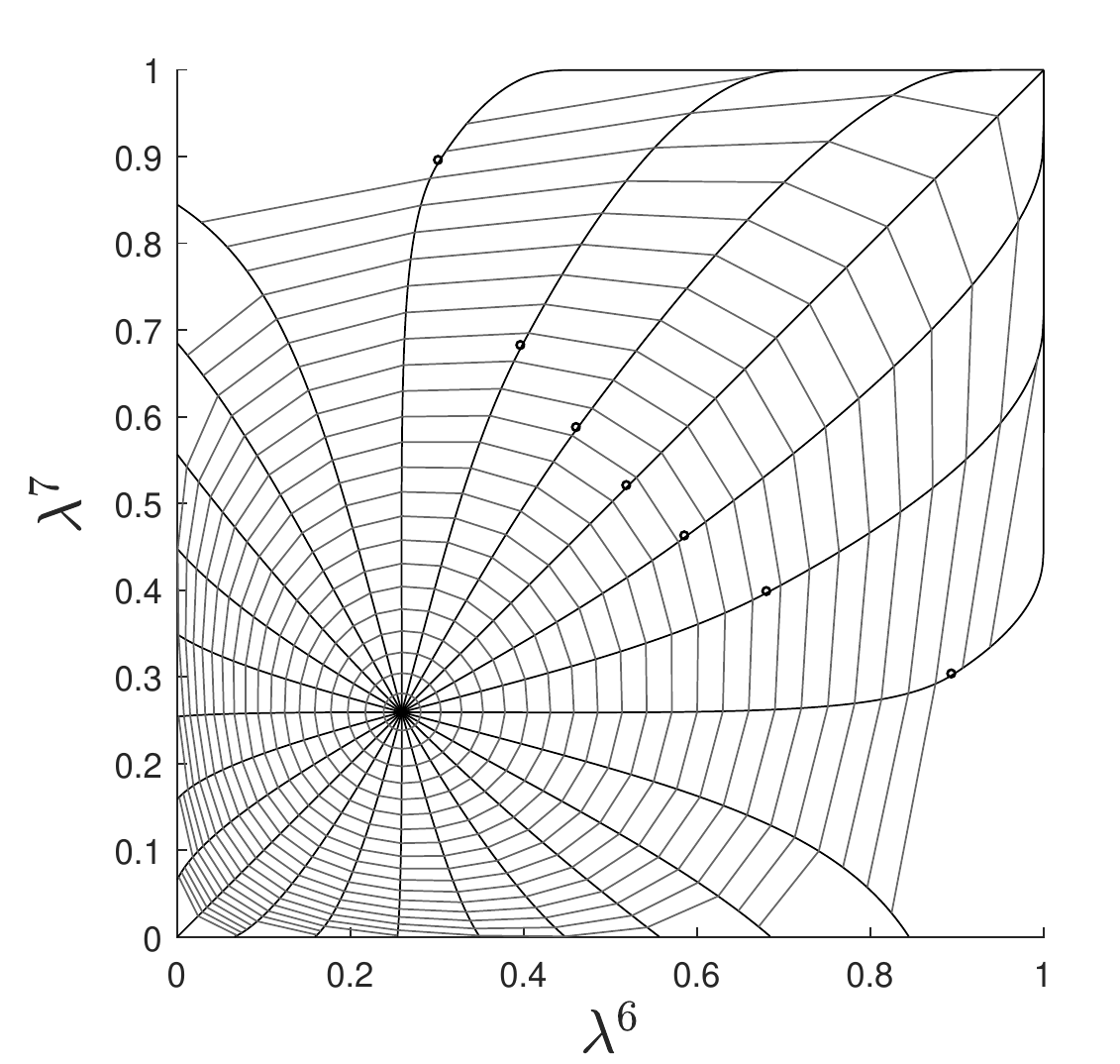} \\
(a) & (b)
\end{tabular}
\begin{caption}{\label{fig:cobweblambda}
Top left and middle left panels on Figure~\ref{fig:TwoStateGeodesics} redrawn in the $\vec{\lambda}$-parametrization.}
\end{caption}
\end{figure}

Figure~\ref{fig:cobweblambda} replots the top left and middle left panels from Figure~\ref{fig:TwoStateGeodesics} using the $\vec{\lambda}$-parametrization, so that the boundary at infinity corresponds to edges of the unit square with weight $1$.  
These plots suggest that trees in which one of the two internal edges is finitely long are `repellant' since the geodesics fired in the North and East compass directions end up passing through the disconnected forest with $\lambda^6=\lambda^7=1$. 
Indeed, the points on the boundary with $\lambda^6=\lambda^7=1$ appear to be `attractive', with geodesics pulled round to pass through these points and arriving in finite time. 

While these results show how the information geometry on $\w{N}$ differs substantially from the BHV geometry, the method for constructing geodesics by integrating the geodesic differential equation suffers from several disadvantages. 
First, it requires summation over the elements of $\{0,1\}^N$ which makes it infeasible for large $N$ (exponential computation time in $N$). 
Secondly, only local geodesics are computed by solving the initial value problem (geodesic ``shooting'' or ``firing'') for the differential equation valid only in maximal orthants. 
Thirdly, for practical applications, an algorithm which takes two points in \wald and joins them by a globally shortest path is more useful, and so it is desirable to solve the boundary value problem for geodesic construction. 
The next section attempts to deal with some of these shortcomings. 


\section{Information geometry for a Gaussian process on trees}\label{sec:gaussianprocess}

In this section we develop the information geometry of a continuous-valued Markov process associated to each tree which is more computationally and analytically tractable than the information geometry for the symmetric two-state Markov process.  
It has the advantage that the geodesic equation~$\eqref{equ:geodesic}$ can be solved numerically much faster than the corresponding equation for the symmetric two-state model, but the solutions for the two models are very similar.

\subsection{Definition of the Gaussian process}\label{sec:defgaussianproc}

Our aim is to construct a Gaussian process which is a continuous-valued analogue of the symmetric two-state process by matching the moments specified in Lemma~\ref{lem:twostatemoments}. 
Consider the Ornstein-Uhlenbeck process $Z(t)$ on $T$ which satisfies
\begin{equation}\label{equ:conttrans}
Z(t_2)\big{|}Z(t_1)=z\ \sim \ N\big(z\,e^{-\ell_{t_1t_2}}, 1-e^{-2\ell_{t_1t_2}}\big)
\end{equation}
where $\ell_{t_1t_2}$ is the path length distance between $t_1$ and $t_2$ on $T$. 
The stationary distribution is the standard normal distribution $N(0,1)$, and we assume the process is stationary over $T$. 
The Markov process satisfies the detailed balance equation, and so is reversible in its stationary state. 
As a result, realizations of the process can be simulated by fixing an arbitrary root $t_0\in T$ for the tree, simulating $Z(t_0)$ from $N(0,1)$ and then using Equation~$\eqref{equ:conttrans}$ to simulate a realization $Z(t)$ for any other $t\in T$. 
Reversibility of the process ensures the distribution obtained is invariant of the choice of root. 
A short calculation shows that the covariance matrix of the random variables $Z_1,\ldots,Z_N$ at the leaves of $T$ is given by $\cov[Z_\ileaf]{Z_\jleaf}=\exp(-\ell_{\ileaf\jleaf})$. 
Since $\ell_{\ileaf\ileaf}=0$ for all $\ileaf=1,\ldots,N$, this gives $\var{Z_\ileaf}=1$. 
Similarly, the conditional distribution of $Z(t_2)$ given $Z(t_1)=z$ has variance $1-\exp(-2\ell_{t_1t_2})$.  
These moments match those in Lemma~\ref{lem:twostatemoments} up to the constant factor of $1/4$, and we will show later that this factor makes no difference to the geometry obtained. 
The process $Z(t)$ can be thought of in two ways. 
First, it approximates the binomial random variables obtained when many independent binary characters evolve under the symmetric two-state process. 
Secondly, it could be regarded as an evolutionary model of a continuous trait for which the observations are standardized to have zero mean and unit variance in each population. 
Mean-reverting Gaussian processes like this are sometimes used to model continuous traits for which there is some constant optimal value for survival \citep{hansen1996}. 
The definition of $Z(t)$ extends from trees to forests by taking the process to be independent on each connected component. 

Given a forest $F\in \prew{N}$, the distribution of the random variables $Z_1,\ldots,Z_N$ at the leaves is multivariate normal with zero mean and covariance matrix $S_F$ where 
\begin{equation}\label{equ:covmtx}
S_F = \Big(\exp(-\ell_{\ileaf\jleaf})\Big)_{\ileaf,\jleaf=1}^N.
\end{equation} 
The path distance $\ell_{\ileaf\jleaf}$ is the sum of the lengths $\ell^e$ on edges $e$ between $\ileaf$ and $\jleaf$, and is taken to be infinite when $\ileaf$ and $\jleaf$ are in different components of $F$. 
This defines a map $F\mapsto \phi_F$ from forests to multivariate normal distributions where $\phi_F$ is the probability density function of $N(\vec{0},S_F)$. 
A similar result to Lemma~\ref{lem:edgeprod} holds: whenever $F_1\sim F_2$, the distributions $\phi_{F_1}$ and $\phi_{F_2}$ are the same, so the map is well-defined on $\w{N}$.

\subsection{A Gaussian process geometry for the \wald}\label{sec:continfogeom}

The information geometry of multivariate normal distributions with zero mean has been studied previously \citep{Lenglet2006} and is analytically tractable. 
The theory described in Section~\ref{sec:twostateinfogeom} needs adapting to account for the change from discrete to continuous characters. 
The Fisher information metric of $\w{N}$ in Equation~$\eqref{equ:TwoStateFisherInfo}$ becomes an integral over $\R^N$ rather than a sum, and the mass function $p_{\vec{\ell}}$ is replaced with the density function $\phi_{\vec{\ell}}$ of the Gaussian corresponding to a fully resolved tree with edge lengths $\vec{\ell}$:
\begin{equation}\label{equ:ContFisherInfo}
g_{\ia\ib}(\vec{\ell}) = \int_{\R^N}\phi_{\vec{\ell}}(s)\Big(\partial_\ia\log \phi_{\vec{\ell}}(s)\Big)\Big(\partial_\ib\log \phi_{\vec{\ell}}(s)\Big)\,ds.
\end{equation}  
The geodesic equation~$\eqref{equ:geodesic}$ remains the same. 
Although $p_T(s)$ and its derivatives could be evaluated exactly for the two-state symmetric model, evaluation of the metric and Christoffel symbols required a sum over all characters. 
For the continuous model, the corresponding integrals have closed form, as we describe below. 

Gaussian distributions with zero mean are parametrized by their covariance matrices, namely by the space of $N\times N$ symmetric positive definite matrices, which we will denote $\spd{N}$. 
The set of covariance matrices associated with forests forms a subset within $\spd{N}$, as determined by the following theorem. 
  
\begin{theorem}\label{thm:posdef}
\begin{enumerate}
\item The covariance matrix $S_F$ associated to any $F\in \prew{N}$, as defined by Equation~$\eqref{equ:covmtx}$, is positive definite so lies in $\spd{N}$, and
\item the map $[F]\mapsto S_F$ from $\w{N}$ to $\spd{N}$ is injective and so it determines a well-defined embedding of $\w{N}$ into $\spd{N}$.
\end{enumerate}
\end{theorem}

The proof is given in the Appendix.

For Gaussians with zero mean, it can be shown that the Fisher information metric at a point with covariance matrix $S$ is
\begin{equation*}
\big\langle X, Y \big\rangle = \frac{1}{2}\trace\left( S^{-1}XS^{-1}Y \right)\!,
\end{equation*}
where $X,Y$ are matrices in the tangent space at $S$, i.e. symmetric matrices \citep{Lenglet2006}. 
This expression is obtained by evaluating the integral in Equation~$\eqref{equ:ContFisherInfo}$.
Working in some fixed maximal orthant parametrized by edge lengths $\vec{\ell}$, let $S_{\vec{\ell}}$ be the corresponding covariance matrix defined in Equation~$\eqref{equ:covmtx}$. 
For each edge $e\in F$ define the \emph{split matrix} $\sigma^e$ by
\begin{equation}\label{equ:splitmatrix}
\sigma^e_{\ileaf\jleaf} = \begin{cases} 1, & \text{if $e$ lies on the path from leaf $\ileaf$ to leaf $\jleaf$, and}\\
0, &\text{otherwise,}
\end{cases}
\end{equation}
for $\ileaf,\jleaf=1,\ldots,N$. 
Then the path length between $\ileaf$ and $\jleaf$ is $\ell_{\ileaf\jleaf}=\sum_e\ell^e\sigma^e_{\ileaf\jleaf}$ where the sum is over all edges in $F$. 
Equation~$\eqref{equ:covmtx}$ becomes
\begin{equation}\label{equ:covHadamard}
S_F = S_{\vec{\ell}} = \bigg(\prod_e \exp\big( -\ell^e\sigma^e_{\ileaf\jleaf} \big)\bigg)_{\ileaf,\jleaf=1}^N.
\end{equation}
An entry above is zero if $\ileaf$ and $\jleaf$ are in different connected components, or equivalently, if they are separated by an infinitely long edge. 

By differentiating Equation~$\eqref{equ:covHadamard}$, it can be seen that the tangent space at $S_{\vec{\ell}}$ is spanned by matrices of the form $\sigma^e\circ S_{\vec{\ell}}$ for each edge $e$, where $\circ$ denotes the Hadamard matrix product. 
The Fisher information metric (\ref{equ:ContFisherInfo}) for $\w{N}$ for a fully resolved tree becomes
\begin{equation}\label{equ:FIMtrace}
g_{\ia\ib}(\vec{\ell}) = \frac{1}{2}\trace\Big( S_{\vec{\ell}}^{-1}\big(S_{\vec{\ell}}\circ\sigma^\ia\big)S_{\vec{\ell}}^{-1}\big(S_{\vec{\ell}}\circ\sigma^\ib\big) \Big)
\end{equation}
where $\ia,\ib=1,\ldots,2N-3$ index edges. 
Algebraic expressions for the first and second derivatives of the Fisher information metric can similarly be obtained, and hence for the Christoffel symbols. 
Note that scaling $S_{\vec{\ell}}$ by some positive constant has no effect on the metric, and so the factor of $1/4$ difference between the covariance matrices obtained from the discrete process $X(t)$ and continuous process $Z(t)$ has no effect on the geometry.  

The inner product and its derivatives can be computed in polynomial time and the paths obtained by integrating the geodesic ODE for the continuous Markov model within orthant interiors in $\w{N}$ resemble those for the two-state model very closely. 
Namely, results for the same initial conditions as Figure~\ref{fig:TwoStateGeodesics} were obtained but omitted, since the plots were almost indistinguishable from those for the two-state model. 
However, the inner products defined using the two different models are not identical: both inner products can be written down explicitly in the case $N=2$, using the transition probabilities for the discrete model and Equation~$\eqref{equ:FIMtrace}$ for the continuous model.  
The two inner products differ when the length of the single edge in the tree is small, but converge as the edge length tends to infinity.  

Using Equation~$\eqref{equ:FIMtrace}$ and its derivatives, we derived algebraic expressions for the Riemannian curvature tensor and the sectional curvatures at any point in $\w{N}$. 
We implemented these expressions in R, and hence evaluated these quantities for certain trees in $\w{5}$. 
We found that at randomly selected points in $\w{5}$, and hence in all spaces with $N\geq 5$, the sectional curvatures had mixed signs. 
As a result, there is no global sign condition on curvature like that for \bhvspace, which is globally non-positively curved.

\section{Geometry via embedding in $\spd{N}$}\label{sec:spdembedding}

The information geometry on Gaussians with zero mean can equivalently be regarded as a geometry for the space  of $N\times N$ symmetric positive definite matrices $\spd{N}$. 
This is a useful viewpoint to adopt, first because it highlights the fact that geometry we develop on $\w{N}$ is based entirely on the covariance between the leaves induced by the Markov process $Z(t)$, and secondly, because other metrics on $\spd{N}$ have been studied \citep{dryden2009}. 
These alternative metrics on $\spd{N}$ could in turn define different metrics on $\w{N}$, although they will not be considered any further in this paper. 
The metric on $\spd{N}$ obtained from the information geometry of Gaussian distributions with zero mean will be denoted $\covmetric$. 
The metric and its associated geodesics in $\spd{N}$ can be computed in polynomial time \citep{Lenglet2006}.
The main idea in this section is to use the analytically tractable geometry in $\spd{N}$, combined with a projection algorithm from the ambient space $\spd{N}$ to the embedded space $\w{N}$, to construct approximate geodesics in $\w{N}$.  

Given the embedding $[F]\mapsto S_F$ of $\w{N}$ within $\spd{N}$, we can consider the intrinsic metric $\inducedcovmetric$ on $\w{N}$ induced by $\covmetric$. 
By construction, the induced metric corresponds to the information geometry on $\w{N}$ for the continuous Markov model considered in Section~\ref{sec:gaussianprocess}.  
The following theorem is analogous to Theorem~\ref{thm:edgeprod-metric-space} for the discrete Markov substitution models. 
A proof is given in the appendix.

\begin{theorem}
    \label{thm:edgeprod-metric-space-cont}
    For any $[F],[G]\in\w{N}$ the induced intrinsic metric $\inducedcovmetric([F],[G])$ is finite and therefore well-defined. 
    Any path which realizes the distance $\inducedcovmetric([F],[G])$ is a solution of Equation~$\eqref{equ:geodesic}$ at any point in the interior of a maximal orthant, where the Riemannian inner product is given by Equation~$\eqref{equ:covHadamard}$. 
\end{theorem}

\citet{Lenglet2006} give formulae for the distance and geodesics between pairs of points in $\spd{N}$. 
The distance between $S_1,S_2\in\spd{N}$ is defined by
\begin{equation}\label{equ:covmetric}
\covmetric(S_1,S_2)^2 = \frac{1}{2}\trace\bigg(\log\Big( S_1^{-1/2}S_2 S_1^{-1/2} \Big)^2 \bigg)
\end{equation}
where $\log$ denotes the matrix logarithm. 
Since the map from $\w{N}$ to $\spd{N}$ is injective, $\covmetric$ pulls back to define an extrinsic metric on $\w{N}$:
\begin{equation}\label{equ:covmetrictrees}
\covmetric\big( [F_1],[F_2] \big) = \covmetric\big( S_{F_1}, S_{F_2} \big). 
\end{equation}
In fact, the space $\spd{N}$ equipped with $\covmetric$ has globally non-positive curvature \citep{Skovgaard1984, ballmann1985} and so there is a unique geodesic between any two points in the ambient space.
The point at proportion $t\in[0,1]$ along the geodesic between $S_1,S_2\in\spd{N}$ is 
\begin{equation}\label{equ:covgeodesic}
\Gamma_{S_1,S_2}(t) = S_1^{1/2}\exp\big(t\,U\big)\,S_1^{1/2}\quad \text{where}\quad U=\log\Big( S_1^{-1/2} S_2 S_1^{-1/2}\Big).
\end{equation}
Equations~$\eqref{equ:covmetric}$ to $\eqref{equ:covgeodesic}$ involve eigen-decompositions of $N\times N$ matrices, and so can be computed in $O(N^4)$ steps. 

\begin{figure} 
\begin{center}
\begin{overpic}[scale=0.7]{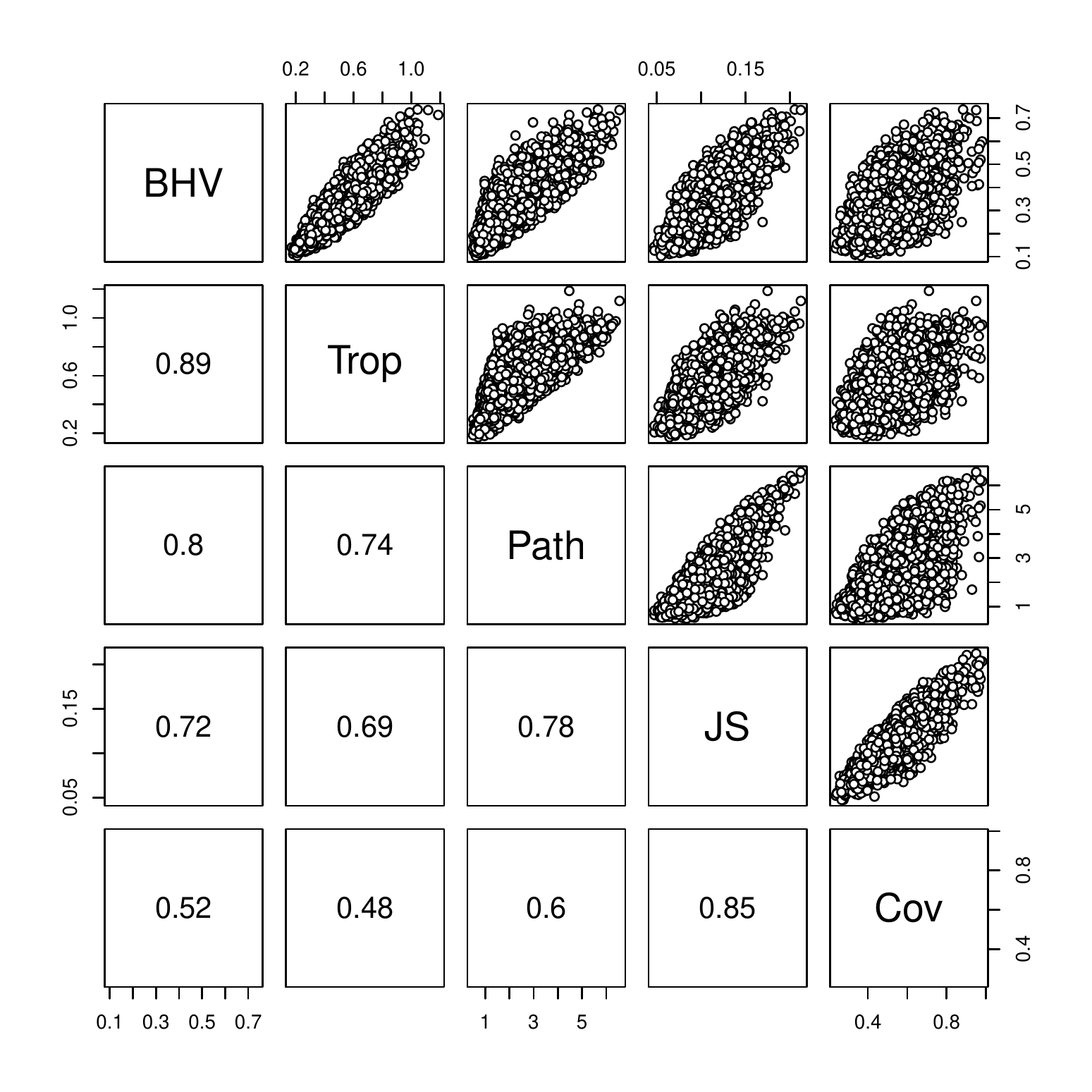}
\end{overpic}
\end{center}
\begin{caption}{\label{fig:metriccomparison}
Comparison of the BHV metric $\bhvmetric$, tropical metric, path difference metric, Jensen-Shannon metric $\jsmetric$ and $\covmetric$ for every pair from a sample of 100 trees obtained from bootstrap replicates during maximum likelihood estimation. 
The trees were inferred from DNA data from 12 species of primate. 
The correlation coefficients are shown in the bottom left panels.
}
\end{caption}
\end{figure}

Figure~\ref{fig:metriccomparison} shows a comparison of BHV metric $\bhvmetric$, tropical metric, path distance metric, Jensen-Shannon metric $\jsmetric$ and $\covmetric$ for every pair of trees in a sample of 100 trees obtained by bootstrap replication during maximum likelihood inference of a phylogenetic tree.   
The trees were inferred using the MrBayes software \citep{huelsenbeck2001mrbayes}, and a sample data set of DNA from 12 primates provided with the software. 
The path difference metric \citep{steel1993} between trees $T,T'$ is $(\sum_{u,v} (\ell_{uv}-\ell_{uv}')^2)^{1/2}$. 
The Jensen-Shannon metric was calculated exactly by summing over all $2^{12}$ characters for the two-state model, as described in Section~\ref{sec:probmetrics}. 
The covariance metric was calculated using Equations~$\ref{equ:covmetric}$ and~$\ref{equ:covmetrictrees}$. 
The BHV and tropical metrics are quite closely correlated, as are the Jenson-Shannon and covariance matrices. 
The path difference metric has a similar correlation with the BHV, tropical and Jensen-Shannon metrics (approximately $0.7$--$0.8$), but is relatively weakly correlated with $\covmetric$. 
This suggests that the BHV and tropical metrics are based on features of the data which are rather distinct from those for the Jensen-Shannon and covariance metrics. 
The covariance metric has the advantage over the Jensen-Shannon metric of being computable in polynomial time. 

\subsection{Projection into \wald}
\label{sec:projection}

To approximate geodesics in the extrinsic covariance metric $\covmetric$ of $\w{N}$, we construct a projection from $\spd{N}$ onto $\w{N}$, that is, given $S_0\in\spd{N}$, we aim to find an element $[F]\in\w{N}$ which minimizes $\covmetric(S_0,S_F)$. 
Suppose that $F$ is a fully resolved tree with edge lengths $\vec{\ell}$. 
The expression for $\covmetric(S_0,S_F)^2$ can be differentiated with respect to edge lengths of $F$ and gives
\begin{equation*}
\partial_i\covmetric\big(S_0,S_F\big)^2 = \trace\bigg(\log\left( S_0^{-1/2}S_F S_0^{-1/2} \right) S_0^{1/2}S_F^{-1}\big( \partial_i S_F \big) S_0^{-1/2} \bigg)
\end{equation*} 
where $\partial_i=\partial/\partial_{\ell^i}$ (e.g. \cite{moakher_differential_2005}). 
Moreover, $\partial_i S_F= S_F\circ\,\sigma^i$ where $\circ$ denotes the Hadamard or element-wise matrix product and $\sigma^i$ is the split matrix associated with edge $i$, as defined in Equation~$\eqref{equ:splitmatrix}$. 

This analytic expression for the derivative can be used to implement a gradient descent algorithm.
Within each maximal orthant $\orthant_\tau$ the edge lengths were updated according to the rule
\begin{equation*}
\vec{\ell}_{k+1} = \vec{\ell}_k-\alpha_k\,\nabla \covmetric\big(S_0,S_{\vec{\ell}_k}\big)^2
\end{equation*}
where $\vec{\ell}_k$ denotes the edge lengths at iteration $k$ and $S_{\vec{\ell}_k}$ the corresponding covariance matrix. 
The step size $\alpha_k$ was determined using the Barzillai-Borwein method. 
Two versions of the algorithm were used. 
In the first, the algorithm was halted whenever an internal edge was assigned a negative length. 
As a result, the algorithm was constrained to lie within the orthant $\orthant_\tau$ containing the initial tree. 
This algorithm was used for $N=5$ by running the algorithm $15$ times, each time with an initial tree in one of the 15 maximal orthants of $\unrooted{5}$, and the overall tree closest to $S_0$ found.  
The second version of the algorithm was able to cross codimension-$1$ BHV boundaries as follows. 
If an edge length was assigned a negative value, then trees in the two corresponding neighbouring orthants to $\tau$ were considered, taking the absolute value of elements in $\vec{\ell}_{k+1}$ as edge lengths.
The tree at step $k+1$ was taken to be whichever of these two trees was closest to $S_0$. 

\begin{figure} 
\begin{center}
\begin{overpic}[scale=0.8]{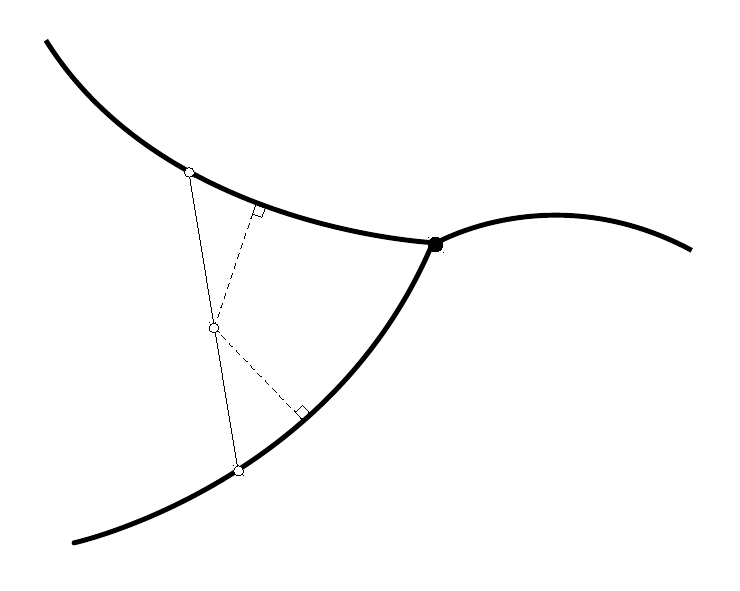}
\put (25,60) {$[F_1]$}
\put (30,9) {$[F_2]$}
\put (35,57) {$[G_1]$}
\put (42,18) {$[G_2]$}
\put (20,36) {$S$}
\end{overpic}
\end{center}
\begin{caption}{\label{fig:discont}
Schematic diagram for $N=4$ showing three neighbouring orthants (curved heavy lines) embedded in $\spd{N}$. 
The black circle is the BHV boundary between the orthants. 
The extrinsic geodesic between trees $[F_1],[F_2]$ is depicted as a straight line segment. 
The projection of this segment consists of a path from $[F_1]$ to $[G_1$] within the orthant, but then jumps to $[G_2]$ in a different orthant. 
The dashed lines show the orthogonal projection of the point $S$ along the extrinsic geodesic, and $S$ is equidistant from $[G_1]$ and $[G_2]$.   
}
\end{caption}
\end{figure}

In general, the closest point in $\w{N}$ to a covariance matrix $S_0\in\spd{N}$ is not necessarily unique as illustrated in Figure~\ref{fig:discont}.
Moreover, the gradient descent algorithm can converge to local minima, and so the result obtained is sensitive to the tree used to initialize the algorithm. 

\subsection{Construction of geodesics in $\w{N}$ via projection of extrinsic geodesics}\label{sec:projgeod}

Since construction of geodesics in $\spd{N}$ between any two given points and projection from $\spd{N}$ into $\w{N}$ can both be performed efficiently, we aim to combine these algorithms to give an efficient means of constructing geodesics within $\w{N}$ between any two given end points. 
A naive approach, given $[F_1],[F_2]\in\w{N}$, is to simply project the extrinsic geodesic between $S_{F_1}$ and $S_{F_2}$ from $\spd{N}$ into $\w{N}$. 
This approach works for the example of the unit sphere $S^2$ embedded within $\R^3$: the projection of the chord between two points in the sphere is a great circle between those two points. 
However, this approach fails for $\w{N}\subseteq\spd{N}$ since the projected paths are often discontinuous and jump between different orthants, as illustrated in Figure~\ref{fig:discont}. 

The following \textbf{recursive algorithm} for constructing an approximate geo-desic in $\w{N}$ is therefore proposed, which is intended to overcome this issue. 
Let $t_i=i/k$ for $i=0,\ldots,k$ and suppose $[F_1],[F_2]\in\w{N}$. 
The algorithm outputs a sequence $[G_0],\ldots,[G_k]\in\w{N}$ where $[G_0]=[F_1]$ and $[G_k]=[F_2]$. 
For each iteration $i=1,\ldots,k-1$ of the algorithm, the following steps are performed. 
\begin{enumerate}
\item Compute the extrinsic geodesic $\Gamma$ from $S_{G_{i-1}}$ to $S_{G_k}$ using Equation~$\eqref{equ:covgeodesic}$. 
\item Find the point $S\in\spd{N}$ at proportion $1/(k-i+1)$ along $\Gamma$. 
\item Let $[G_{i}]$ be the projection of $S$ into $\w{N}$. 
\end{enumerate}

The idea is that at each iteration, a new extrinsic geodesic is constructed from the previous point $[G_{i-1}]$ to the destination $[F_2]$, a small step is taken along that geodesic, and that point is projected into $\w{N}$ to give $[G_i]$. 
For the results in this paper, the projection at Step 3 was performed using the second version of the projection algorithm described in Section~\ref{sec:projection}, rather than using the less efficient search over all orthants. 
The gradient descent for the projection to obtain $[G_i]$ at Step 3 was initialized using the edge lengths from the forest $[G_{i-1}]$. 

This algorithm has the disadvantage that it is not symmetric under swapping the end points $[F_1],[F_2]$, whereas the image of the geodesic should be invariant under this operation.

The following \textbf{symmetrized version} of the algorithm overcomes this issue.
The algorithm produces a sequence $[G_0],\ldots,[G_k],[H_k],\ldots,[H_0]\in\w{N}$ where the initial values are taken to be $[G_0]=[F_1]$ and $[H_0]=[F_2]$. 
For each iteration $i=1,\ldots,k-1$ of the algorithm, the following steps are performed. 
\vspace{-1ex}
\begin{enumerate}
\item Compute the extrinsic geodesic $\Gamma$ from $S_{G_{i-1}}$ to $S_{H_{i-1}}$ using Equation~$\eqref{equ:covgeodesic}$. 
\item Find the points $R,S\in\spd{N}$ at proportions $1/(k-i+1)$ and $1-1/(k-i+1)$ along $\Gamma$. 
\item Let $[G_{i}]$ and $[H_i]$ be the projections of $R$ and $S$ into $\w{N}$ respectively. 
\end{enumerate}

\begin{figure}
\begin{tabular}{cc}
\includegraphics[clip, trim= 4.5cm 9.4cm 5cm 9.4cm, scale=0.5]{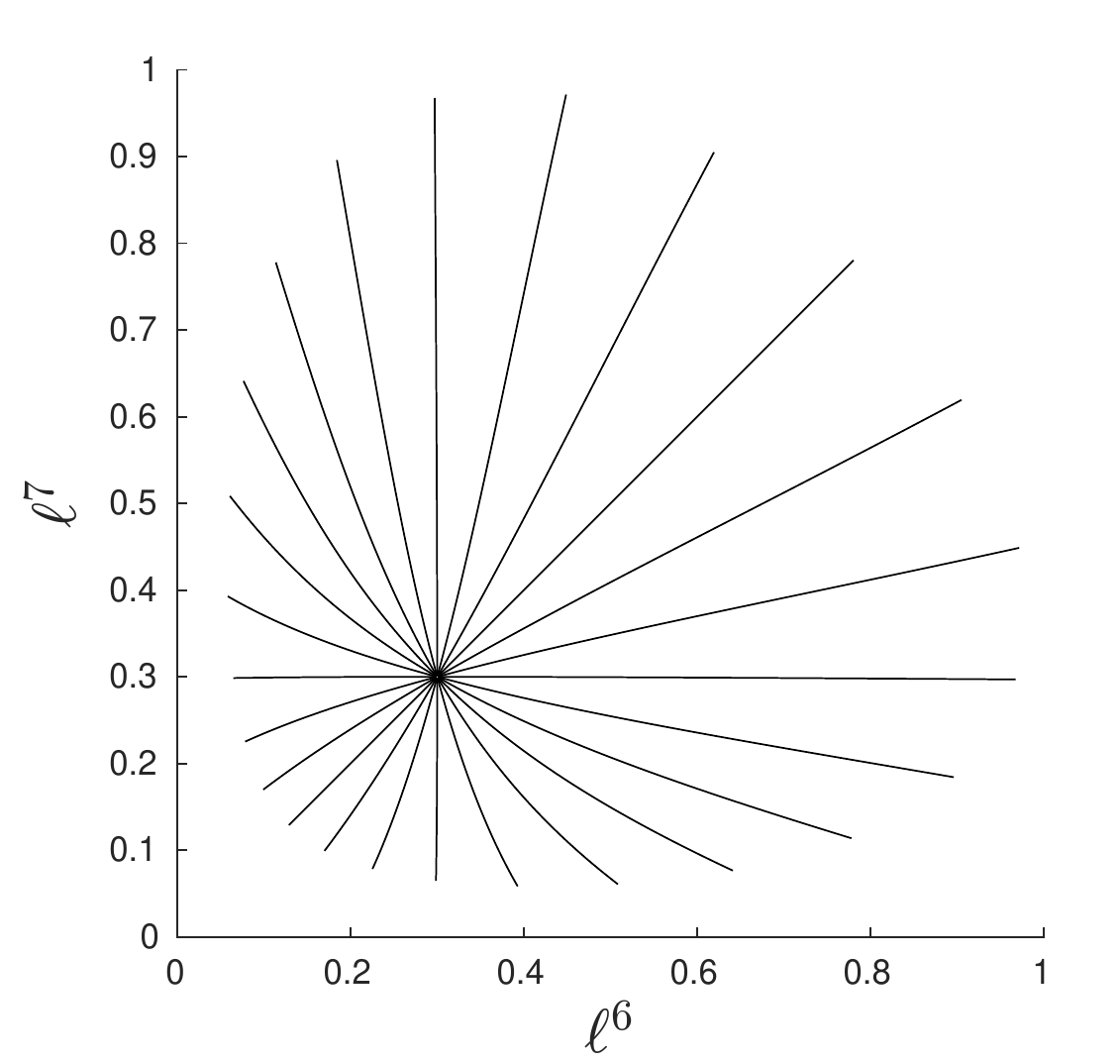} & \includegraphics[clip, trim= 4.5cm 9.4cm 5cm 9.4cm, scale=0.5]{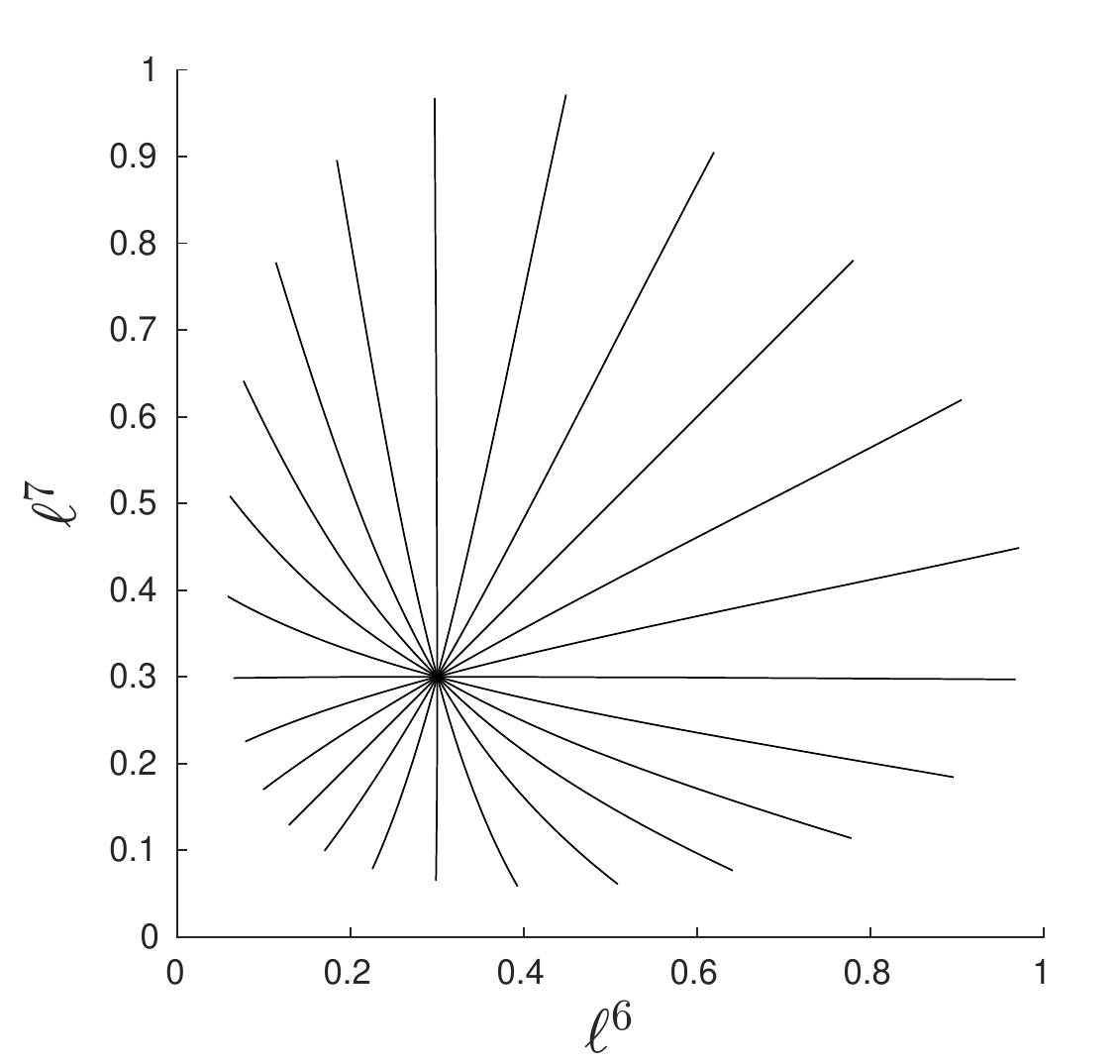} \\
(a) & (b)
\end{tabular}
\begin{caption}{\label{fig:single_orthant_comparison}
Comparison of paths obtained by (a) integrating the geodesic ODE for the Gaussian process model and (b) by applying the symmetrized projection method to the end points obtained in (a). }
\end{caption}
\end{figure}

The quality of the approximate geodesics produced by the symmetrized algorithm can be assessed by comparing them to geodesics in a single orthant constructed by integrating the geodesic equation as described in Section~\ref{sec:continfogeom}. 
Given an initial tree $[F_1]\in\w{5}$ and an initial velocity for $\vec{\ell}$, the geodesic equation was integrated until the path obtained reached some specified length. 
The final point reached was taken to be $[F_2]$, and the symmetrized algorithm was then used to obtain an approximate geodesic between $[F_1]$ and $[F_2]$. 
In all cases, the paths obtained with the two methods matched very closely, with the quality of the match improving for shorter internal edge lengths. 
Figure~\ref{fig:single_orthant_comparison} shows typical results.

\begin{figure} 
\begin{center}
\includegraphics[scale=0.45]{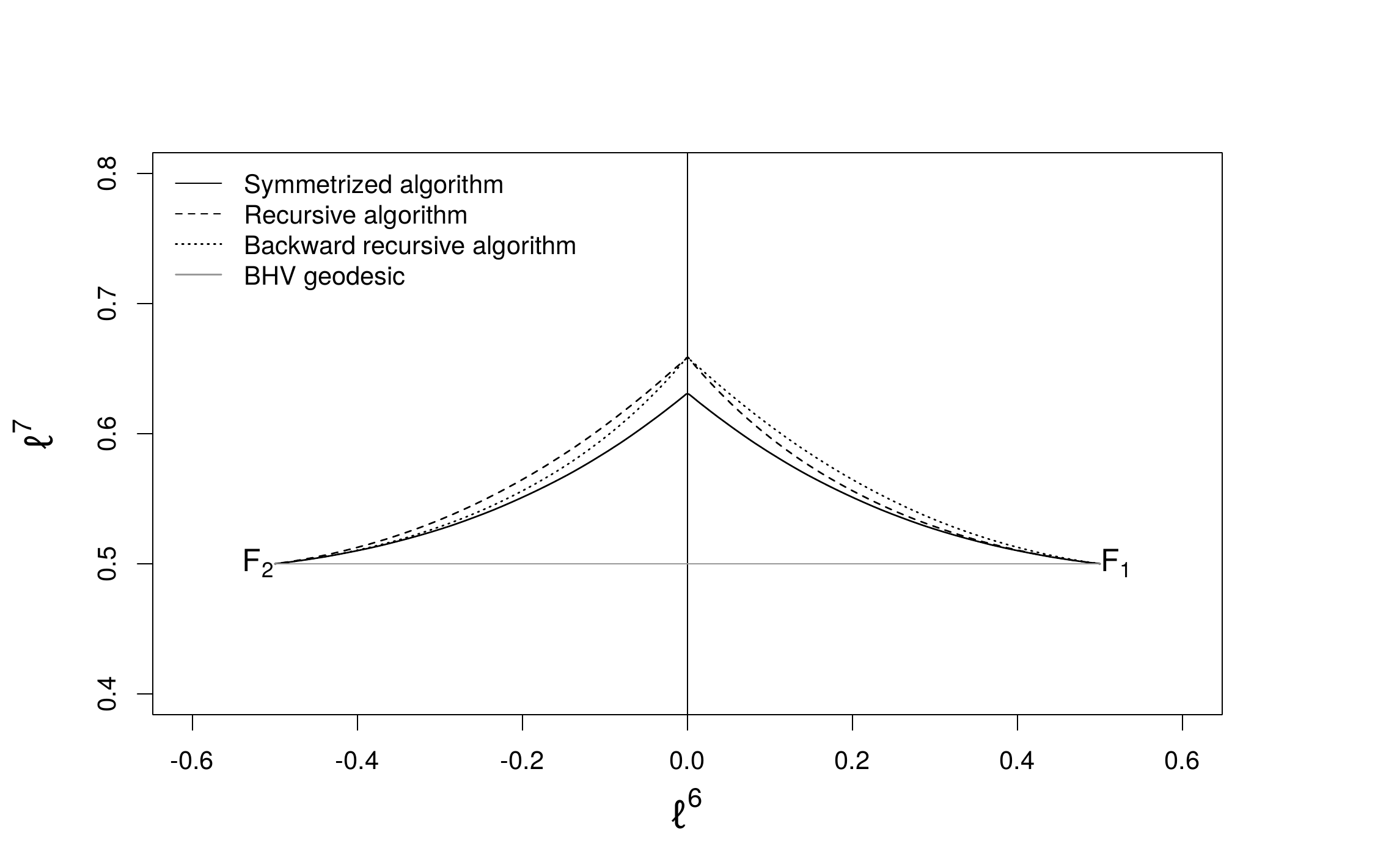}
\end{center}
\begin{caption}{\label{fig:proj_alg_2orthants}
Comparison of approximate geodesics in $\w{5}$ constructed beween trees $F_1$ and $F_2$ from (\ref{eq:two-trees-example}) in neighbouring orthants.   
The vertical axis $\ell^7$ represents a codimension-1 BHV boundary between two orthants. When, due to a nearest neighbor interchange, crossing it, $\ell^6$ tends to zero, another edge appears, with negative length corresponding to the negative values on the $\ell^6$ axis. 
Three approximate geodesics are shown: (i) construction via the recursive algorithm from $F_1$ and $F_2$, (ii) using the same algorithm but reversing the end-points, and (iii) construction via the symmetrized algorithm. 
}
\end{caption}
\end{figure}

\begin{figure} 
\begin{center}
\includegraphics[scale=0.55]{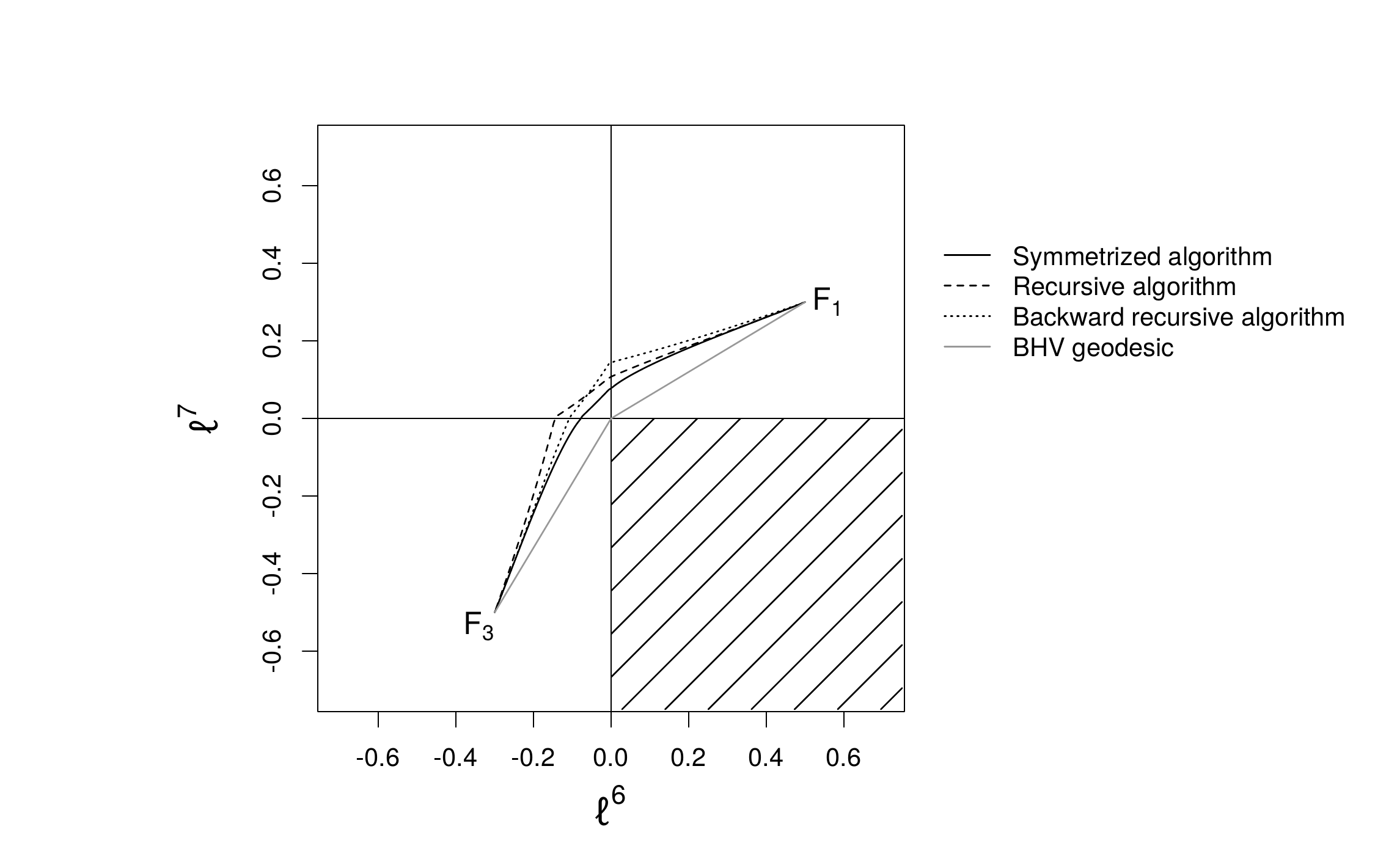}
\end{center}
\begin{caption}{\label{fig:proj_alg_3orthants}
Comparison of approximate geodesics constructed between trees $F_1$ and $F_3$ from (\ref{eq:two-trees-example}) in orthants separated by two nearest neighbour interchanges. 
Three neighbouring orthants in $\w{5}$ are shown, and the bottom right-hand orthant does not correspond to a valid tree topology. As in Figure \ref{fig:proj_alg_2orthants}, negative values on axes correspond to negative lengths of new edges. 
}
\end{caption}
\end{figure}

In contrast to the methods presented in Sections~\ref{sec:twostateinfogeom} and~\ref{sec:gaussianprocess}, these algorithms are not based on `firing' geodesics, and can produce approximate geodesics between end points in different orthants. 
Figures~\ref{fig:proj_alg_2orthants} and~\ref{fig:proj_alg_3orthants} show results obtained when the end points are separated by 1 or 2 nearest neighbour interchange operations respectively in $\w{5}$. 
More precisely, we consider trees corresponding to Newick strings
\begin{equation}\label{eq:two-trees-example}
\begin{array}{rclc}
 F_1&:&\big((1:0.1,2:0.1):\ell^6,3:0.1,(4:0.1,5:0.1):\ell^7\big),\\[0.3em]
 F_2&:&\big((2:0.1,3:0.1):\ell^6,1:0.1,(4:0.1,5:0.1):\ell^7\big)&\mbox{ and }\\[0.3em]
 F_3&:&\big((2:0.1,3:0.1):\ell^6,5:0.1,(1:0.1,4:0.1):\ell^7\big).\\ 
 \end{array}
\end{equation}
In both figures, the approximate geodesics constructed using the recursive algorithm are not symmetric under interchange of the end points, and differ from the paths obtained using the symmetric algorithm. 
The lengths of the paths can be computed by using large $k$ and summing $\covmetric$ between successive points in the output. 
In all the examples we explored, the symmetrized algorithm produced shorter paths. 
The approximate geodesics shown in the figures are significantly different from BHV geodesics, which consist of straight (or once broken) lines between the end points in both plots. 

\begin{figure}
\begin{minipage}{\linewidth}
\centering
\begin{minipage}{0.48\linewidth}
\includegraphics[clip, trim= 2cm 7cm 3cm 7cm, width=\linewidth]{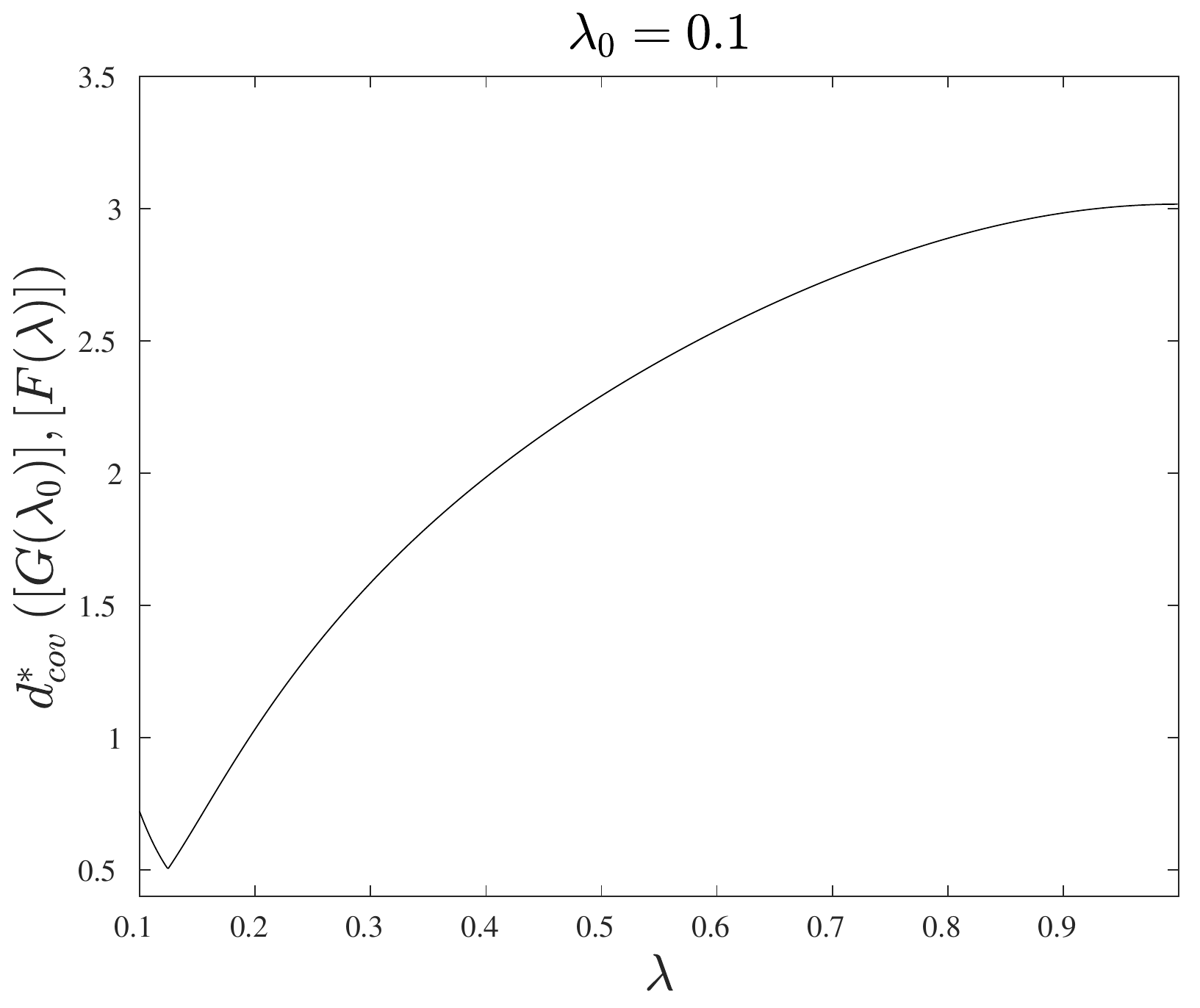}
\end{minipage}
\begin{minipage}{0.48\linewidth}
\includegraphics[clip, trim= 2cm 7cm 3cm 7cm, width=\linewidth]{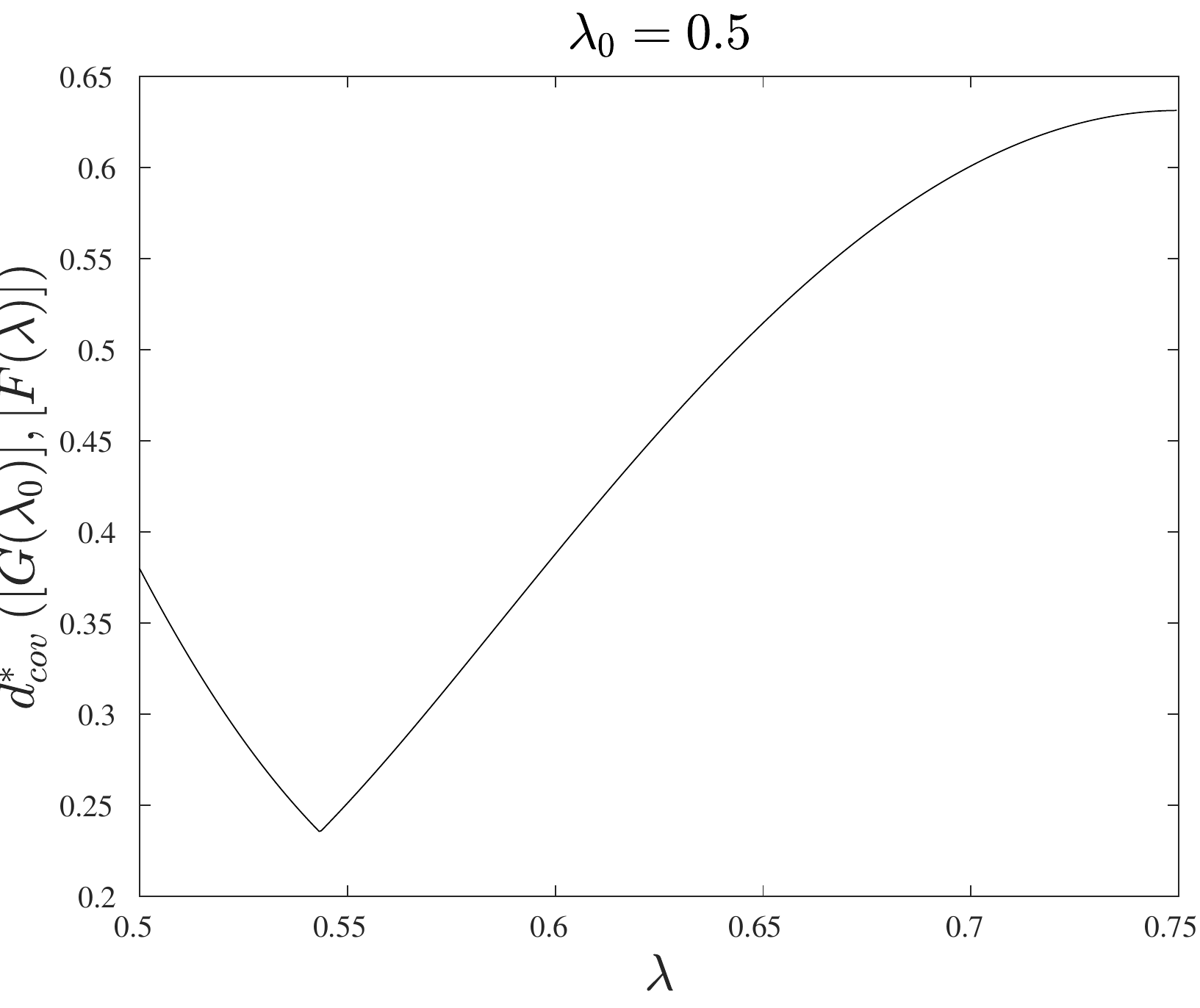}
\end{minipage}
\end{minipage}
\begin{minipage}{\linewidth}
\centering
\begin{minipage}{0.48\linewidth}
\includegraphics[clip, trim= 2cm 7cm 3cm 7cm, width=\linewidth]{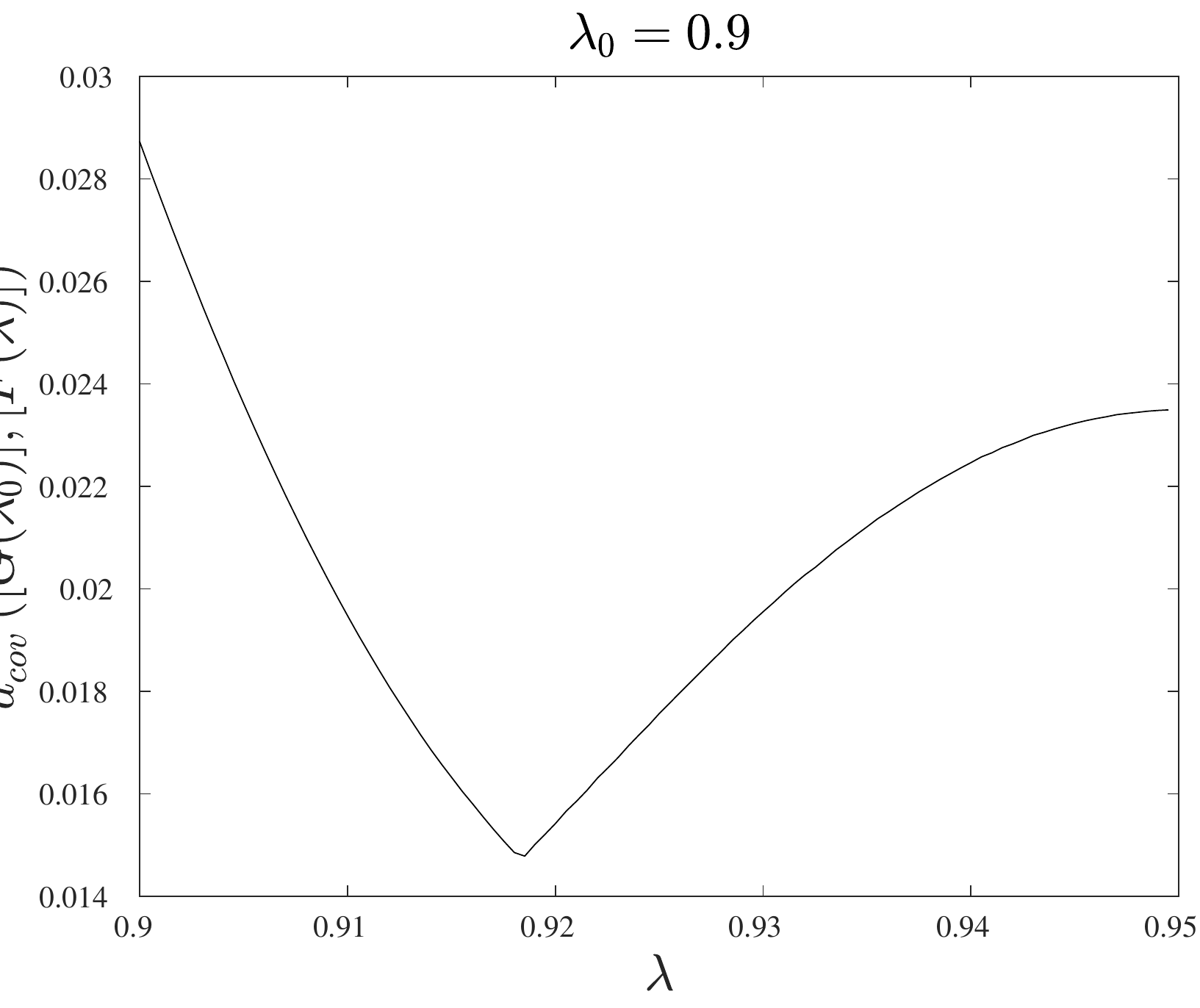}
\end{minipage}
\begin{minipage}{0.48\linewidth}
\includegraphics[clip, trim= 2cm 7cm 3cm 7cm, width=\linewidth]{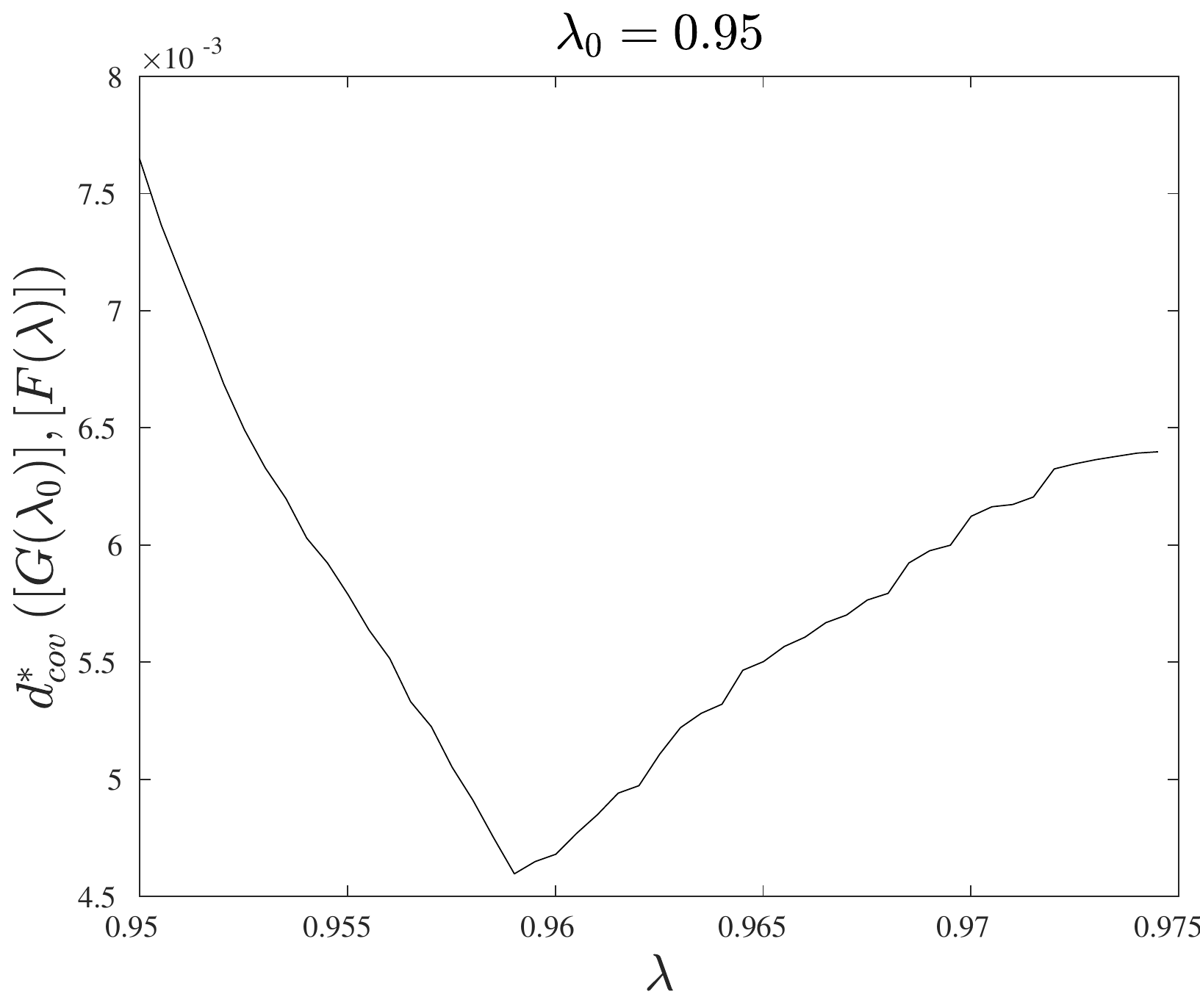}
\end{minipage}
\end{minipage}
\begin{caption}{
\label{fig:approx_distances_to_star_trees}
Approximated distances $\lambda \mapsto \inducedcovmetric\big([G(\lambda_0)], [F(\lambda)]\big)$ for different values of $\lambda_0$. 
}
\end{caption}
\end{figure}

We apply the symmetrized algorithm to investigate the distance from trees in the interior of a maximal orthant to the star stratum on the boundary of that orthant. 
If $[G_1],[G_2],\dots,[G_k]\in\w{N}$ is an approximated geodesic between $[G_1]$ and $[G_k]$ computed by the symmetrized algorithm, the intrinsic distance between $[G_1]$ and $[G_k]$ can be approximated by $\inducedcovmetric\big([G_1],[G_k]\big) \approx \sum_{i=1}^{k-1} \covmetric\big([G_i],[G_{i+1}]\big)$. 
Consider the following setup. 
For $\lambda_0\in (0,1]$, let $G_1 = G(\lambda_0)$ be the forest corresponding to the Newick string $((1: \lambda_0, 2: \lambda_0):\lambda_0, (3:\lambda_0, 4:\lambda_0))$ in $\vec{\lambda}$-parametrization. 
This is a fully resolved 4-taxon tree on which each edge has weight $\lambda_0$. 
By symmetry, the edges on the tree in the star stratum closest to $G(\lambda_0)$ must all have equal weight $\lambda\in (0,1]$. Thus, let $G_k=F(\lambda)$ be the star tree corresponding to the Newick string $(1: \lambda, 2: \lambda, 3: \lambda, 4: \lambda)$, again in $\vec{\lambda}$-parametrization. 
Figure \ref{fig:approx_distances_to_star_trees} shows for each $\lambda_0\in\{0.1, 0.5, 0.9, 0.95\}$ the approximated values of $\inducedcovmetric\big([G(\lambda_0)], [F(\lambda)]\big)$ against $\lambda$. 
Obviously, $F(\lambda)$ is closest for $\lambda$ slightly larger than $\lambda_0$ with distance decreasing as $\lambda_0 \to 1$. 
This suggests that the star stratum is closer to the tree $G(\lambda_0)$ than the forest consisting of 4 isolated points (obtained from $F(\lambda)$ as $\lambda \rightarrow 1$), even for $\lambda_0$ values close to 1. Note, though, that the forest is a boundary point of the star stratum.   
For any $G(\lambda_0)$ the distance to $F(\lambda)$ tends to infinity as $\lambda \to 0$. Indeed, $S_{F(0)}\notin\spd{4}$ is not of full rank.

\section{Discussion}\label{sec:discussion}

In order to do statistical inference on data sets of phylogenetic trees one needs a structure rich enough to enable the use of geometric statistical methods. 
Recent research has produced geometries such as the BHV and tropical tree spaces and statistical methods adapted to these geometries. 
Based on a more principled set of underlying assumptions by regarding phylogenetic trees as probability models for genetic sequence data, we have developed a canonical and biologically motivated geometry on tree space by applying tools from information geometry, giving the \wald.
In particular, unlike previous related work \citep{Garba2018} in which various extrinsic metrics were considered, in this paper we have focused on developing intrinsic metrics and their associated geodesics to explore and to enable accessing the geometry, for this is a key ingredient for statistical inference on non-Euclidean spaces.

There are two main difficulties with achieving our aim. 
First, the discrete-valued Markov process on trees with genetic alphabet $\Omega$ characterizes trees as probability models with sample space $\Omega^N$ where $N$ is the number of phylogenetic taxa. 
Therefore, calculations of distances and construction of geodesics involve summations over $|\Omega|^N$ terms, which becomes infeasible for large $N$. 
In order to establish computational tractability, we generalized the discrete-valued probability model to a continuous-valued Markov process in a canonical way and applied the information geometry again.
 
Secondly, information geometry is formulated for parametrized probability models that are a manifold, whereas tree space is a union of manifolds having different dimensions due to the orthants (representing forests with different number of edges) being glued together in a certain way. 
One has to be careful to compare the structure to the one defined in \cite{Moulton2004}, for example, as we are not including forests with coincident leaves and furthermore we consider a different topology induced by the Fisher information metric. 
We tackled this issue of not having a single connected parametrized manifold by using the continuous-valued Markov models to embed \wald in the ambient space of symmetric positive definite matrices, which has an analytically tractable geometry and thus allows for approximation of geodesics in the embedded space $\w{N}$.
Our computational results show that the geometry obtained is significantly different from the BHV and tropical geometries, partly due to the inclusion of trees with infinitely long edges in \wald. 

Several questions about the geometry of the \wald $\w{N}$ remain. 
While we have shown that trees with infinitely long edges are a finite distance away from other trees (Theorems~\ref{thm:edgeprod-metric-space} and~\ref{thm:edgeprod-metric-space-cont}), computational results suggest that parts of this subspace are repulsive and are avoided by geodesics (see Figure~\ref{fig:cobweblambda} and Figure~\ref{fig:approx_distances_to_star_trees}). 
An explanation for this behaviour might be obtained via calculations or results about curvature for such points of $\w{N}$, but further investigation is required. 
Secondly, Theorems~\ref{thm:edgeprod-metric-space} and~\ref{thm:edgeprod-metric-space-cont} establish $\w{N}$ as length spaces for the two induced intrinsic metrics we study. 
It is desirable to strengthen these results and prove that the distance between every pair of points in the space is realized by at least one path, so that our \wald is a geodesic metric space as opposed to a length space. 
It appears that such a proof requires thorough analysis of the condition on edge weights which excludes trees with coincident leaves. 
Furthermore, the methods and results presented in Section~\ref{sec:projgeod} represent a first step towards the development of more sophisticated and efficient algorithms for the construction of information geodesics in \wald via the embedding in the space of covariance matrices. 
A more thorough evaluation of the computational cost as $N$ increases could be carried out, and a more rigorous treatment might establish convergence properties for the symmetrized algorithm. 
Alternatively, existing algorithms taken from computational Riemannian geometry could be adapted to work in \wald (see~\citet{schmidt2006} for example) and might offer better performance.  

The underlying motivation for this work has been to obtain a novel geometric framework for the space of phylogenetic trees which has more principled biological justification than existing geometries, and which can be used to develop statistical methods for analysing data sets of trees. 
Ultimately, realizing this aim is still some way off. 
For example, given a sample of points $\{x_1,\ldots,x_n\}\subseteq X$ in a metric space $(X,d)$, the Fr\'echet mean $\bar{x}\in X$ is a point which minimises the sum of squared distances to the data:
\begin{equation*}
\bar{x} = \argmin_{x \in X}\sum_{i=1}^n d(x,x_i)^2.
\end{equation*}
In general, the Fr\'echet mean does not always exist, or it can fail to be unique, but in globally non-positively curved spaces such as $(\spd{N},\covmetric)$ and $(\unrooted{N},\bhvmetric)$ there exists a unique Fr\'echet mean \citep{BH1999}. 
Development of methods for calculation of a Fr\'echet mean using an intrinsic information metric in \wald seems very challenging, and the curvature calculations in Section~\ref{sec:continfogeom} have implications for the existence and uniqueness of Fr\'echet means. 
On the other hand, given any sample of trees in $\w{N}$, there is a unique inrinsic Fr\'echet mean in $\spd{N}$ and an algorithm for computing the mean is given by \citet{Lenglet2006}. 
Our projection algorithm could be used to project this to an extrinsic mean back into $\w{N}$. 
Properties of the projected Fr\'echet mean tree could be investigated. 

In comparison to the BHV and tropical metrics, the intrinsic information metrics have the advantage of interpretability in terms of genetic substitutions and the distributions of characters represented by two trees. 
This suggests the information metrics might be better suited for statistical tasks such as hypothesis testing. 
In the BHV and tropical geometries, contraction and expansion of edges offer the means of moving between different topologies. 
In the \wald, additional topological transformations are possible via expanding edges to infinite length, and these correspond to tree bisection and reconnection (TBR) operations. 
Many applications in phylogenetics require searches over the space of phylogenetic trees, and movement along information geodesics in the \wald might have advantages over existing methods. 

\section*{Acknowledgement} The second and the last author express their thanks to the Oberwolfach 1804 meeting ``Statistics for Data with Geometric Structure'' in which \wald was first discussed. The last two authors gratefully acknowledge support from DFG GRK 2088. The last author was supported by the Niedersachsen Vorab of the Volkswagen Foundation.

\section*{Appendix A: Calculation of $p_T(s)$ and its derivatives.}

The probability $p_T(s)$ of any binary character $s$ on a tree $T\in \preunrooted{N}$ can be computed efficiently via the following algorithm \citep{Semple2003}, often called the \emph{Felsenstein pruning algorithm}. 
First, an arbitrary internal vertex $v_0\in T$ is chosen and used to root the tree. 
The two-state symmetric model is a reversible Markov process, and so the choice of the root does not affect the value of $p_T(s)$.
The root determines ancestral relations on the tree, and we let $T_v$ denote the subtree of $T$ descended from vertex $v$. 
We let $L_v$ denote the leaves of $T_v$, and given a binary characer $s$, let $s_v$ denote the restriction of $s$ to $L_v$. 
Finally, we let $p_{T_v}(s_v|\,\omega)$ be the probability of $s_v$ on $T_v$ given the letter $\omega\in\{0,1\}$ at $v$:
\begin{equation*}
p_{T_v}(s_v|\,\omega) = \pr{\bigcap_{u\in L_v} X(u)=s(u) \,\bigg|\, X(v)=\omega},
\end{equation*}
since $s_{v}(u)=s(u)$ for all $u\in L_v$. 
The theorem of total probability gives
\begin{equation}\label{equ:likelihood}
p_T(s) = \frac{1}{2} \sum_{\omega\in\{0,1\}} p_{T_{v_0}}\!(s_{v_0}|\,\omega).
\end{equation}
For an interior vertex $v$, if we let $v_i$, $i=1,\ldots,m$ be the vertices immediately descended from $v$ via edges of length $\ell^i$, then the transition probabilities in Equation~$\eqref{equ:twostatetrans}$ give
\begin{equation}
    \label{equ:felsenstein}
    p_{T_{v}}\!(s_{v}|\,\omega) = \prod_{i=1}^{m}\; \frac{1}{2} 
    \bigg( 
    \big(1+e^{-\ell^i}\big)\,p_{T_{v_i}}\!(s_{v_i}|\,\omega) \,+\,  \big(1-e^{-\ell^i}\big)\,p_{T_{v_i}}\!(s_{v_i}|\,\bar{\omega}) 
    \bigg)
\end{equation}
where $\bar{\omega} = 1-\omega$. For a leaf $u$, we have $p_{T_u}(s_u|\,\omega) = p_{T_u}(s(u)|\,\omega)=1$ if $s(u)=\omega$ and zero otherwise. 
This and Equation~(\ref{equ:felsenstein}) can be applied recursively to compute the terms $p_{T_v}(s_v|\omega)$ for each vertex $v\in T$, starting at the leaves and working up the tree to the root $v_0$. 
Finally, $p_T(s)$ can be computed using Equation~(\ref{equ:likelihood}), and it follows from the recursion that $p_T(s)$ is a multivariate polynomial with arguments of the form $1+e^{-\ell^k}$ and $1-e^{-\ell^k}$, where $k$ ranges over the edges of $T$. 
The coefficients of the polynomial depend on the topology of $T$. 

Equations~$\eqref{equ:likelihood}$ and~$\eqref{equ:felsenstein}$ can also be used to differentiate $p_T(s)$ with respect to any edge length parameter.
These derivatives are required in Section~\ref{sec:twostateinfogeom}. 
Suppose $e$ is an edge of $T$ and we wish to compute the derivative $\partial p_T(s)/\partial\ell^e$. 
Since we are free to choose $v_0$, the calculation is simplified if we let $v_0$ be an internal vertex at one end of edge $e$. 
We can order the vertices $v_i$ attached to $v_0$ so that $e=(v_0,v_1)$. 
Equation~$\eqref{equ:felsenstein}$ then gives
\begin{align*}
\frac{\partial p_{T_{v_0}}\!(s_{v_0}|\,\omega)}{\partial\ell^e} 
&= \frac{1}{2}\, e^{-\ell^e} 
\Big(
p_{T_{v_1}}\!\big(s_{v_1}|\,\bar{\omega}\big) - p_{T_{v_1}}\!\big(s_{v_1}|\,{\omega}\big)
\Big)\\ 
&\times \prod_{i=2}^{\deg{v_0}}\frac{1}{2}
\bigg( 
\big(1+e^{-\ell^i}\big)\,p_{T_{v_i}}\!(s_{v_i}|\,\omega) + \big(1-e^{-\ell^i}\big)\,p_{T_{v_i}}\!(s_{v_i}|\,\bar{\omega}) 
\bigg),
\end{align*}
where the $p_{T_{v_i}}\!$ terms can be calculated recursively. 
Second derivatives of the mass function can be calculated analytically in a similar way.

\section*{Appendix B: Proof of Lemma~\ref{lem:edgeprod}}

First, suppose $F_1\sim F_2$. 
The BHV boundary rule does not affect the distribution on characters induced by a tree, because the same distribution is obtained whether an edge of length zero is present in a tree or not. 
Similarly, if $F_1,F_2$ are equal modulo an application of the boundary rule at infinity, then $p_{F_1}\!(s)=p_{F_2}\!(s)$ since an edge with weight $1$ results in independence between the letters at leaves at either side of the edge under the transition probabilities in Equation~$\eqref{equ:twostatetrans}$. 
Specifically, if $v_0,v_1$ are vertices at the ends of an edge $e$ with $\lambda^e=1$, then 
\begin{equation*}
X(v_1) \,\,\big|\, X(v_0)=\omega ~\ \sim\ Bern(1/2)
\end{equation*} 
where $\omega\in\{0,1\}$, so the conditional distribution of $X(v_1)$ is the same as its marginal. 
The map $[F]\mapsto p_F$ from elements of $\w{N}$ to distributions on characters is therefore well-defined. 
In fact, the work of \citet{Allman2008} shows the map is injective, and this establishes the lemma.

\section*{Appendix C: Proof of Lemma~\ref{lem:KLDfromIP}}

The Riemannian metric in Equation~$\eqref{equ:TwoStateFisherInfo}$ can be expanded 
as
\begin{align*}
\delta\ell^\ia g_{\ia\ib}(\vec{\ell})\delta\ell^\ib 
&= \sum_s p_{\vec{\ell}}(s)
\bigg( \delta\ell^\ia \frac{\partial}{\partial_{\ell^\ia}}\log p_{\vec{\ell}}(s) \bigg) \bigg( \delta\ell^\ib \frac{\partial}{\partial_{\ell^\ib}}\log p_{\vec{\ell}}(s) \bigg)\\
&=\sum_s p_{\vec{\ell}}(s)
\bigg( \delta\ell^\ia\frac{1}{p_{\vec{\ell}}(s)}\frac{\partial p_{\vec{\ell}}(s)}{\partial_{\ell^\ia}} \bigg) 
\bigg( \delta\ell^\ib\frac{1}{p_{\vec{\ell}}(s)}\frac{\partial p_{\vec{\ell}}(s)}{\partial_{\ell^\ib}} \bigg) \\
&=\sum_s \frac{1}{p_{\vec{\ell}}(s)}
\bigg( \delta\ell^\ia\frac{\partial p_{\vec{\ell}}(s)}{\partial_{\ell^\ia}} \bigg) 
\bigg( \delta\ell^\ib\frac{\partial p_{\vec{\ell}}(s)}{\partial_{\ell^\ib}} \bigg).
\end{align*}
The Taylor expansion of $p_{\vec{\ell}}(s)$ is 
\begin{align*}
p_{\vec{\ell}+\delta\vec{\ell}}(s)-p_{\vec{\ell}}(s) = \sum_{\ia}\delta\ell^\ia\frac{\partial p_{\vec{\ell}}(s)}{\partial_{\ell^\ia}} +O\big(|\delta\vec{\ell}|^2\big).
\end{align*}
Substituting this into the expression for the Riemannian metric gives
\begin{align*}
\sum_{\ia,\ib}\delta\ell^\ia g_{\ia\ib}(\vec{\ell})\delta\ell^\ib 
&= \sum_s 
\frac{\big( p_{\vec{\ell}+\delta\vec{\ell}}(s) - p_{\vec{\ell}}(s) +O\big(|\delta\vec{\ell}|^2\big)\big)^2}{p_{\vec{\ell}}(s)}\\ 
&= \sum_s
\frac{\big(p_{\vec{\ell}+\delta\vec{\ell}}(s) - p_{\vec{\ell}}(s)\big)^2}{p_{\vec{\ell}}(s)}+O\big(|\delta\vec{\ell}|^3\big)
\end{align*}
since $p_{\vec{\ell}+\delta\vec{\ell}}(s)-p_{\vec{\ell}}(s)$ is $O\big(|\delta\vec{\ell}|\big)$. 
On the other hand, a Taylor expansion of $f$ around $1$ gives
\begin{align*}
D_f\big(p_{\vec{\ell}+\delta\vec{\ell}};p_{\vec{\ell}}\big) 
&= 
\sum_s p_{\vec{\ell}}(s)\,f\!\left( \frac{ p_{\vec{\ell}+\delta\vec{\ell}}(s) }{ p_{\vec{\ell}}(s) } \right)\\
&\;\begin{aligned}
= \sum_s p_{\vec{\ell}}(s)\, \bigg(f(1) &+ f'(1)\, \frac{p_{\vec{\ell}+\delta\vec{\ell}}(s)-p_{\vec{\ell}}(s)}{p_{\vec{\ell}}(s)}\\
 &+\frac{1}{2}f''(1)\bigg( \frac{ p_{\vec{\ell}+\delta\vec{\ell}}(s)-p_{\vec{\ell}}(s)}{ p_{\vec{\ell}}(s) }\bigg)^2\\
 &+O\bigg( \bigg| \frac{p_{\vec{\ell}+\delta\vec{\ell}}(s)-p_{\vec{\ell}}(s)}{p_{\vec{\ell}}(s)} \bigg|^3 \bigg)\bigg)
\end{aligned}
\\[0.8em]
&\;\begin{aligned}
= f(1) &+ f'(1)\sum_s p_{\vec{\ell}+\delta\vec{\ell}}(s) - f'(1) \sum_s  p_{\vec{\ell}}(s)\\
&+ \frac{1}{2}f''(1)\sum_s \frac{\big( p_{\vec{\ell}+\delta\vec{\ell}}(s)-p_{\vec{\ell}}(s)\big)^2}{ p_{\vec{\ell}}(s) }   +O\left( |\delta\vec{\ell}|^3 \right).
\end{aligned}
\end{align*}
The first three terms vanish since $f(1)=0$ and since $\sum p_{\vec{\ell}}(s)=1$ for all $\vec{\ell}$. 
It follows that 
\begin{equation*}
\sum_{\ia,\ib}\delta\ell^\ia g_{\ia\ib}(\vec{\ell})\delta\ell^\ib = \frac{2}{f''(1)}D_f\big(p_{\vec{\ell}+\delta\vec{\ell}};p_{\vec{\ell}}\big) +O\big(|\delta\vec{\ell}|^3\big)
\end{equation*} 
and the lemma is established.

\section*{Appendix D: Proof of Theorem~\ref{thm:edgeprod-metric-space}}

We need to show that the induced intrinsic metric $d^*\big([F_1],[F_2]\big)$ is finite for any $[F_1],[F_2]\in\w{N}$. 
Start by choosing the representative $F_1\in \prew{N}$ for the equivalence class $[F_1]$ to be a connected tree with edge weights $\vec{\lambda}_1$, some elements of which might have value $1$. 
The tree $F_1$ can be continuously deformed within the orthant corresponding to its topology, to the star tree $F_\ast$ on which all pendant edges have weight $\lambda=1/2$, by changing $\vec{\lambda}$ along the obvious linear path. 
If the path has finite length, then the first part of the theorem has been established, since any $[F_1]$ and $[F_2]$ can be joined to $[F_\ast]$ in this way. 
As shown in Remark \ref{rmk:polynomial}, each $p_{\vec{\lambda}}(s)$ for $\vec{\lambda}$ along the path from $F_1$ to $F_*$ is a polynomial in   $\vec{\lambda}$.
It follows that if $\vec{\lambda}$ and $\vec{\lambda}+\delta\vec{\lambda}$ represent the edge weights at two nearby points on the path then
\begin{equation*}
p_{\vec{\lambda}+\delta\vec{\lambda}}(s)-p_{\vec{\lambda}}(s) =  \pi_\ia(s,\vec{\lambda})\,\delta\lambda^\ia+O\big(|\delta\vec{\lambda}|^2\big)
\end{equation*}
where for each character $s$ and $\ia=1,\ldots,2N-3$, $\pi_\ia(s,\vec{\lambda})$ is a polynomial in $\vec{\lambda}$. 
Then 
\begin{equation*}
\Big(p_{\vec{\lambda}+\delta\vec{\lambda}}(s)-p_{\vec{\lambda}}(s)\Big)^2 =  \Big(\pi_\ia(s,\vec{\lambda})\,\delta\lambda^\ia\Big)^2+O\big(|\delta\vec{\lambda}|^3\big). 
\end{equation*}
The path distance between any pair of leaves is continuous along this path, and is strictly positive since the pendant edge lengths are non-zero at all points on the path, apart from potentially at $F_1$. 
Pendant edge lengths can be zero on $F_1$, but by the definition of $\prew{N}$, the path distance between leaves is non-zero. 
It follows that $p_{\vec{\lambda}}(s)$ is also bound away from zero. 
Thus there is a constant $C(s)$ such that 
\begin{align*}
\frac{\big(p_{\vec{\lambda}+\delta\vec{\lambda}}(s)-p_{\vec{\lambda}}(s)\big)^2}{p_{\vec{\lambda}}(s)} 
&\leq C(s)\Big( \pi_\ia(s,\vec{\lambda})\,\delta\lambda^i \Big)^2\\
&\leq C(s)\bigg( \sum_\ia \pi_\ia(s,\vec{\lambda})^2 \bigg)\big\| \delta\vec{\lambda}\big\|^2
\end{align*}
where the second line comes from the Cauchy-Schwarz inequality and the norm is the Euclidean norm. 
Since the $\pi_\ia(s,\vec{\lambda})$ are polynomials in $\vec{\lambda}$ and the elements of $\vec{\lambda}$ lie between 0 and 1, the $\pi_\ia(s,\vec{\lambda})$ are bounded from above and we obtain 
\begin{equation*}
\frac{\big(p_{\vec{\lambda}+\delta\vec{\lambda}}(s)-p_{\vec{\lambda}}(s)\big)^2}{p_{\vec{\lambda}}(s)} \leq
B(s) \big\| \delta\vec{\lambda} \big\|^2
\end{equation*}
for some constant $B(s)$. 
Now suppose that $D_f=d_f^2$ is a $f$-divergence, where $d_f$ is a metric.
Applying Equation~$\eqref{equ:infinitesimal-fdiv}$ from Lemma~\ref{lem:KLDfromIP} with the  $\vec{\lambda}$-parametrization gives
\begin{align}\label{eqn:metric-Riemann}
d_f^2\big(p_{\vec{\lambda}+\delta\vec{\lambda}}, p_{\vec{\lambda}}\big) &= \frac{1}{2}f''(1)
\sum_s \frac{\big( p_{\vec{\lambda}+\delta\vec{\lambda}}(s)-p_{\vec{\lambda}}(s) \big)^2}{p_{\vec{\lambda}}(s)} + O\big(|\delta\vec{\lambda}|^3\big)\\ \nonumber
& \leq \frac{1}{2}f''(1)\sum_s B(s) \big\|\delta\vec{\lambda}\big\|^2 +O(|\delta\vec{\lambda}|^3)\\ \nonumber
& \leq K \big\|\delta\vec{\lambda}\big\|^2
\end{align}
for some constant $K$. 
Thus the infinitesimal path length in $\w{N}$ as measured by the metric $d$ is bounded by some multiple of the Euclidean path length on $\vec{\lambda}$. 
The length of the linear path from $F_1$ to $F_\ast$ measured with $d$ is therefore finite, since the Euclidean length of this path is finite, and hence $d^*\big([F_1],[F_2]\big)$ is finite. 

For the second part of the theorem, suppose that $D_{f_0}=d_{f_0}^2$ is a $f_0$-divergence, where $d_{f_0}$ is a metric. Further, suppose that $F_1$ and $F_2$ are given by $\vec{\lambda}_1$ and $\vec{\lambda}_2$, respectively. It suffices to consider $\vec{\lambda}_1, \vec{\lambda}_2$ from the same topology and sufficiently close such that the image $t\mapsto \vec{\lambda}(t)$, $\vec{\lambda}(0) = \vec{\lambda}_1$ and $\vec{\lambda}(1) = \vec{\lambda}_2$ of the geodesic from $p_{\vec{\lambda}_1}$ to $p_{\vec{\lambda}_2}$ in the metric induced by the Riemannian metric $g_{ij}$ lies fully in a convex $\vec{\lambda}$ coordinate patch and has finite Euclidean length there, say $L$. Hence, for every $n\in \mathbb{N}$, there are $\delta \vec{\lambda}_{(j)}$ with $\|\delta \vec{\lambda}_{(j)}\| \leq L/n$, $j\in \{1,\ldots,n\}$ such that $\vec{\lambda}_2 =\vec{\lambda}_1 + \sum_{j=1}^{n} \delta\vec{\lambda}_{(j)}$.  Then the second assertion of the theorem follows from (\ref{eqn:metric-Riemann}), setting $c= f''(1)/f''_0(1)$, as $n\to \infty$, because
\begin{align*}
\Big|d_f\big(p_{\vec{\lambda}_1}, p_{\vec{\lambda}_2}\big)^2 - c\cdot d_{f_0}\big(p_{\vec{\lambda}_1}, p_{\vec{\lambda}_2}\big)^2\Big| 
&=  \bigg|\sum_{j=1}^{n-1} O\big(|\delta\vec{\lambda}_{(j)}|^3\big)\bigg| ~=~O\left(\frac{L}{n^2}\right)\,.
\end{align*}
The second equality sign holds because the constants in the individual summands $O\big(|\delta\vec{\lambda}_{(j)}|^3\big)$ ($1\leq j \leq n$) can be bounded by the supremum of absolute values of the gradient of $p_{\vec{\lambda}}$ with respect to $\vec{\lambda}$, in the coordinate patch, 
as can be seen from the last lines of the proof of Lemma \ref{lem:KLDfromIP}.  

%

The third part of the theorem, which states that minimal length paths satisfy the geodesic equation locally, is part of the standard theory for Riemannian geometry on manifolds, e.g. \cite[Section 4]{Lee1997}.

\section*{Appendix E: Proof of Theorem~\ref{thm:posdef}}

The theorem is trivial for  $N = 2$, so suppose $F\in \prew{N}$ with $N\geq 3$ and that the first assertion holds for all $G\in \prew{N-1}$. 
The matrix $S_F$ is not changed by inserting edges $e$ with $\lambda^e = 1$ to connect trees in $F$, or by adding edges with $\lambda^e=0$, so without loss of generality we may assume $F$ is a fully resolved tree. 
We may also assume there is a  cherry between leaves $N-1$ and $N$ since each bifurcating tree with $N\geq 3$ has a cherry and permuting the labels of $F$ results in a tree with covariance $P^T S_F P$ (with permutation matrix $P$), where positive definiteness is preserved.

Let $e_{N-1}$ and $e_{N}$ be the edges incident to leaves $N-1$ and $N$, respectively. Since $F$ is bifurcating, there is exactly one edge, say $e_0$, incident to $e_{N-1}$ and $e_{N}$. 
Let $S_F = (s_{uv})_{u,v=1}^N$. 
The tree $G\in\prew{N-1}$ obtained by deleting $e_{N}$ and leaf $N$ and merging $e_0$ and $e_{N-1}$ to $\tilde{e}$ with weight $\lambda_{\tilde{e}} = 1 - (1 - \lambda_{e_0})(1 - \lambda_{e_{N-1}})$ has covariance $S_{G} = (s_{uv})_{u,v=1}^{N-1}$, which is by induction positive definite.
Using this and Sylvester's criterion (that a matrix is positive if and only if all principal minors have positive determinant) it suffices to show $\det(S_F) > 0$. 
We have for all $1\leq u \leq N-1$ that $s_{uN} = s_{Nu} = (1 - \lambda_{e_{N}})c_u$, where 
\begin{equation*}
c_u = \begin{cases}
1 - \lambda_{e_{N-1}} & \textrm{when $u=N-1$, and}\\ 
\prod_{e\neq e_N, e_{N-1}} \big(1 - \lambda^e\big)^{\sigma_{Nu}^e} & \textrm{for $1\leq u \leq N-2$}
\end{cases}
\end{equation*}
and $\sigma_{Nu}^e$ is defined by Equation~$\eqref{equ:splitmatrix}$. 
Note that for $u,v\leq N-1$, $c_u$ and $s_{uv}$ do not involve $\lambda_{e_{N}}$, and that $s_{NN} = 1$.
If $\perm_N$ denotes the set of permutations of $\{1,\ldots,N\}$, then the Leibniz formula for determinants gives
\begin{align*}
&\det(S_F) = \bigg(
 \sum_{\substack{\tau\in\perm_N\\\tau(N) = N}} \mathrm{sgn}(\tau)\prod_{u = 1}^N s_{u\tau(u)}\bigg) + 
 \bigg(\sum_{\substack{\tau\in\perm_N\\\tau(N) \neq N}} \mathrm{sgn}(\tau) \prod_{u = 1}^N s_{u\tau(u)}\bigg)\\
 &= \det\!\Big((s_{uv})_{u,v=1}^{N-1}\Big) + 
  (1 - \lambda_{e_N})^2\bigg(\sum_{\substack{\tau\in\perm_N\\\tau(N) \neq N}}\!\! \mathrm{sgn}(\tau) c_{\tau(N)} c_{\tau^{-1}(N)}\prod_{\substack{u = 1\\u\neq \tau^{-1}(N)}}^{N-1}\!\!\! s_{u\tau(u)}\bigg),
\end{align*}
so $\det(S_F)$ is linear in $x \coloneqq (1 - \lambda_{e_N})^2$. 
By symmetry of the cherry, $\det(S_F)$ is also linear in $y\coloneqq (1 - \lambda_{e_{N-1}})^2$ as well. 
We write $g(x, y) = \det(S_F)$. 
For $x = 0$, we have $s_{Nu} = 0$ for $u<N$ and $s_{NN} = 1$, so $g(0, y) = \det(S_F) = \det(S_{G}) > 0$ for all $y\in[0,1]$, and similarly $g(x, 0) > 0$ for all $x\in[0,1]$. Furthermore, $g(1, 1) = 0$, since in that case the last two rows of $S_F$ coincide. Since $g$ is linear in $x$ and in $y$, respectively, we have $g(x,y) > 0$ for all $(x, y) \in [0,1]^2\setminus\{(1,1)\}$, so that $\det(S_F) > 0$ for all $(\lambda_{e_{N-1}},\lambda_{e_N})\in[0,1]^2\setminus\{(0,0)\}$. 
If $\lambda_{e_{N-1}} = \lambda_{e_{N}} = 0$, we would have $d_{N(N-1)} = 0$, but this is not allowed by the definition of $\prew{N}$. 

We also need to show that the map $[F]\mapsto S_F$ is injective on $\w{N}$ where $[F]$ denotes the equivalence class of $F\in \prew{N}$.
This is trivial, however, since whenever $F_1,F_2\in \prew{N}$ are in different equivalence classes, the matrix of distances between the leaves is different.

\section*{Appendix F: Proof of Theorem~\ref{thm:edgeprod-metric-space-cont}}

The proof is similar to that for Theorem~\ref{thm:edgeprod-metric-space}, and so we give a brief sketch. 
We consider the same path between the trees $[F_1],[F_\ast]\in\w{N}$, and show that each element of $g_{ij}(\vec{\lambda})$ is bound from above along the path.  
Working in the $\vec{\lambda}$-parametrization of an orthant, Equation~$\eqref{equ:covHadamard}$ becomes 
\begin{equation*}
S_{\vec{\lambda}} = \bigg(\prod_e\big( 1-\lambda^e \big)^{\sigma^e_{uv}}\bigg)_{u,v=1}^N.
\end{equation*}
Each element of the matrix is therefore a polynomial in the elements of $\vec{\lambda}$, and their derivatives with respect to $\vec{\lambda}$ are also polynomials. 
Recalling that the tangent space of $ \w{N}$ at $S_{\vec{\lambda}}$ in a maximal orthant is spanned by $\partial_\ia S_{\vec{\lambda}}$, where $i \in \{1,\ldots,2N-3\}$ ranges over the edges in that maximal orthant,  Equation~$\eqref{equ:FIMtrace}$ becomes
\begin{equation*}
g_{\ia\ib}(\vec{\lambda}) = \frac{1}{2}\,\trace\Big( S_{\vec{\lambda}}^{-1} \big(\partial_\ia S_{\vec{\lambda}}\big)S_{\vec{\lambda}}^{-1} \big(\partial_\ib S_{\vec{\lambda}}\big) \Big)\,,
\end{equation*}
 $i,j \in \{1,\ldots,2N-3\}$.
Applying the Cauchy-Schwartz inequality $|\trace(A^TB)|^2\leq \trace(A^TA)\trace(B^TB)$ gives 
\begin{align*}
\big|g_{\ia\ib}(\vec{\lambda})\big|^2 &\leq \frac{1}{4}\,
\trace\Big( \big(\partial_\ia S_{\vec{\lambda}}\big)^2S_{\vec{\lambda}}^{-2} \Big)\,
\trace\Big( \big(\partial_\ib S_{\vec{\lambda}}\big)^2S_{\vec{\lambda}}^{-2} \Big)\\
&\leq \frac{1}{4}\,\trace\Big( S_{\vec{\lambda}}^{-4} \Big)\,
\trace\Big( \big( \partial_\ia S_{\vec{\lambda}}\big)^4\Big)^\frac{1}{2}
\trace\Big( \big( \partial_\ib S_{\vec{\lambda}}\big)^4\Big)^\frac{1}{2}.
\end{align*}
The first term in this product is bounded on a geodesic path from $F_1$  to $F_\ast$, since $S_{\vec{\lambda}}$ is positive definite and its eigenvalues are bound away from zero. 
The other two terms are also bounded from above, because the derivatives of $S_{\vec{\lambda}}$ are polynomials in $\vec{\lambda}$, 
Thus $|g_{\ia\ib}(\vec{\lambda})|\leq C$ for some constant $C$ at all points along that path, and the same argument as for Theorem~\ref{thm:edgeprod-metric-space} shows that $\inducedcovmetric\big([F_1],[F_\ast]\big)$ is finite. 

%


\bibliographystyle{apa}
\bibliography{information-geometry-revision}   

\begin{thebibliography}{}

\bibitem[\protect\astroncite{Adams and Castoe}{2019}]{adams2019}
Adams, R.~H. and Castoe, T.~A. (2019).
\newblock Probabilistic species tree distances: implementing the multispecies
  coalescent to compare species trees within the same model-based framework
  used to estimate them.
\newblock {\em Syst. Bio.}

\bibitem[\protect\astroncite{Allen and Steel}{2001}]{allen2001}
Allen, B.~L. and Steel, M. (2001).
\newblock Subtree transfer operations and their induced metrics on evolutionary
  trees.
\newblock {\em Ann. Comb.}, 5(1):1--15.

\bibitem[\protect\astroncite{Allman et~al.}{2008}]{Allman2008}
Allman, E.~S., An\'e, C., and Rhodes, J.~A. (2008).
\newblock Identifiability of a {M}arkovian model of molecular evolution with
  gamma-distributed rates.
\newblock {\em Adv. Appl. Probab.}, 40(1):229--249.

\bibitem[\protect\astroncite{Ballmann et~al.}{1985}]{ballmann1985}
Ballmann, W., Gromov, M., and Schroeder, V. (1985).
\newblock {\em Manifolds of nonpositive curvature}, volume~61 of {\em Progress
  in mathematics}.
\newblock Birkh{\"a}user.

\bibitem[\protect\astroncite{{Ba\v{c}\'ak}}{2014}]{Bacak2014}
{Ba\v{c}\'ak}, M. (2014).
\newblock Computing medians and means in {H}adamard spaces.
\newblock {\em SIAM J. Optim.}, 24(3):1542--1566.

\bibitem[\protect\astroncite{Billera et~al.}{2001}]{BHV2001}
Billera, L., Holmes, S., and Vogtman, K. (2001).
\newblock Geometry of the space of phylogenetic trees.
\newblock {\em Adv. Appl. Math.}, 27:733--767.

\bibitem[\protect\astroncite{Bridson and Haefliger}{2011}]{BH1999}
Bridson, M.~R. and Haefliger, A. (2011).
\newblock {\em Metric Spaces of Non-Positive Curvature}.
\newblock Springer, Berlin.

\bibitem[\protect\astroncite{Bryant et~al.}{2005}]{bryant05}
Bryant, D., Galtier, N., and Poursat, M.-A. (2005).
\newblock Likelihood calculation in molecular phylogenetics.
\newblock In Gascuel, O., editor, {\em Mathematics of Evolution and Phylogeny},
  pages 33--62. Oxford University Press.

\bibitem[\protect\astroncite{Dryden et~al.}{2009}]{dryden2009}
Dryden, I.~L., Koloydenko, A., Zhou, D., et~al. (2009).
\newblock Non-{E}uclidean statistics for covariance matrices, with applications
  to diffusion tensor imaging.
\newblock {\em Ann. Appl. Stat.}, 3(3):1102--1123.

\bibitem[\protect\astroncite{Engstr{\"o}m et~al.}{2013}]{engstrom2013toric}
Engstr{\"o}m, A., Hersh, P., and Sturmfels, B. (2013).
\newblock Toric cubes.
\newblock {\em Rendiconti del Circolo Matematico di Palermo}, 62(1):67--78.

\bibitem[\protect\astroncite{Feragen et~al.}{2013}]{Feragen2013}
Feragen, A., Owen, M., Petersen, J., Wille, M., Thomsen, L., Dirksen, A., and
  de~Bruijne~M. (2013).
\newblock Tree-space statistics and approximations for large-scale analysis of
  anatomical trees.
\newblock In {\em 23rd biennial {I}nternational {C}onference on {I}nformation
  {P}rocessing in {M}edical {I}maging (IPMI)}.

\bibitem[\protect\astroncite{Garba}{2019}]{phdthesis}
Garba, M.~K. (2019).
\newblock {\em Information geometry for phylogenetic trees}.
\newblock PhD thesis, School of Mathematics, Statistics and Physics, Newcastle
  University.

\bibitem[\protect\astroncite{Garba et~al.}{2018}]{Garba2018}
Garba, M.~K., Nye, T. M.~W., and Boys, R.~J. (2018).
\newblock Probabilistic distances between trees.
\newblock {\em Syst. Bio.}, 67(2):320--327.

\bibitem[\protect\astroncite{Gill et~al.}{2008}]{Gill2008}
Gill, J., Linusson, S., Moulton, V., and Steel, M. (2008).
\newblock A regular decomposition of the edge-product space of phylogenetic
  trees.
\newblock {\em Adv. Appl. Math.}, 41(2):158--176.

\bibitem[\protect\astroncite{Hansen and Martins}{1996}]{hansen1996}
Hansen, T.~F. and Martins, E.~P. (1996).
\newblock Translating between microevolutionary process and macroevolutionary
  patterns: the correlation structure of interspecific data.
\newblock {\em Evolution}, 50(4):1404--1417.

\bibitem[\protect\astroncite{Huelsenbeck and
  Ronquist}{2001}]{huelsenbeck2001mrbayes}
Huelsenbeck, J.~P. and Ronquist, F. (2001).
\newblock {MrBayes}: {B}ayesian inference of phylogenetic trees.
\newblock {\em Bioinformatics}, 17(8):754--755.

\bibitem[\protect\astroncite{Kim}{2000}]{Kim2000}
Kim, J. (2000).
\newblock Slicing hyperdimensional oranges: The geometry of phylogenetic
  estimation.
\newblock {\em Mol. Phylogenet. Evol.}, 17(1):58--75.

\bibitem[\protect\astroncite{Lee}{1997}]{Lee1997}
Lee, J.~M. (1997).
\newblock {\em Riemannian manifolds: an introduction to curvature}, volume 176.
\newblock Springer.

\bibitem[\protect\astroncite{Lenglet et~al.}{2006}]{Lenglet2006}
Lenglet, C., Rousson, M., Deriche, R., and Faugeras, O. (2006).
\newblock Statistics on the manifold of multivariate normal distributions:
  Theory and application to diffusion tensor \text{MRI} processing.
\newblock {\em J. Math. Imaging Vis.}, 25(3):423--444.

\bibitem[\protect\astroncite{Lin et~al.}{2018}]{Lin2019}
Lin, B., Monod, A., and Yoshida, R. (2018).
\newblock Tropical foundations for probability and statistics on phylogenetic
  tree space.
\newblock {\em arXiv preprint arXiv:1805.12400}.

\bibitem[\protect\astroncite{Lin and Yoshida}{2018}]{Lin2018}
Lin, B. and Yoshida, R. (2018).
\newblock Tropical {F}ermat--{W}eber points.
\newblock {\em SIAM J. Discrete Math.}, 32(2):1229--1245.

\bibitem[\protect\astroncite{Miller et~al.}{2015}]{MOP2015}
Miller, E., Owen, M., and Provan, J.~S. (2015).
\newblock Polyhedral computational geometry for averaging metric phylogenetic
  trees.
\newblock {\em Adv. Appl. Math.}, 68:51--91.

\bibitem[\protect\astroncite{Moakher}{2005}]{moakher_differential_2005}
Moakher, M. (2005).
\newblock A {Differential} {Geometric} {Approach} to the {Geometric} {Mean} of
  {Symmetric} {Positive}-{Definite} {Matrices}.
\newblock {\em SIAM J. Matrix Analysis Applications}, 26:735--747.

\bibitem[\protect\astroncite{Moulton and Steel}{2004}]{Moulton2004}
Moulton, V. and Steel, M. (2004).
\newblock Peeling phylogenetic oranges.
\newblock {\em Adv. Appl. Math.}, 33(4):710--727.

\bibitem[\protect\astroncite{Nye}{2014}]{Nye2014}
Nye, T. (2014).
\newblock An algorithm for constructing principal geodesics in phylogenetic
  treespace.
\newblock {\em IEEE ACM T. Comput. Bi.}, 11(2):304--315.

\bibitem[\protect\astroncite{Nye}{2011}]{Nye2011}
Nye, T. M.~W. (2011).
\newblock Principal components analysis in the space of phylogenetic trees.
\newblock {\em Ann. Statist.}, 39(5):2716--2739.

\bibitem[\protect\astroncite{Nye et~al.}{2017}]{Nye2017}
Nye, T. M.~W., Tang, X., Weyenberg, G., and Yoshida, R. (2017).
\newblock Principal component analysis and the locus of the {F}r\'echet mean in
  the space of phylogenetic trees.
\newblock {\em Biometrika}, 104(4):901--922.

\bibitem[\protect\astroncite{Owen and Provan}{2011}]{Owen2011}
Owen, M. and Provan, J.~S. (2011).
\newblock A fast algorithm for computing geodesic distances in tree space.
\newblock {\em IEEE ACM T. Comput. Bi.}, 8(1):2--13.

\bibitem[\protect\astroncite{Rogers}{1997}]{rogers1997consistency}
Rogers, J.~S. (1997).
\newblock On the consistency of maximum likelihood estimation of phylogenetic
  trees from nucleotide sequences.
\newblock {\em Systematic biology}, 46(2):354--357.

\bibitem[\protect\astroncite{Sason and Verdu}{2016}]{sason2016f}
Sason, I. and Verdu, S. (2016).
\newblock $ f $-divergence inequalities.
\newblock {\em IEEE Transactions on Information Theory}, 62(11):5973--6006.

\bibitem[\protect\astroncite{Schmidt et~al.}{2006}]{schmidt2006}
Schmidt, F.~R., Clausen, M., and Cremers, D. (2006).
\newblock Shape matching by variational computation of geodesics on a manifold.
\newblock In {\em Joint Pattern Recognition Symposium}, pages 142--151.
  Springer.

\bibitem[\protect\astroncite{Semple and Steel}{2003}]{Semple2003}
Semple, C. and Steel, M. (2003).
\newblock {\em Phylogenetics, Oxford Lecture Series in Mathematics and its
  Applications, 24}.
\newblock Oxford University Press.

\bibitem[\protect\astroncite{Skovgaard}{1984}]{Skovgaard1984}
Skovgaard, L.~T. (1984).
\newblock A {R}iemannian geometry of the multivariate normal model.
\newblock {\em Scand. J. Stat.}, 11(4):211--223.

\bibitem[\protect\astroncite{Speyer and Sturmfels}{2004}]{Speyer2004}
Speyer, D. and Sturmfels, B. (2004).
\newblock The tropical {G}rassmannian.
\newblock {\em Adv. Geom.}, 4(3):389--411.

\bibitem[\protect\astroncite{Steel and Penny}{1993}]{steel1993}
Steel, M.~A. and Penny, D. (1993).
\newblock Distributions of tree comparison metrics -- some new results.
\newblock {\em Syst. Biol.}, 42(2):126--141.

\bibitem[\protect\astroncite{Willis}{2019}]{Willis2019}
Willis, A. (2019).
\newblock Confidence sets for phylogenetic trees.
\newblock {\em J. Am. Stat. Assoc.}, 114(525):235--244.

\bibitem[\protect\astroncite{Yang}{2006}]{Yang2006}
Yang, Z. (2006).
\newblock {\em Computational molecular evolution}.
\newblock Oxford University Press.

\bibitem[\protect\astroncite{Yoshida et~al.}{2019}]{Yoshida2019}
Yoshida, R., Zhang, L., and Zhang, X. (2019).
\newblock Tropical principal component analysis and its application to
  phylogenetics.
\newblock {\em B. Math. Biol.}, 81(2):568--597.

\bibitem[\protect\astroncite{Zwiernik and Smith}{2012}]{Zwiernik2012}
Zwiernik, P. and Smith, J.~Q. (2012).
\newblock Tree cumulants and the geometry of binary tree models.
\newblock {\em Bernoulli}, 18(1):290--321.

\end{thebibliography}

\end{document}